\documentclass[10pt]{article}
\usepackage{graphicx}
\usepackage{amsmath}
\usepackage{amssymb}
\usepackage{amscd}
\usepackage{hhline}
\usepackage{epic}
\usepackage{mathrsfs}

\newtheorem{theorem}{Theorem}
\newtheorem{theorema}{Theorem}
\newtheorem{theoremb}{Theorem}
\newtheorem{theoremc}{Theorem}
\newtheorem{theoremd}{Theorem}
\newtheorem{theoreme}{Theorem}
\newtheorem{dfn}[theoremb]{Definition}
\newtheorem{rk}[theoremc]{Remark}
\newtheorem{cor}[theoremd]{Corollary}

\newtheorem{examp}[theoreme]{Example}
\newtheorem{prop}[theorem]{Proposition}
\newenvironment{proof}[1][Proof]{\textbf{#1.} }{\qed \vspace{5pt}}
\newcommand\bib[1]{\bibitem[#1]{#1}}

\newcommand\abz{\hspace{13.5pt}}
\newcommand\qed{\phantom{\underline{y}}\hfill\hfill$\square$}
\newcommand{\comm}[1]{}

\newcommand\1{{\bf 1}}
\renewcommand\a{\alpha}
\renewcommand\b{\beta}

\newcommand\C{{\mathbb C}}

\renewcommand\d{\delta}
\newcommand\D{{\mathcal D}}

\newcommand\E{{\mathcal E}}

\newcommand\F{{\frak F}}

\newcommand\g{\gamma}

\newcommand\hps{\hskip-16pt . \hskip2pt}
\newcommand\hpss{\hskip-13.5pt . \hskip2pt}

\newcommand\I{{\mathcal I}}
\renewcommand\k{\varkappa}
\renewcommand\l{\lambda}

\newcommand\La{\Lambda}

\renewcommand\O{\Omega}
\newcommand\oo{\omega}
\newcommand\op[1]{\mathop{\rm #1}\nolimits}
\newcommand\ot{\otimes}
\newcommand\p{\partial}

\newcommand\po{$\!\!\!{\text{\bf.}}$ }

\newcommand\R{{\mathbb R}}
\renewcommand\t{\tau}

\newcommand\tg{\tilde g}
\newcommand\ti{\tilde}

\newcommand\vp{\varphi}
\newcommand\we{\wedge}
\newcommand\x{\xi}
\newcommand\z{\sigma}
\newcommand\Z{{\mathbb Z}}

\DeclareFontFamily{U}{wncyr}{}
\DeclareFontShape{U}{wncyr}{m}{it}{%
  <5><6><7><8><9>gen*wncyi%
  <10><10.95><12><14.4><17.28><20.74><24.88>wncyi10}{}
\DeclareSymbolFont{MathRussLetters}{U}{wncyr}{m}{it}
\DeclareMathSymbol{\re}{\mathalpha}{MathRussLetters}{3}

\makeatletter
\renewcommand{\@oddhead}{\hfil Compatibility of PDEs via multi-brackets\hfil}
\makeatother \makeatletter
\renewcommand{\@evenhead}{\hfil Boris Kruglikov, Valentin Lychagin\hfil}
\makeatother

\newcommand\Cc{\let\mathcal\mathscr\mathcal C}

\begin{document}

 \title{Compatibility, multi-brackets and \\
 integrability of systems of PDEs}
 \author{Boris Kruglikov, Valentin Lychagin}
 \date{}
 \maketitle

 \vspace{-14.5pt}
 \begin{abstract}
We establish an efficient compatibility criterion for a system of
generalized complete intersection type in terms of certain
multi-brackets of differential operators. These multi-brackets
generalize the higher Jacobi-Mayer brackets, important in the
study of evolutionary equations and the integrability problem. We
also calculate Spencer $\d$-cohomology of generalized complete
intersections and evaluate the formal functional dimension of the
solutions space. The results are applied to establish new
integration methods and solve several differential-geometric
problems.
 \end{abstract}

\section*{Introduction and main results}

 \abz
In this paper we introduce multi-brackets of non-linear vector
differential operators. In the case of bi-brackets they coincide
with the well-known Jacobi bracket, which is a generalization of the
classical Lagrange-Jacobi bracket important in the theory of 1st
order differential equations. These latter brackets become the usual
commutators in the case of linear differential operators and are
widely used in mathematical physics and non-linear analysis. We
apply multi-bracket of differential operators to establish a
criterion of formal integrability of systems of PDEs.

\subsection{\hpss Multi-brackets of non-linear differential operators}\label{S01}

 \abz
Our multi-bracket $\{F_1,\dots,F_{m+1}\}$ is defined for
differential operators on sections of the trivial $m$-dimensional
bundle over a manifold $M$ (notice that trivialization assumption is
usually implicit for commutators or bi-brackets) and its value
is an operator of the same kind. When $F_i$ are linear vector
differential operators $\nabla_i:m\cdot C^\infty_\text{loc}(M)\to
C^\infty_\text{loc}(M)$, represented as rows
$(\nabla_i^1,\dots,\nabla_i^m)$ of scalar linear differential
operators, the multi-bracket has the form:
 $$
\{\nabla_1,\dots,\nabla_{m+1}\}=\sum_{k=1}^{m+1}(-1)^k
\op{Ndet}[\nabla_i^j]_{i\ne k}^{1\le j\le m}\cdot\nabla_k,
 $$
where $\op{Ndet}$ is a version of non-commutative determinant and
$\cdot$ is the product (one can perceive any determinant and
product for a while, but we will discuss various versions in the
sequel). For non-linear operators the bracket is obtained via
linearization.

If $\pi$ is a trivial vector bundle and
$\ell(F)=(\ell_1(F),\dots,\ell_m(F))$ is a linearization of the
operator $F$, then:
 $$
\hspace{-0.05in}\{F_{1},\dots,F_{m+1}\}=\dfrac{1}{m!}\hspace{-0.25in}%
\sum_{\vphantom{\frac{2^2_2}2}\ \alpha\in S_m,\beta\in
S_{m+1}}\hspace{-0.25in}\left( -1\right) ^{\alpha}\left( -1\right)
^{\beta}\ell_{\alpha(1)}(F_{\beta \left(  1\right)
})\circ\ldots\circ\ell_{\alpha(m)}(F_{\beta\left( m\right)
})\left(F_{\beta(m+1)}\right).
 $$

We will need restriction of this bracket to the system $\E$ of
PDEs $F_i=0$. Let $\op{ord}(F_i)=l(i)$ be the orders of the
considered operators. Denote by $\op{diff}(\pi,\1)$ the algebra of
all (non-linear) scalar differential operators on the sections of
$\pi$.

Define $\mathcal{J}_s(F_1,\dots,F_k)=\langle\D_\t
F_i\,\bigl|\,l(i)+|\t|\le s\rangle$ to be the subalgebra of the
differential ideal generated by $F_1,\dots,F_k$ in
$\op{diff}(\pi,\1)$, where $\D_\t$ is the total derivative
operator with respect to a multi-index $\t$ (formula in
\S\ref{S11}). We will explore the following reduced multi-bracket
(further discussion in \S\ref{S31d}):
 $$
\{F_1,\dots,F_{m+1}\}\mod\mathcal{J}_{l(1)+\dots+l(m+1)-1}(F_1,\dots,F_{m+1}).
 $$

The above equivalence class can be given by other multi-brackets,
more convenient for calculations. In the canonical coordinates
$(x^i,p^j_\z)$ on the jet-space $J^\infty(\pi)$ of a rank $m$
(vector) bundle $\pi$ over the base $M$ with $\op{dim}M=n$ the
reduced multi-bracket has the following representative:
 $$
[F_1,\dots,F_{m+1}]'=\dfrac1{m!}\sum_{\z\in
S_{m+1}}\!\!\!\op{sgn}(\z)\sum_{\t_i}
\left|\dfrac{\p(F_{\z(1)},\dots,F_{\z(m)})}
{\p(\,p^1_{\t_1}\,,\,\dots\,,\,p^m_{\t_m}\,)}\right|
\D_{\t_1+\dots+\t_m}F_{\z(m+1)},
 $$
where
 $$
\left|\dfrac{\p(f_1,\dots,f_m)}{\p(u_1,\dots,u_m)}\right|=\det\left\|\dfrac{\p
f_i}{\p u_j}\right\|_{m\times m}.
 $$
is the Jacobian. In other words we have (see \S\ref{S31b} for
details):
 $$
[F_1,\dots,F_{m+1}]'\equiv
\{F_1,\dots,F_{m+1}\}\mod\mathcal{J}_{l(1)+\dots+l(m+1)-1}(F_1,\dots,F_{m+1}).
 $$

When $m=1$ we obtain the Mayer bracket $[F,G]$ of scalar
differential operators. This bracket coincides with the classical
Lie-Mayer bracket for the first order equations and is closely
related to the Jacobi bracket on higher jets. We applied it in the
previous works (\cite{KL$_1$,KL$_2$,KL$_3$}) to establish a
compatibility criterion for overdetermined systems of scalar
equations of a certain type.

Namely, we considered a system of complete intersection type, i.e.
given by $r\le n=\op{dim}M$ equations which have transversal
characteristic varieties on regular strata. In other words, the
ideal generated by the symbols of the operators is an algebraic
complete intersection.

\subsection{\hpss Main results}\label{S02}

 \abz
In this paper we extend the compatibility result to the case of
systems of PDEs with arbitrary number of unknown functions. To
characterize the systems for which the criterion is sufficient
(necessity holds always) we introduce a new class of systems
generalizing the complete intersection for the scalar case.

The conditions informally say the system is not too overdetermined
(we will also discuss the opposite case) and is of general kind
(transversality condition).

 \begin{dfn}\po\label{dfn1}
We say a system $\E\subset J^k(\pi)$ of $r$ differential equations
on $m$ unknowns is of {\em generalized complete intersection
type\/} if
 \begin{enumerate}
\item $m<r<n+m$;
\item The characteristic variety has $\dim_\C\op{Char}^\C_{x_k}(\E)=n+m-r-2$
at each point $x_k\in\E$ $($we assume $\dim\emptyset=-1)$;
\item The characteristic sheaf $\mathcal{K}$ over
$\op{Char}^\C_{x_k}(\E)\subset P^\C T^*$ has fibers of dimension 1
everywhere (see \S\ref{S11}-\S\ref{S13} for details of the
involved objects).
 \end{enumerate}
The case $r=m$ corresponds to determined systems, where the
compatibility conditions are void, but all the statements
hold for this case as well.
 \end{dfn}

The class of systems, introduced above, is included into the systems
of Cohen-Macaulay type, introduced in \cite{KL$_2$}, see also the
discussion of complete intersection for PDEs there. Note that the
number $r$ of equations, called codimension of the system $\E$, is
defined invariantly and is calculated via the Spencer
$\d$-cohomology by the formula $r=\op{codim}(\E)=\dim H^{*,1}(\E)$,
see \cite{KL$_2$}.

Define the reduced multi-bracket due to the system $\langle
F_1,\dots,F_r\rangle$ by the formula
 $$
[F_{i_1},\dots,F_{i_{m+1}}]_\E=\{F_{i_1},\dots,F_{i_{m+1}}\}
\mod\mathcal{J}_{l(i_1)+\dots+l(i_{m+1})-1}(F_1,\dots,F_r).
 $$

 \begin{theorema}\po\label{th A}
Consider a system of PDEs
 $$
\E=\Bigl\{F_i\Bigl(x^1,\dots,x^n,u^1,\dots,u^m,
\frac{\p^{|\z|}u^j}{\p x^\z}\Bigr)=0\,|\,1\le i\le r\Bigr\}, \quad
\op{ord}(F_i)=l(i).
 $$
1. If the system $\E$ is formally integrable, then the
multi-bracket vani\-shes due to the system, i.e. for every
collection $1\le i_1<\dots<i_{m+1}\le r$
 $$
[F_{i_1},\dots,F_{i_{m+1}}]_\E=0.
 $$
2. Let $\E$ be a system of PDEs of generalized complete
intersection type. Then the system $\E$ is formally integrable if
and only if the multi-bracket vani\-shes due to the system:
 $$
[F_{i_1},\dots,F_{i_{m+1}}]_\E=0.
 $$
 \end{theorema}

In particular, we deduce the following compatibility criterion for
scalar PDEs:

 \begin{cor}\po
Let $\E$ be a system of complete intersection type, i.e. given by
$r\le n$ differential equations on one unknown function
$F_1[u]=0,\dots,F_r[u]=0$ of orders $l_1,\dots,l_n$. Then the
system $\E$ is formally integrable iff the Mayer bracket vanishes
due to the system, i.e. the Jacobi bracket satisfies:
 $$
 \hspace{40pt}
\{F_i,F_j\}=0\mod\mathcal{J}_{l_i+l_j-1}(F_1,\dots,F_r),
\quad\forall\ 1\le i<j\le r.
 \hspace{40pt}\square
 $$
 \end{cor}

Theorem \ref{th A} was announced in \cite{KL$_4$}. The corollary was
established in \cite{KL$_3$} and its particular cases for $n=2$ and
$r=2$ appeared in \cite{KL$_1$} and \cite{KL$_2$} respectively. We
notice however that the technique used in these papers was different
and we required an additional assumption that at least one of the
equations has no multiple components in the characteristic variety.
We remarked then that this condition is superfluous, but proved the
claim only for the second order equations. Now we can remove this
technical assumption completely.

Recall (\cite{S,GS,Go}) that the obstructions to integrability
belong to the second Spencer $\d$-cohomology group $H^{*,2}(\E)$
(we recall the definition in \S\ref{S21}). Thus it is important to
calculate this bi-graded cohomology
$H^{*,j}(\E)=\oplus_iH^{i,j}(\E)$.

 \begin{theorema}\po\label{th B}
Let $\E$ be a system of differential equations defined by a set of
$r=\op{codim}(\E)$ differential operators
$\Delta=(\Delta_1,\dots,\Delta_r):C^\infty(\pi)\to C^\infty(\nu)$
(can be of different orders). If $\E$ is a generalized complete
intersection, then the only non-zero Spencer $\d$-cohomology are
given by the formula:
 $$
H^{*,j}(\E)=\left\{\begin{array}{ll}
 \pi &\text{\rm for }\ j=0,\\
 \nu &\text{\rm for }\ j=1,\\
 S^{j-2}\pi^*\ot\La^{m+j-1}\nu &\text{\rm for }\ 2\le j\le r+1-m\, (\le n).
\end{array}\right.
 $$
 \end{theorema}
In the above formula we describe $H^{*,j}(\E)$ as a usual
(non-graded) vector space. See \S\ref{S41} for more information
about grading.

For the case of scalar systems $m=\dim\pi=1$ we have:
$H^{*,j}(\E)\simeq\La^j\nu$, $0\le j\le r$. This corresponds to the
following well-known algebraic result: Algebra $g^*$ of
$\op{codim}g^*=r$ is a complete intersection iff its Koszul homology
forms the exterior algebra $H_i(g^*)=\La^i H_1(g^*)$, $0\le i\le r$
(\cite{BH}).

The precise obstructions to formal integrability $W_i(\E)$ are
certain curvature-type invariants called Weyl tensors \cite{L$_1$}.
In \cite{KL$_1$,KL$_2$} we calculated them for $\op{codim}(\E)=2$
complete intersections in terms of Jacobi-Mayer brackets. Now we can
generalize this result in terms of our multi-brackets:

 \begin{cor}\po\label{cor2}
There is a basis $e_1,\dots,e_s$ in $H^{*,2}(\E)$,
$s=\binom{r}{m+1}$, and a bijection $\psi$ between the set of
power $(m+1)$ subsets of $\{1,\dots,r\}$ and the set
$\{1,\dots,s\}$ such that the graded Weyl tensor
$W(\E)=\oplus_iW_i(\E)$ of the system of equations
$\E=\{F_i=0\,|\,1\le i\le r\}$ with $l(i)=\op{ord}F_i$ equals
 $$
 \hspace{50pt}
W(\E)=\sum_{1\le i_1<\dots<i_{m+1}\le
r}[F_{i_1},\dots,F_{i_{m+1}}]_\E\cdot
e_{\psi(i_1,\dots,i_{m+1})}.
 \hspace{48pt}\square
 $$
 \end{cor}

This follows directly from theorems \ref{th A} and \ref{th B}. What
is more interesting is the precise form of the basis. We calculated
it for the case of 2 scalar equations in \cite{KL$_1$,KL$_2$}. The
result immediately generalizes to arbitrary complete intersections.
The case $m>1$ is more involved and we do not discuss it here.

Finally we give a result on the space
$\mathcal{S}_\E=\mathcal{S}ol_\E$ of local/formal solutions of the
system $\E\subset J^k(\pi)$ of generalized complete intersection
type. As before $n$ is dimension of the base $M$ of $\pi$ and $m$
its rank. Let $r$ be the formal codimension (see \S\ref{S21}) of
the system $\E$, the same number as in definition \ref{dfn1}.

In classical textbooks the solutions space is characterized as
follows: a general solution (a generic point of $\mathcal{S}ol_\E$)
depends on $s_p$ functions of $p$ variables, $s_{p-1}$ functions of
$(p-1)$ variables, \dots and $s_0$ constants, where $s_i$ are Cartan
characters (introduced by E.Cartan \cite{C}; we adapt notations from
\cite{BCG$^3$}).

Here $p$ (called genre of $\E$) is the maximal number, such that
$s_p\ne0$: only this character $s_p$ has absolute meaning (citing
\cite{C}). We call the number $p$ {\em formal functional
dimension\/} of the solutions space $\mathcal{S}_\E$ and the
number $d=s_p$ {\em formal functional rank\/}.

The above numbers are well-defined in analytic category, i.e. when
the PDEs and the solutions are considered analytic, see Cartan's
test \cite{BCG$^3$}. Cartan-K\"ahler theorem guarantees
integrability. For smooth equations we need to impose additional
requirements on the system to ensure that the space
$\mathcal{S}ol_\E$ is non-empty and regular (see \cite{Ho,S,M}). In
general we take $p$ and $d$ to be the formal functional dimension
and rank of the space of formal solutions.

In abstract terms the number $d$ equals $P_\E^{(p-1)}(t)$, where
$P_\E(t)$ is the Hilbert polynomial of the symbolic module of $\E$
and $p=\deg P_\E(t)+1$ (see more in \S\ref{S43+e}; for the detailed
discussion of this subject we refer to \cite{KL$_5$}).

 \begin{theorema}\po\label{th C}
Let $\E$ be a formally integrable system of generalized complete
intersection type. Denote its orders by $k_1,\dots,k_r$ and the
corresponding $l$-th symmetric polynomials by
$S_l(k_1,\dots,k_r)=\sum_{i_1<\dots<i_l}k_{i_1}\cdots k_{i_l}$.
Then the space $\mathcal{S}ol_\E$ has formal functional dimension
and rank equal respectively
 $$
p=m+n-r-1,\qquad d=S_{r-m+1}(k_1,\dots,k_r).
 $$
 \end{theorema}

Thus in our case $p$ is dimension of the affine characteristic
variety and when $r=m+n-1$ the space $\mathcal{S}ol_\E$ is a
$d$-dimensional smooth manifold. When $m<r<m+n-1$ and the system
$\E$ is analytical, the space $\mathcal{S}ol_\E$ is
infinite-dimensional and a general analytic solution depends on
precisely $d$ arbitrary functions of $m+n-r-1$ variables.

In smooth category the above formula for the functional rank $d$
is important for formulation of well-posed boundary value
problems. Note also that due to Cauchy-Kovalevskaya theorem the
above theorem holds true in the case $r=m$ of determined system of
PDEs.

\subsection{\hpss Discussion and plan of the paper}\label{S03}

 \abz
The main result (Theorem \ref{th A}) provides an explicit
compatibility criterion. To our knowledge there were only two such
criteria before. One is a particular case of our theorem for the
first order scalar systems of PDEs -- this was one of the
motivations for the appearance of the brackets (see the historical
note in \cite{KL$_1$}) and the base for Lagrange-Charpit method (see
\cite{Gou} and \S\ref{S42}).

Another classical result concerns the system of linear evolution
equations and the compatibility is expressed via commutators,
being thus also a special case of our general result. In fact, all
known integrability methods use these simple compatibility
criteria, see \S\ref{S43}.

All other methods are algorithmic, but non-explicit, and are based
on the Cartan's prolongation-projection idea. We mention two, which
apply in the non-linear situation. One is the Spencer theory
\cite{Go,S} and the Weyl tensors in the 2$^\text{nd}$
$\d$-cohomology groups \cite{L$_1$}. Another uses the differential
Gr\"obner basis and is being implemented into computer algebra
systems now \cite{Ma,Hu}. However neither of them give precise
formulas and from computational point of view our criterion is more
effective \cite{K$_1$}.

Theorem \ref{th B} can be specified to bi-degrees, see
\S\ref{S41}. This is important, since it yields the place, where
the system becomes involutive. In fact, we think that the
generalized complete intersections represent the class of systems,
where the amount of prolongations to achieve involutivity is
maximal. This gives a possibility to reduce the estimate in the
Poincar\'e $\delta$-lemma (see \cite{Sw}, but this estimate is
accepted to be too large).

Theorem \ref{th C} gives an asymptotic estimate for the Hilbert
polynomial of the symbolic module of the system. The dimension
formula is important for Lagrange-Charpit method of establishing
exact solutions of PDEs.

The paper is organized as follows. In Section 1 we collect the
background on the jet-spaces and linear differential operators and
establish a machinery to check the formal integrability. In
Section 2 we review the algebraic machinery and develop the
commutative algebra concepts finishing with a resolvent for
generalized complete intersections. In Section 3 we introduce
multi-brackets and discuss their properties. Non-linear
differential equations are treated geometrically (as in
\S\ref{S11}) and we refer the reader to \cite{KLV} for more
details.

In Section 4 we prove Theorem \ref{th A} for linear systems and
then extend the methods to the non-linear situation. We apply our
results to construct the compatibility complex and non-linear
Spencer cohomology. Theorems \ref{th B} and \ref{th C} are proved
in Section 5, where we also relate our results to classical
integrability methods and multi-Poisson geometry.

In Section 6 we apply the compatibility criterion to solve some
classical problems in differential geometry. We discuss invariant
characterization of Liouville metrics on surfaces and the
generalized Bonnet problem. Previously the compatibility criterion
was applied to the plane web-geometry to solve the Blashke
conjecture and to count Abelian relations \cite{GL$_1$,GL$_2$}. This
illustrates efficiency of our main result.

\section{\hps Jet-spaces and linear differential operators}

 \abz
In this sections we collect the basic knowledge of the geometric
theory of differential equations required for our goals.

\subsection{\hpss Systems of PDEs}\label{S11}

 \abz
Let $M$ be a smooth $n$-dimensional manifold and $\pi:E_\pi\to M$
a (vector) bundle of rank $m$. Two local sections $s_1,s_2\in
C^\infty_\text{loc}(\pi)$ having tangency of order $\ge k$ are
said to have the same $k$-jet at $x\in M$ and the equivalence
class is called the $k$-jet $x_k=[s]^k_x$.

Thus we obtain the jet-bundle $\pi_k:J^k(\pi)\to M$ and there are
natural projections $\pi_{k,l}:J^k(\pi)\to J^l(\pi)$ for $l<k$. We
denote $x_l=\pi_{k,l}(x_k)$. Any smooth section $s\in
C^\infty_\text{loc}(\pi)$ induces the local section
$j_ks:x\mapsto[s]^k_x$ of the bundle $\pi_k$.

A system of PDEs $\E$ is represented as a collection of subsets
$\E_k\subset J^k(\pi)$, $k\ge0$, satisfying certain conditions.
The first one, regularity, is that $\E_k$ with restricted map
$\pi_{k,l}$ is a (fiber) bundle. To formulate the second condition
let us define for a submanifold $\E\subset J^k(\pi)$ its
$i^{\text{th}}$ prolongation by the formula
 $$
\E^{(i)}\!=\!\{x_{k+i}=[s]^{k+i}_x\!\in\!J^{k+i}(\pi):j_ks(M)
\text{ is tangent to }\E\text{ at }x_k\text{ with order}\!\ge i\}.
 $$
Thus we form $\E$ by a collection of some given equations $\E_k$
and the other $\E_l$ are obtained via the prolongation.

So a system of different order PDEs is the following collection of
submanifolds: $\E_i=J^i(\pi)$ until a certain order $l_0$, at
which we add some PDEs and get $\E_{l_0}\subset J^{l_0}(\pi)$,
then $\E_i=\E_{l_0}^{(i-l_0)}$ for $l_0<i<l_1$, whereupon we add
new equations, obtain $\E_{l_1}$, prolong this system until
jet-level $l_2$ etc.

Following Cartan's prolongation-projection scheme we consider
$\pi_{i+s,i}(\E_{i+s})$ and if this is a proper subset of $\E_i$,
the system becomes inconsistent in the sense that we need to add
some equations not specified in the original system.

If we wish to exclude this we obtain: The system $\E$ is said to
be {\em compatible\/} on the level $k$ if
$\emptyset\ne\E_{k+1}\subset\E_k^{(1)}$ and
$\pi_{k+1,k}:\E_{k+1}\to\E_k$ is surjective. In the regular case
the last map is a bundle projection (submersion).

The system $\E$ is said to be integrable to order $k$ if it is
compatible on every level $l\le k$.
System $\E$ is called {\em formally integrable\/} if it is
integrable to order $\infty$ (we usually assume regularity).

We always assume there are no functional equations in $\E$, i.e.\
$\E_0=J^0(\pi)=E_\pi$. The minimal $l$ such that $\E_l\ne
J^l(\pi)$ is called the {\em minimal order\/} $l_0$ of the system.
Every number $l$ with the property $\E_l\ne\E_{l-1}^{(1)}$ is
called an {\em order\/} and codimension of $\E_l$ in
$\E_{l-1}^{(1)}$ is called its {\em multiplicity\/}.

Due to Cartan-Kuranishi theorem on prolongations (in the regular
case) the set of orders $\op{ord}(\E)\subset\mathbb{N}$ is finite,
i.e.\ there exists a {\em maximal order\/} $l_\text{max}$ starting
from which $\E_l^{(1)}=\E_{l+1}$.

Every local coordinate system $(x^i,u^j)$ on the bundle $\pi$
induces coordinates $(x^i,p^j_\z)$ on $J^k(\pi)$ (multiindex
$\z=(i_1,\dots,i_n)$ has length $|\z|=\sum_{s=1}^ni_s\le k$),
where $p^j_\z\bigl([s]^k_x\bigr)=\dfrac{\p^{|\z|}s^j}{\p
x^\z}(x)$. We call them {\em canonical coordinates\/}.

In a sequel we will need the operator of {\em total differential\/}
(also denoted $\hat d$):
 $$
\D:C^\infty(J^k(\pi))\to\O^1(M)\ot_{C^\infty(M)}C^\infty(J^{k+1}(\pi)).
 $$
To define $\D$ we note that every function on $J^k(\pi)$ is a
scalar differential operator $\square:C^\infty(\pi)\to\
C^\infty(M)$ of order $k$. Post-composing it with a vector field
$X\in\D(M)$ we get a differential operator
$\square':C^\infty(\pi)\to\ C^\infty(M)$ of order $k+1$, producing
the needed $1^\text{st}$ order differential operator
$\D_X=i_X\circ\D :C^\infty(J^k(\pi))\to C^\infty(J^{k+1}(\pi))$.

If we write in local coordinates $X=\sum X^i\p_{x^i}$, then
$\D_X=\sum X^i\D_i$ in the corresponding canonical coordinates,
where the operator of {\em total derivative\/} $\D_i=\D_{\p_{x^i}}$
is given by the infinite series (when applied, only finite number of
terms act non-trivially):
 \begin{equation}\label{D}
\D_i=\p_{x^i}+\sum p^j_{\z+1_i}\p_{p^j_\z}.
 \end{equation}
Similarly for $X\in S^l\D(M)$ we get the operator
$\D_X:C^\infty(J^k(\pi))\to C^\infty(J^{k+l}(\pi))$. For instance,
if $\z=(i_1,\dots,i_n)$ is a multiindex, we obtain
$\D_\z=\D_1^{i_1}\cdots\D_n^{i_n}$.

\subsection{\hpss Linear differential operators}\label{S12}

 \abz
Denote by $\1$ the trivial one-dimensional bundle over $M$. Let
$\mathcal{A}_k=\op{Diff}_k(\1,\1)$ be the $C^\infty(M)$-module of
scalar linear differential operators of order $\le k$ and
$\mathcal{A}=\cup_k\mathcal{A}_k$ be the corresponding filtered
algebra,
$\mathcal{A}_k\circ\mathcal{A}_l\subset\mathcal{A}_{k+l}$.

Consider two linear vector bundles $\pi,\nu$. Denote by
$\op{Diff}(\pi,\nu)=\cup_k\op{Diff}_k(\pi,\nu)$ the filtered
module of all differential operators from $C^\infty(\pi)$ to
$C^\infty(\nu)$. We have the natural pairing
 $$
\op{Diff}_k(\rho,\nu)\times\op{Diff}_l(\pi,\rho)\to
\op{Diff}_{k+l}(\pi,\nu)
 $$
given by the composition of differential operators.

In particular, $\op{Diff}(\pi,\1)$ is a filtered left
$\mathcal{A}$-module, $\op{Diff}(\1,\pi)$ is a filtered right
$\mathcal{A}$-module and they have an $\mathcal{A}$-valued
$\mathcal{A}$-linear pairing
 $$
\Delta\in\op{Diff}_l(\pi,\1),\ \nabla\in\op{Diff}_k(\1,\pi)\mapsto
\langle\Delta,\nabla\rangle=\Delta\circ\nabla\in\mathcal{A}_{k+l},
 $$
with
$\langle\a\Delta,\nabla\rangle=\a\langle\Delta,\nabla\rangle$,
$\langle\Delta,\nabla\a\rangle=\langle\Delta,\nabla\rangle\a$ for
$\a\in\mathcal{A}$.

Each differential operator $\Delta:C^\infty(\pi)\to C^\infty(\nu)$
of order $l$ induces an $\mathcal{A}$-homomorphism
$\phi_\Delta:\op{Diff}(\1,\pi)\to\op{Diff}(\1,\nu)$ by the
formula:
 $$
\op{Diff}_k(\1,\pi)\ni\nabla\mapsto\Delta\circ\nabla
\in\op{Diff}_{k+l}(\1,\nu).
 $$
Its $\langle\,,\rangle$-dual is the $\mathcal{A}$-homomorphism
$\phi^\Delta:\op{Diff}(\nu,\1)\to\op{Diff}(\pi,\1)$ given by
 $$
\op{Diff}_k(\nu,\1)\ni\square\mapsto\square\circ\Delta
\in\op{Diff}_{k+l}(\pi,\1).
 $$

By the very definitions of jets with
$\mathcal{J}^k(\pi)=C^\infty(\pi_k)$ we have:
 \begin{equation}\label{krem}
\op{Diff}_k(\pi,\nu)=\op{Hom}_{C^\infty(M)}(\mathcal{J}^k(\pi),C^\infty(\nu)),
 \end{equation}
and differential operators $\Delta$ are in bijective correspondence
with morphisms $\psi^\Delta:J^l(\pi)\to\nu$ via the formula
$\Delta=\psi^\Delta\circ j_l$, where $j_l:C^\infty(\pi)\to
\mathcal{J}^l(\pi)$ is the jet-section operator.

The prolongation $\psi^\Delta_k:J^{k+l}(\pi)\to J^k(\nu)$ of
$\psi^\Delta=\psi^\Delta_0$ is conjugated to the
$\mathcal{A}$-homomorphism
$\phi^\Delta:\op{Diff}_k(\nu,\1)\to\op{Diff}_{k+l}(\pi,\1)$ via
isomorphism (\ref{krem}). This makes a geometric interpretation of
the differential operator $\Delta$ as the bundle morphism.

Similarly, the homomorphism
$\phi_\Delta:\op{Diff}_k(\1,\pi)\to\op{Diff}_{k+l}(\1,\nu)$ is
conjugated to the following morphism:
 $$
 \begin{CD}
\op{Hom}(J^k(\1),\pi)@>{\psi_\Delta^k}>>\op{Hom}(J^{k+l}(\1),\nu)\\
@VV{\text{prolongation}}V  @AA{i^*}A\\
\op{Hom}(J^l(J^k(\1)),J^l(\pi))@>{\psi^\Delta}>>\op{Hom}(J^l(J^k(\1)),\nu),
 \end{CD}
 $$
where $i:J^{k+l}(\1)\to J^l(J^k(\1))$ is the natural embedding.

\subsection{\hpss Systems of differential equations as modules}\label{S13}

 \abz
A system $\E$ of PDEs of order $l$ associated to $\Delta$ is, by
definition, the subbundle $\E_l=\op{Ker}(\psi^\Delta)\subset
J^l(\pi)$. Its prolongation is
$\E_{k+l}=\E_l^{(k)}=\op{Ker}(\psi^\Delta_k)\subset J^{k+l}(\pi)$.

We define the dual $\E^*$ as the cokernel of the operator
$\phi^\Delta$:
 $$
\op{Diff}_k(\nu,\1)\stackrel{\phi^\Delta_k}
\longrightarrow\op{Diff}_{k+l}(\pi,\1)\to\E_{k+l}^*\to0.
 $$
So $\E^*=\{\E_i^*\}$. We have a natural map
$\pi_{i+1,i}^*:\E_i^*\to\E^*_{i+1}$. Then we define the inductive
limit $\E^\Delta=\underrightarrow\lim\E_i^*$. Notice that it is a
filtered left $\mathcal{A}$-module. Thus we can consider a system
as a module over differential operators ($\mathcal{D}$-module).

The dual $\E_\Delta=\op{Ker}(\phi_\Delta)\subset\op{Diff}(\1,\pi)$
is a right $\mathcal{A}$-module and we have a pairing
$\E^\Delta\times\E_\Delta\to\mathcal{A}$. This pairing is however
non-degenerate only for formally-integrable systems. This follows
from the following statement:

 \begin{prop}\po
A system $\E=\op{Ker}(\psi^\Delta)$ is formally integrable iff
$\E_i^*$ are projective $C^\infty(M)$-modules and the maps
$\pi_{i+1,i}^*$ are injective.
 \end{prop}

 \begin{proof}
The projectivity condition is equivalent to regularity (constancy
of rank), while invjectivity of $\pi_{i+1,i}^*$ is equivalent to
surjectivity of $\pi_{i+1,i}$.
 \end{proof}


We can associate to the above modules their symbolic analogs.
Namely, since $ST^*M\ot\pi=\oplus S^iT^*M\ot\pi$ is the graded
module associated to the filtrated $C^\infty(M)$-module
$\op{Diff}(\pi,\1)=\cup\op{Diff}_i(\pi,\1)$, the bundle morphism
$\phi^\Delta$ produces the homomorphisms (symbols)
$\z_k^\Delta:S^{k+l}T^*M\ot\pi\to S^kT^*M\ot\nu$ of our
differential operator $\Delta$. Its dual is the graded degree $l$
morphism
 $$
\z_\Delta:STM\ot\nu^*\to STM\ot\pi^*.
 $$
The value $\z_{\Delta,x}$ of $\z_\Delta$ at $x\in M$ is a
homomorphism of $ST_xM$-modules.

The $ST_xM$-module $\mathcal{M}_\Delta=\op{Coker}(\z_{\Delta,x})$
is called the {\em symbolic module\/} at $x\in M$. Its annihilator
is the called {\em characteristic ideal\/} $I(\Delta)$ and the set
of its zeros is the {\em characteristic variety\/}
$\op{Char}(\Delta)$. We will always consider projectivization of
this conical affine variety.

Moreover in this paper we shall complexify the symbolic module and
work with complex characteristics. In particular, the
characteristic variety becomes $\op{Char}_x^\C(\Delta)\subset P^\C
T_x^*M$.

 \begin{prop}{\bf \cite{Go,S}.}
For $p\in T_x^*M\setminus\{0\}$ let ${\frak m}(p)\subset
S(T_xM)=\oplus_{i\ge0}S^iT_xM$ be the maximal ideal of homogeneous
polynomials vanishing at $p$. Then localization
$(\mathcal{M}_\Delta)_{{\frak m}(p)}\ne0$ iff the covector $p$ is
characteristic.
 \end{prop}

The set of the localizations $(\mathcal{M}_\Delta)_{{\frak
m}(p)}\ne0$ for characteristic covectors $p$ form the {\em
characteristic sheaf\/} $\mathcal{K}$ over the characteristic
variety $\op{Char}_x^\C(\Delta)$.

If we have several differential operators
$\Delta_i\in\op{Diff}(\pi,\nu_i)$ of different orders $l_i$, $1\le
i\le t$, then their sum is no longer a differential operator of
pure order $\Delta:C^\infty(\pi)\to C^\infty(\nu)$,
$\nu=\oplus\nu_i$. Then $\phi^\Delta$ is not an
$\mathcal{A}$-morphism, unless we put certain weights to the
graded components $\nu_i$. Still we have the bundle morphism
$\psi^\Delta$ and the symbol map $\z_\Delta:STM\ot\nu^*\to
STM\ot\pi^*$, which becomes a homomorphism after a suitable
weighting (in \S\ref{S23}). This will be used in the next section
to pursuit the theory into the general setting of various orders
systems.

\section{\hps Algebra of differential equations}

 \abz
Here we review the basics of symbolic theory, establish the
preparatory material and extend it to the general case of
non-linear differential equations.

\subsection{\hpss Spencer cohomology}\label{S21}

 \abz
We consider at first the symbolic theory. Let $T=T_xM$ be the
tangent space to the base and $N=T_{x_0}\pi^{-1}(x)$, $x_0\in
E_\pi$, the tangent space to the fiber of $\pi$. We can identify
$F_k(x_k)=T_{x_k}[\pi^{-1}_{k,k-1}(x_{k-1})]$ with $S^kT^*\ot N$
and let $g_k=g(x_k)=T_{x_k}\E_k\cap F_k(x_k)$ be the {\em
symbol\/} of differential equation $\E$. Clearly $g_0=N$ and
$g_i=S^iT^*\ot N$ for $i<l_1$ -- the minimal order of the system.

The symbol of the de Rham operator is called Spencer $\d$-operator
 $$
\d:S^kT^*\ot N\to S^{k-1}T^*\ot N\ot T^*
 $$
and it maps $g_k$ to $g_{k-1}\ot T^*$. In other words, if
$$g_k^{(1)}=\{p\in S^{k+1}T^*\ot N\,|\,\d p\in g_k\ot T^*\}$$ is the
{\em first prolongation\/}, which in the regular case equals the
symbol of the equation $\E_k^{(1)}$, then $g_k\subset
g_{k-1}^{(1)}$.

 \begin{dfn}\po
A sequence of subspaces $g_k\subset S^kT^*\ot N$, $k\ge0$, is
called a {\em symbolic system\/} if $g_{k+1}\subset g_k^{(1)}$.
 \end{dfn}

Thus symbols of a PDE system form a symbolic system. With every
such a system we associate its Spencer $\d$-complex of order $k$:
 $$
0\to g_k\stackrel{\d}\to g_{k-1}\ot T^*\stackrel{\d}\to
g_{k-2}\ot\La^2T^*\to\dots\stackrel{\d}\to g_{k-n}\ot\La^nT^*\to0.
 $$
The cohomology group at the term $g_i\ot\La^jT^*$ is denoted by
$H^{i,j}(\E,x_k)$, though we usually omit reference to the point
and can also write $H^{i,j}(g)$.

In terms of this cohomology $l$ is an order of the system $\E$ if
$H^{l-1,1}(\E)\ne0$ and multiplicity of this order is equal to
$m(l)=\dim g_{l-1}^{(1)}/g_l=\dim H^{l-1,1}(\E)$.

Now the number $r$ of equations defining $\E$, which we called in
\cite{KL$_2$} {\em formal codimension\/}, is the total number of
orders counted with multiplicities, i.e.
 $$
r=\op{codim}(\E)=\sum_{k=1}^\infty\dim H^{k-1,1}(\E).
 $$
By Hilbert basis theorem this number is finite.

Consider a symbolic system $g=\{g_k\subset S^kT^*\ot
N\,|\,k\ge0\}$ and let $g^*=\oplus g_k^*$ be its graded dual over
$\R$ (or possibly $\C$) . Then $g^*$ is an $ST$-module with the
structure operation given by the formula
 $$
\langle w\cdot\k,p\rangle=\langle \k,\d_wp\rangle,\quad w\in
S^kT,\ \k\in g^*_l,\ p\in g_{k+l},
 $$
where $\d_w=\d_{w_1}\cdots\d_{w_k}$ for $w=w_1\cdots w_k\in S^kT$,
$w_j\in T$, and $\d_{w_j}=i_{w_j}\circ\d:g_t\to g_{t-1}$. This $g^*$
is called the {\em symbolic module\/}. It coincides with the module
$\mathcal{M}_\Delta$ introduced in \S\ref{S13} in the case of linear
equations of the same order.

This module is Noetherian and the Spencer cohomology of $g$
dualizes to the Koszul homology of $g^*$.

\subsection{\hpss Characteristic variety and Fitting ideals}\label{S22}

 \abz
Define the {\em characteristic ideal\/} by
$I(g)=\op{ann}(g^*)\subset ST$ in the ring of polynomials $R=ST$
and the {\em characteristic variety\/} as the set of non-zero
covectors $v\in T^*$ such that for every $k$ there exists a vector
$w\in N\setminus\{0\}$ with $v^k\ot w\in g_k$. This is a punctured
conical affine variety. We projectivize its complexification and
denote it by $\op{Char}^\C(g)\subset P^\C T^*$. When $g=g(x_k)$ is
the symbol of the system at a point $x_k\in\E$, we also denote the
characteristic variety by $\op{Char}^\C_{x_k}(\E)$.

Another definition of characteristic variety is given via the
homogeneous characteristic ideal $I(g)$ graded by the degree:
$I=\oplus I_k$.
 \begin{prop}{\bf \cite{S}.}
$\op{Char}^\C(g)=\{p\in P^\C T^*\,|\,f(p^k)=0\,\forall f\in I_k,
\forall k\}$.
 \end{prop}

Consider the symbolic $R$-module $g^*$. Its {\em dimension\/}
$\dim_Rg^*$ is the minimal number $d$ of homogeneous elements
$f_1,\dots,f_d\in R$ of positive degree such that the quotient
$g^*/(f_1,\dots,f_d)g^*$ is a finite-dimensional vector space.
Thus due to equality $\dim_Rg^*=\dim[R/\op{ann}(g^*)=R/I(g)]$, we
can interpret
 $$
\dim_Rg^*=\dim_\C\op{Char}^\C(g)+1
 $$
as dimension of the affine characteristic variety.

A sequence of elements $f_1,\dots,f_s\in R$ is called regular if
$f_i$ is not a zero divisor in the $R$-module
$g^*/(f_1,\dots,f_{i-1})g^*$. {\em Depth\/} of the module $g^*$
is the maximal number of elements in a regular sequence. The depth
and dimension of a module $g^*$ are related by the following
inequality:
 $$
\op{depth}g^*\le\dim g^*
 $$
(we shall omit sometimes the subscript $R$ in $\dim_R$). Now $g^*$
is called a Cohen-Macaulay module if $\op{depth}g^*=\dim g^*$ (see
\cite{BH} for details). In such a case we call the system $\E$ and
the corresponding symbolic system $g$ Cohen-Macaulay.

For an ideal $I\subset R$ and an $R$-module $G$ the length of
maximal $G$-regular sequence in $I$ is denoted $\op{depth}(I,G)$.
{\em Depth\/} of the ideal $I$ is $\op{depth}I=\op{depth}(I,R)$
(this quantity is also called the {\em grade\/} of the ideal $I$).
In these terms the depth of the module is
 $\op{depth}g^*=\op{depth}(\mathfrak{m},R/\op{Ann}g^*)$,
where $\mathfrak{m}=\oplus_{i>0}S^iT$ is the maximal ideal in $R$
of homogeneous polynomials with positive degree.

We shall also use {\em codimension\/} of the ideal $I$, which is
defined as $\op{codim}I=\min\dim R_\mathfrak{p}$, the lower bound
being taken over all primes $\mathfrak{p}\supset I$ in $R$
\cite{E} (in other sources it is called {\em height\/}
\cite{BH,BV}).

Both the depth and the codimension are geometric quantities, i.e.
they are defined by the conical affine locus of the ideal $I$ over
$\C$: If $\op{Rad}(I)$ is the radical of $I$, then
$\op{depth}\op{Rad}(I)=\op{depth}I$ and
$\op{codim}\op{Rad}(I)=\op{codim}I$. Moreover, since our ring $R$
is polynomial, for any ideal $I\subset R$ we have the equality
 $$
\op{depth}I=\op{codim}I.
 $$

For a homomorphism of free $R$-modules $\vp:U\to V$ with
$\op{rank}U=r$ and $\op{rank}V=m$ denote by $J_j(\vp)$ the image
of the map $\La^j U\ot\La^j V^\star\to R$ induced by the map
$\La^j\vp$, where $\star$ means the functor $\op{Hom}_R(\cdot,R)$.
If we choose bases for $U$ and $V$, i.e.\ identify $U\simeq R^r$
and $V\simeq R^m$, then the map $\vp$ is represented by an
$m\times r$ matrix $A$ and the ideal $J_j(\vp)\subset R$ is
generated by all $j\times j$ minors of $A$. For the case of pure
order differential operator $\Delta$ of \S\ref{S13} we mean:
$U=ST\ot\nu^*$, $V=ST\ot\pi^*$, $\vp=\z_\Delta$.

Let $G=\op{Coker}(\vp)$. By the Fitting lemma the ideal
$\op{Fitt}_i(G)=J_{m-i}(\vp)$ does not depend on representation
$U\stackrel{\vp}\to V\to G\to 0$ and is called the $i$-th Fitting
invariant of $G$. It is known that
$\op{Fitt}_0(G)\subset\op{ann}(G)$, the two terms have the same
radicals and the equality $\op{Fitt}_0(G)=\op{ann}(G)$ holds if
$\op{depth}\op{ann}(G)=r-m+1$ (see \cite{E, BV}).

We will be interested in the dual over $\R$ map $\vp^*:V^*\to
U^*$, which is the symbol of the collection of differential
operators determining the system $\E$. Thus in this case
$g=\op{Ker}(\vp^*)=\op{Coker}(\vp)^*$, whence $G=g^*$ and
$\op{ann}(G)=I(g)$.

Notice that the characteristic variety of the symbolic system $g$
can be written via $\vp^*:ST^*\ot\pi\to ST^*\ot\nu$ as
$\op{Char}^\C(g)=\{p\in P^\C T^*\,|\,\op{rank}\vp^*(p)<m\}$, where
by $\vp^*(p):\pi\to\nu$ we understand the value of $\vp^*$ at the
covector $p$. Then we define also the {\em characteristic sheaf\/}
(or kernel sheaf; actually it is not a sheaf, but just a family of
vector spaces) over $\op{Char}^\C(g)$ by associating to the
covector $p$ the subspace $\mathcal{K}_p=\op{Ker}\vp^*(p)\subset
\pi$.

 \begin{rk}\po
The last requirement of definition \ref{dfn1}, $\dim K_p=1$ $\forall
p\in\op{Char}^\C(\E)$, means that the system is similar to a system
of scalar PDEs and thus be treated via the usual Koszul complex as
in \cite{KL$_2$}. This, however, cannot be fully formalized, so that
we use another approach with generalized Koszul complexes.
 \end{rk}

\subsection{\hpss Application of the Buchsbaum-Rim complex}\label{S23}

 \abz
Let $F_i$ be some (not necessary linear) differential operators
from a bundle $\pi$ to a bundle $\nu_i$ of $\op{ord}F_i=l(i)$,
$1\le i\le t$. Denote $f_i=\z(F_i)$ the dual symbols (i.e.
$\z_\Delta$) of these operators. These $f_i:U_i\to V$ are
$ST$-homomorphisms of degree $l(i)$, where $U_i=ST\ot\nu_i^*$,
$V=ST\ot\pi^*$. Let $U=\oplus U_i$, $\dim U=\sum \dim U_i=r$.

Consider the map $\vp^*=\oplus f_i^*:V^*\to U^*$, which is the
symbol of differential operator $F=(F_1,\dots,F_t)$. Its $\R$-dual
$\vp=\sum f_i: U\to V$ is a morphism of $R$-modules, but it is not
a graded homomorphism unless the system is of pure order, i.e.
$l(i)\equiv k$. However it becomes homogeneous of degree 1 if we
consider the weighted grading $U\simeq\oplus U_i$, with the weight
$l(i)^{-1}$ for the $i$-th summand.

We wish to find all $R$-relations between the homomorphisms $f_i$.
In other words, we seek to determine the 1-syzygy of the map $\vp:
U\to V$. It is given by the Fitting ideal, but we better describe
the whole free resolution.

This resolution belongs to the Buchsbaum-Eisenbud family of
complexes $\mathcal{C}^i$ (\cite{E}), from which we are interested
in the Buchsbaum-Rim complex $\mathcal{C}^1$:
 $$
0\to S^{r-m-1}V^\star\ot\La^r U\stackrel{\p}\to
S^{r-m-2}V^\star\ot\La^{r-1}U
\stackrel{\p}\to\dots\stackrel{\p}\to
\La^{m+1}U\stackrel{\varepsilon}\to U\stackrel{\vp}\to V.
 $$
Here $\p$ is the multiplication by the trace element $e\in V\ot
V^\star\subset SV\ot\La V^\star$ ($\La V^\star$ acts on $\La U$
via the map $\La\vp^\star$), corresponding to $1\in
R\hookrightarrow V\ot V^\star$:
 $$
\p(a^k\ot b_1\we\dots\we b_t)=\sum (-1)^ik\langle\vp(b_i),a\rangle
\,a^{k-1}\ot b_1\we\dots\check b_i\dots\we b_t.
 $$
The splice map $\varepsilon:\La^{m+s}U\to\La^sU$ is the action of
$\La^m\vp^\star\O\in\La^mU^\star$, where $\O\in\La^mV^\star\simeq
R$ is a volume element (generator).

As proved in \cite{BR} the complex $\mathcal{C}^1$ is exact iff
the map $\vp$ satisfies the condition $\op{depth}J_m(\vp)\ge
r-m+1$.

 \begin{dfn}\po\label{dfn2}
Let us call an $R$-module $G$ {\em generalized complete
intersection\/} if $\op{codim}\op{ann}(G)\ge r-m+1$ (for a minimal
free resolution $U\stackrel{\vp}\to V\to G\to0$).
 \end{dfn}

Note that the usual complete intersections $G=R/I$ satisfy this
property.

 \begin{rk}\po
By the generalized principal ideal theorem of Macaulay (\cite{E})
we have: $\op{codim}\op{ann}(G)=\op{codim}\op{Fitt}_0(G)\le
r-m+1$, so that in fact we have an equality above. In addition, as
we shall see, the module $G$ is Cohen-Macaulay, whenever it is a
generalized complete intersection.
 \end{rk}

 \begin{prop}\po
If an $R$-module $G$ is a generalized complete intersection, then
the Buchsbaum-Rim complex is exact.
 \end{prop}

 \begin{proof}
Let $\vp:U\to V$ be the 1-syzygy map for $G$, $r=\dim U$, $m=\dim
V$. Then we have:
 $$
\op{depth}J_m(\vp)=\op{depth}\op{Fitt}_0(G)=\op{depth}\op{ann}(G)
=\op{codim}\op{ann}(G)=r-m+1,
 $$
where the first and third equalities are general properties of
Noetherian modules, the last one is part of the definition and the
second is a property of Fitting ideal, mentioned in \S\ref{S22}.
Therefore the Buchsbaum-Rim complex is exact.
 \end{proof}

 \begin{rk}\po
Since the polynomial ring $R$ is an affine domain, we have
\cite{BH,E}:
 $$
\op{dim}\bigl(R/\op{ann}(G)\bigr)=\dim
R-\op{codim}\op{ann}(G)=n+m-r-1.
 $$
 \end{rk}

Recall that a ring $P$ is called {\em determinental\/} if
$P=S/Q_s$, where $S$ is a regular Cohen-Macaulay ring and $Q_s$ is
the ideal generated by $s\times s$ minors of an $m\times r$ matrix
$A$ such that the codimension of $Q_s$ in $S$ is exactly
$(m-s+1)(r-s+1)$. By a theorem of Eagon and Hochster such rings
are Cohen-Macaulay \cite{BV}. Let us also call the ideal $Q_s$
itself determinental, if this makes no confusion.

 \begin{theorem}\po\label{090}
Let the symbolic module $G=g^*$ be a generalized complete
intersection in the sense of definition \ref{dfn1} and
$\vp:ST\ot(\oplus\nu_i^*)\to ST\ot\pi^*$  be the corresponding
$R$-homomorphism.  Then we have:
 \begin{enumerate}
 \item The ideal $J_m(\vp)$ is determinental;
 \item $\op{Fitt}_0(G)=I(g)=\op{ann}(G)$;
 \item $G$ is a generalized complete intersection in the sense of
definition \ref{dfn2}.
 \end{enumerate}
 \end{theorem}

 \begin{proof}
Let the conditions of definition \ref{dfn1} be satisfied. It was
shown in \cite{KL$_2$} that if the characteristic sheaf
$\mathcal{K}$ over $\op{Char}^\C(g)$ has fibers of constant
dimension $k$, then $\op{codim}\op{Char}^\C(g)\le l=k(r-m+k)$. When
$k=1$ we get $l=r-m+1$ and this is exactly the codimension of the
characteristic variety of $g$
 $$\op{codim}\op{Char}^\C(g) =r-m+1, $$
determined by $m\times m$ minors, or equivalently by the Fitting
ideal $J_m(\vp)$. Thus we see that the ideal $J_m(\vp)$ is
determinental and $\op{codim}J_m(\vp)=r-m+1$.

This implies that the ring $R/J_m(\vp)$ is Cohen-Macaulay and
$\op{depth}J_m(\vp)=\op{codim}J_m(\vp)$. Since
$\op{Fitt}_0(g^*)=J_m(\vp)$ and $I(g)=\op{ann}(g^*)$ have the same
radicals we have:
 $$
\op{codim}I(g)=\op{codim}J_m(\vp)=\op{depth}J_m(\vp)=r-m+1.
 $$
Thus by the results of \S\ref{S22} we conclude that
$\op{Fitt}_0(g^*)=I(g)$.
 \end{proof}

 \begin{cor}\po
The Buchsbaum-Rim complex $\mathcal{C}^1$ is a resolution of the
symbolic module $g^*$ if the latter is a generalized complete
intersection.
 \end{cor}

 \begin{rk}\po
In \cite{KL$_2$} we also obtained a criterion when the ideal
$\op{Fitt}_0(g^*)=J_m(\psi)$ is a topological complete intersection.
 \end{rk}

\subsection{\hpss Non-linear differential equations}\label{S24}

 \abz
In this section we study non-linear differential equations $\E$. A
system of such equations can be considered as sequence of
submanifolds $\E_k\subset J^k(M)$ with the property
$\E_k^{(1)}\supset\E_{k+1}$ (we assume regularity, but do not
require formal integrability of $\E$).

Let $\F=C^\infty(J^\infty\pi)$ be the filtered algebra of smooth
functions depending on finite jets of $\pi$, i.e. $\F=\cup_i\F_i$
with $\F_i=C^\infty(J^i\pi)$.

Denote $\F^\E_i=C^\infty(\E_i)$. The projections
$\pi_{i+1,i}:\E_{i+1}\to\E_i$ induce the maps
$\pi_{i+1,i}^*:\F_i^\E\to\F_{i+1}^\E$, so that we can form the space
$\F^\E=\cup\F_i^\E$, the points of which are infinite sequences
$(f_i,f_{i+1},\dots)$ with $f_i\in\F_i^\E$ and
$\pi_{i+1,i}^*(f_i)=f_{i+1}$. This $\F^\E$ is a
$C^\infty(M)$-algebra. If the system $\E$ is not formally
integrable, the set of infinite sequences can be void, and the
algebra $\F^\E$ can be trivial. To detect formal integrability, we
investigate the finite level jets algebras $\F^\E_i$ via the
following algebraic approach.

Let $\E$ be defined by a collection $F=(F_1,\dots,F_r)$ of
non-linear scalar differential operators of (possibly different)
orders $l(1),\dots,l(r)$. Each $F_i$ determines a sequence of
smooth maps $J^k(\pi)\to J^{k-l(i)}(\1)$ and so their collection
yields a map $J^\infty(\pi)\to J^\infty(\nu)$, where
$\nu=\oplus\nu_s$ with $\dim\nu_s=m(s)=\{\#i:l(i)=s\}$, $\sum
m(s)=r$.

Pre-composition of our differential operator $F:C^\infty(\pi)\to
C^\infty(\nu)$ with other non-linear differential operators gives
the following exact sequence of $C^\infty(M)$-modules
 \begin{equation}\label{diffE}
\op{diff}(\nu,\1)\stackrel{F}\longrightarrow\op{diff}(\pi,\1)\to\F^\E\to0
 \end{equation}

Note that
$\mathcal{J}_s(F_1,\dots,F_r)=\op{Im}(F)_s\subset\op{diff}_s(\pi,\1)$
is the submodule described in Introduction, and
 \begin{equation}\label{JdE}
\F^\E_i=\op{diff}_i(\pi,\1)/\mathcal{J}_i(F_1,\dots,F_r).
 \end{equation}

On the level of finite jets, the map $F$ of (\ref{diffE})
decreases the order appropriately, but is not homogenous. However
we can adjust this by imposing weights to the vector bundles
$\nu_i$ as we did in \S\ref{S23}. Thus we can assume for
simplicity that the operator $f$ has pure order $k$.

It is important that the terms of (\ref{diffE}) are modules over the
algebra of scalar $\Cc$-differential operators
$\Cc\op{Diff}(\1,\1)$, which are total derivative operators and have
the following form in local coordinates \cite{KLV}: $\Delta=\sum
f_\z\D_\z$, with $f_\z\in C^\infty(J^\infty(M))$. We can identify
$\Cc\op{Diff}(\1,\1)=\cup\F_i\ot\op{Diff}_j(\1,\1)$ with the twisted
tensor product of the algebras $\F$ and $\op{Diff}(\1,\1)$ over the
action
 $$
\hat\Delta:\ \F_i\to\F_{i+j}\quad\text{ for }\quad
\Delta\in\op{Diff}_j(\1,\1).
 $$
This $\Cc\op{Diff}(\1,\1)$ is a non-commutative
$C^\infty(M)$-algebra. We need a more general $\F$-module of
$\Cc$-differential operators $\Cc\op{Diff}(\pi,\1)=\cup
\Cc\op{Diff}_i(\pi,\1)$, where
 $$
\Cc\op{Diff}_i(\pi,\1)=\F_i\ot_{C^\infty(M)}\op{Diff}_i(\pi,\1).
 $$
Remark that $\Cc\op{Diff}(\pi,\1)$ is a filtered
$\Cc\op{Diff}(\1,\1)$-module, i.e.
$\Cc\op{Diff}_i(\1,\1)\cdot\Cc\op{Diff}_j(\pi,\1)\subset
\Cc\op{Diff}_{i+j}(\pi,\1)$.

Define now the filtered $\F^\E$-module $\Cc\op{Diff}^\E(\pi,\1)$
with $\Cc\op{Diff}^\E_i(\pi,\1)=\F^\E_i\ot\op{Diff}_i(\pi,\1)$.
Since the module $\op{Diff}(\pi,\1)$ is projective and we can
identify $\op{diff}(\pi,\1)$ with $\F$, we have from (\ref{JdE})
the following exact sequence
 \begin{equation}\label{0JEF0}
0\to\mathcal{J}_s(F)\ot\op{Diff}_s(\pi,\1)\to
\Cc\op{Diff}_i(\pi,\1)\to\Cc\op{Diff}^\E_i(\pi,\1)\to0.
 \end{equation}

Similar modules can be defined for the vector bundle $\nu$ and
they determine the $\F^\E$-module $\E^*=\cup\E^*_i$ by the
following sequence:
 \begin{equation}\label{diffC}
\Cc\op{Diff}^\E_i(\nu,\1)\stackrel{\ell_F}\longrightarrow
\Cc\op{Diff}^\E_{i+k}(\pi,\1)\to\E^*_{i+k}\to0,
 \end{equation}
where
$\ell:\op{diff}(\pi,\nu)\to\F\ot_{C^\infty(M)}\op{Diff}(\pi,\nu)$
is the operator of universal linearization \cite{KLV},
$\ell_F=\ell(F)$.

This sequence is not exact in the usual sense, but it becomes
exact in the following one. The space to the left is an
$\F^\E_i$-module, the middle term is an $\F^\E_{i+k}$-module. The
image $\ell_F(\Cc\op{Diff}^\E_i(\nu,\1))$ is an $\F^\E_i$-module,
but we generate by it an $\F^\E_{i+k}$-submodule in the middle
term. With this understanding of the image the term $\E^*_{i+k}$
of (\ref{diffC}) is an $\F^\E_{i+k}$-module and the sequence
is exact. In other words
 $$
\E_s^*=\Cc\op{Diff}^\E_s(\pi,\1)/(\F^\E_s\cdot\op{Im}\ell_F).
 $$

Sequences (\ref{diffC}) are nested (i.e. their union is filtered)
and so we have the sequence
 \begin{equation}\label{EEex}
\E_{s-1}^*\to\E^*_s\to\mathcal{F}g_s^*\to0,
 \end{equation}
which becomes exact if we treat the image of the first arrow as
the corresponding generated $\F_s^\E$-module. Thus
$\mathcal{F}g_s^*$ is an $\F_s^\E$-module with support on $\E_s$
and its value at a point $x_s\in\E_s$ is dual to the $s$-symbol of
the system $\E$:
 $$
(\mathcal{F}g_s^*)_{x_s}=g_s^*(x_s);\qquad
g_s(x_s)=\op{Ker}[T_{x_s}\pi_{s,s-1}:T_{x_s}\E_s\to
T_{x_{s-1}}\E_{s-1}].
 $$

This is a geometric definition of the symbol. Equivalently we can
use the algebraic approach as in \S\ref{S13}: Graded space
$g=\oplus g_s$ is dual to cokernel of the symbol $\z_F$ of $F$,
considered as an $ST$-homomorphism $ST\ot\nu^*\to ST\ot\pi^*$,
which depends on the point of equation $\E$.

Our weight-convention apply here and hence we describe the
situation on the level of finite jets $x_i\in\E_i$ for a pure
order $k$ operator $F$, which is the case represented by the
following exact sequence ($x=\pi_i(x_i)\in M$) with the dual
symbol map:
 $$
0\to g_i(x_i)\to
S^iT_x^*M\ot\pi_x\stackrel{\z^F_i(x_i)}\longrightarrow
S^{i-k}T^*_xM\ot\nu_x.
 $$

 \begin{rk}\po
We interpret $\F^\E$ as the algebra of all smooth functions on our
equation $\E$. Define $\E^*=\cup\E^*_i$ in the same manner as
$\F^\E$, taking into account that the map $\E_s^*\to\E^*_{s+1}$ in
our approach is coupled with the change of rings ($\F_s^\E$ to
$\F^\E_{s+1}$). So we can think of $\E^*$ as of sections of the
symbolic bundle $g$ over $\E$ with $\F^\E$-coefficients. Thus
(note linearization in (\ref{diffC})) we interpret $\E^*$ as the
space $\Omega^1(\E)$ of differential forms on $\E$.
 \end{rk}

This remark gives us a way to treat formal integrability of the
system $\E$ as possibility of augmenting exact sequence
(\ref{EEex}) with 0 from the left (injectivity of
$\E_s^*\to\E^*_{s+1}$ under the change of rings). This is reduced
to the question of finding a left resolution of complex
(\ref{diffC}), which will boil down onto the symbolic level as we
shall show.

Thus linearization of the system of generalized complete
intersection type and methods from \S\ref{S23} will lead to the
proof of our compatibility criterion.

\section{\hps Multi-bracket of vector differential operators}

\subsection{\hpss Non-commutative determinants}\label{S31a}

 \abz
Consider the algebra $\mathcal{A}=\op{Diff}(\1,\1)$ of linear
scalar differential operators. It is non-commutative, so no direct
generalization of the determinant function
$\det:\La^m\mathcal{A}^m\to\mathcal{A}$ exists (tensor product is
taken over scalars, not over $\mathcal{A}$). We view the elements
of the space $\mathcal{A}^m$ as rows
$\Delta=(\Delta_1,\dots,\Delta_m)$, which act on columns of
functions $s=(s_1,\dots,s_m)^t\in C^\infty(M)^m$.

We define non-commutative determinant
$\op{Ndet}:\La^m\mathcal{A}^m\to\mathcal{A}$ via the standard
formula, where order of multiplication of matrix elements
corresponds to the order of columns:
 $$
\op{Ndet}\begin{pmatrix}
 \nabla_{11} & \nabla_{12} & \dots & \nabla_{1m}\\
 \vdots      & \vdots      & \ddots & \vdots \\
 \nabla_{m1} & \nabla_{m2} & \dots & \nabla_{mm}
\end{pmatrix}=\sum_{\a\in
S_m}(-1)^\a\nabla_{\a(1)1}\nabla_{\a(2)2}\cdots\nabla_{\a(m)m}.
 $$
In other words, we define non-commutative determinant via
decomposition by columns, i.e. if $C_i(B)$ is the $i^\text{th}$
column of $B\in\op{Mat}_{m\times m}(\mathcal{A})$ and $M_{ij}(B)$
is the minor obtained by removing row $i$ and column $j$, then we
have:
 $$
\!\!\! \op{Ndet}(B)=\sum_{i=1}^m(-1)^{i-1}C_1(B)_i\,
\op{Ndet}(M_{i1}(B))= \sum_{j=1}^m(-1)^{n-j}\op{Ndet}(M_{jm}(B))\,
C_m(B)_j
 $$
(it is however embarrassing to write decomposition via a
mid-column). We obviously have skew-symmetry by rows and
$\R$-linearity, but we lack $\mathcal{A}$-linearity and
skew-symmetry by columns. Thus we can write the non-commutative
determinant in the form
 $$
\op{Ndet}(\nabla_1\we\dots\we\nabla_m).
 $$

Note that the symbol of the non-commutative determinant is the
standard determinant
 $$
\z\bigl(\op{Ndet}(\nabla_1\we\dots\we\nabla_m)\bigr)=
\op{det}\bigl(\z(\nabla_1)\we\dots\we\z(\nabla_m)\bigr),
 $$
where the symbol of order $l$ differential operator
$\nabla_i=(\nabla_{i1},\dots,\nabla_{im})\in\mathcal{A}_l^m$ is
 $$
\z(\nabla_i)=\z_l(\nabla_i)=(\z_l(\nabla_{i1}),\dots,\z_l(\nabla_{im})).
 $$
Beware that since the components of the operator can have smaller
order, it is not true that
$\z(\nabla_i)=(\z(\nabla_{i1}),\dots,\z(\nabla_{im}))$. In other
words, we consider the grading of $\mathcal{A}^m$ corresponding to
increasing filtration $\{\mathcal{A}^m_l\}_{l=0}^\infty$.

Denoting $U^*=\op{Hom}_\R(U,\mathcal{A})$ we have an $\R$-linear map
 $$
\Xi:\La^{m+1}\mathcal{A}^m\to\mathcal{A}^{(m+1)}{}^*
 $$
given
by the formula
 $$
\Xi(\nabla_1\we\dots\we\nabla_{m+1})
=\op{Ndet}\begin{pmatrix}
 \begin{array}{c}\text{\fbox{\ $\hskip8.3pt\nabla_1\hskip8.3pt$\ }} \\ \dots \\
 \text{\fbox{\hskip5pt$\nabla_{m+1}$\hskip5pt}} \end{array} &
 \hskip-8pt\text{\fbox{\vphantom{$\begin{array}{c}\nabla_1\\ \dots\\
\nabla_{m+1}\end{array}$}\hskip6pt}}\
\end{pmatrix},
 $$
where the last column serves as a place-holder, though the result
(image of $\Xi$) we write as as a row.

Notice that the map $\Xi(\nabla_1\we\dots\we\nabla_{m+1})$ is a
right $\mathcal{A}$-homomorphism for all $\nabla_i\in\mathcal{A}$,
i.e.
$\op{Im}\Xi\subset\op{Hom}^\text{right}_\mathcal{A}(\mathcal{A}^{m+1},\mathcal{A})
\subset\mathcal{A}^{(m+1)}{}^*$.

 \begin{rk}\po
Since our constructions are algebraic they can be generalized to
other operator algebras, like pseudo-differential operators, Fourier
operators etc. Then multi-brackets of the next section lead to the
compatibility conditions for the corresponding overdetermined
problems.

One can use the theory of quasi-determinants by Gelfand et al
{\sl\cite{G$^2$RW}} to define other multi-brackets via similar formulas.
However this requires division and extends the class of differential
operators to non-local operators. It could be an exciting relation
between local and global aspects of compatibility.
 \end{rk}

  \comm{
 \begin{rk}\po
Differential operators form a non-commutative algebra without
division. We can introduce a partial division by extending to
pseudo-differential operators (this does not work for all
operators, but many important determinants, like Capelli
determinant will be partial cases).

Then we can use the theory of quasi-determinants
{\sl\cite{G$^2$RW}} to define multi-brackets of the next section.
In this context our theory appeals however to be extended to the
compatibility problem for integral-differential equations, which
is yet a completely unexplored topic.
 \end{rk}
 }

\subsection{\hpss Multi-brackets}\label{S31b}

 \abz
At first we define multi-brackets in the linear case.

Let
$\Upsilon:\underbrace{\mathcal{A}^m\times\dots\times\mathcal{A}^m}_{m+1}\to
\op{Mat}_{(m+1)\times m}(\mathcal{A})$ denote the matrix formed by
$m+1$ vectors-rows from $\mathcal{A}^m$. Then we define the
multi-bracket
 $$
\La^{m+1}\mathcal{A}^m\to\mathcal{A}^m
 $$
of $m+1$ vector differential operators $\nabla_i\in\mathcal{A}^m$
via the operation of the last section and the multiplication
action $\mathcal{A}^{(m+1)}{}^*\times\op{Mat}_{(m+1)\times
m}(\mathcal{A})\to\mathcal{A}^m$ on columns of matrices:
 $$
\{\nabla_1,\dots,\nabla_{m+1}\}=\Xi(\nabla_1\we\dots\we\nabla_{m+1})
\cdot\Upsilon(\nabla_1,\dots,\nabla_{m+1}).
 $$
The $i^\text{th}$ component of the multi-bracket is given by
 $$
\{\nabla_1,\dots,\nabla_{m+1}\}_i=\Xi(\nabla_1\we\dots\we\nabla_{m+1})
\bigl(C_i(\Upsilon(\nabla_1,\dots,\nabla_{m+1}))\bigr).
 $$
It is easy to check that this multi-bracket coincides with the
multi-bracket defined in the introduction. This bracket is
skew-symmetric by its entries and is $\R$-linear. It is not
however $\mathcal{A}$-linear and does not commute with $\R$-linear
transformations of $\mathcal{A}^m$.

To formulate properties of this bracket we will need later an
opposite multi-bracket
$\{\cdots\}^\dag:\La^{m+1}\mathcal{A}^m\to\mathcal{A}^m$, which is
defined by the same formula except that the map $\Xi$ is changed
to $\Xi^\dag$, with the place-holder in non-commutative
determinant being put to the first column.

Note that
 $$
\{\mathcal{A}_{l_1},\dots,\mathcal{A}_{l_{m+1}}\}\subset
\mathcal{A}_{l_1+\dots+l_{m+1}-1},
 $$
where $\mathcal{A}_i$ is the $i$-th subalgebra of the filtered
algebra $\mathcal{A}$.

Let now $F_i\in\op{diff}(m\cdot\1,\1)$, $i=1,\dots,m+1$, be
non-linear differential operators of orders $\op{ord}(F_i)=l(i)$,
which we can identify with smooth functions of the jet-space
space, i.e. elements of $C^\infty(J^{l(i)}(M;\R^m))\subset
\mathcal{F}(J^\infty(M;\R^m))$. Then we can define the
multi-bracket $\{F_1,\dots,F_{m+1}\}$ as an operation
 $$
\La^{m+1}\mathcal{F}(J^\infty(M;\R^m))\to\mathcal{F}(J^\infty(M;\R^m))
 $$
via the linearization operator $\ell: \op{diff}(m\cdot\1,\1)\to
C^\infty(J^\infty(M;\R^m))\ot_{C^\infty(M)}\mathcal{A}^m$, see
\cite{KLV}.

Namely, exploring the formula for the linear case, we let:
 $$
\{F_1,\dots,F_{m+1}\}=\Xi(\ell_{F_1}\we\dots\we\ell_{F_{m+1}})
\cdot\Upsilon(F_1,\dots,F_{m+1}).
 $$
This multi-bracket is related to the multi-bracket of linear
differential operators via the formula
 \begin{equation}\label{lin-mult}
\ell_{\{F_1,\dots,F_{m+1}\}}=\{\ell_{F_1},\dots,\ell_{F_{m+1}}\}.
 \end{equation}

Similarly we can define the opposite multi-bracket
$\{\cdots\}^\dag$ for non-linear differential operators.

\subsection{\hpss Non-commutative "Pl\"ucker identities"}\label{S31c}

 \abz
The multi-bracket we introduced does not satisfy the Jacobi identity
of Nambu \cite{N} (or generalized Poisson) multi-bracket. Neither
does it satisfy the axioms of SH-algebras \cite{LS} (because the
background is different: If we change the length of the
multi-bracket, the functional space changes as well).

However there are certain properties, these brackets do satisfy.
Later in this section we explain that they should be viewed as a
kind of generalized Jacobi identity. For simplicity we begin with
the formulation in the linear case.

 \begin{theorem}\po\label{Plucker}
Let $\nabla_i\in\mathcal{A}^m$ be linear vector differential
operators, $1\le i\le m+2$, and let $\nabla_{i,j}$ denote
component $j$ of $\nabla_i$. Then we have the identities (where as
usual check means absence of argument) relating the multi-bracket
and the opposite multi-bracket for $1\le i\le m$:
 $$
\sum(-1)^k\{\nabla_1,\dots,\check\nabla_k,\dots,\nabla_{m+2}\}^\dag_i
\nabla_k= \sum(-1)^k \nabla_{k,i}
\{\nabla_1,\dots,\check\nabla_k,\dots,\nabla_{m+2}\}.
 $$
 \end{theorem}

 \begin{proof}
Indeed let $\Upsilon_i=C_i(\hat\Upsilon)$ be column $i$ of the
matrix $\hat\Upsilon(\nabla_1,\dots,\nabla_{m+2})$, the map
$\hat\Upsilon:\underbrace{\mathcal{A}^m\times\dots\times\mathcal{A}^m}_{m+2}\to
\op{Mat}_{(m+2)\times m}(\mathcal{A})$ being given by the same
rule as $\Upsilon$ (but with one more row).

Then the right hand side of the identity is obtained by decomposing
the determinant
 $$
\op{Ndet}\begin{pmatrix}
 \hskip2pt\text{\fbox{$\begin{array}{c} \\ \!\!\!\Upsilon_i \\ {}
\end{array}$\hskip-7pt}}\hskip-7pt &
 \begin{array}{c}
 \text{\fbox{\ $\hskip8.3pt\nabla_1\hskip8.3pt$\ }} \\ \dots \\
 \text{\fbox{\hskip5pt$\nabla_{m+2}$\hskip5pt}} \end{array} &
 \hskip-8pt\text{\fbox{$\begin{array}{c} \\ \!\!\!\Upsilon_j \\ {}
\end{array}$\hskip-7pt}}\
\end{pmatrix}
 $$
via the first and then the last column, while the left hand side of
the identity is the result of decomposition by the last and then the
first column. But these operations commute. Now we unite the results
by $j$ to a row.
 \end{proof}

Let $\z\in S_m$ be a permutation and
$T_\z:\mathcal{A}^m\to\mathcal{A}^m$ be the corresponding linear
transformation, $\nabla_i=(\nabla_{i,1},\dots,\nabla_{i,m})\mapsto
T_\z(\nabla_i)=(\nabla_{i,\z(1)},\dots,\nabla_{i,\z(m)})$.

This action leads to conjugated multi-brackets given by
 $$
\{\nabla_1,\dots,\nabla_{m+1}\}^\z=
T_\z^{-1}\{T_\z(\nabla_1),\dots,T_\z(\nabla_{m+1})\}.
 $$
It's easy to see that
 $$
\{\nabla_1,\dots,\nabla_{m+1}\}^\dag_m
=\{T_\tau(\nabla_1),\dots,T_\tau(\nabla_{m+1})\}_1
=\{\nabla_1,\dots,\nabla_{m+1}\}^\tau_m
 $$
for $\t=\left(\begin{array}{cccc} 1 & 2 & \dots & m \\
m & 1 & \dots & m-1 \end{array}\right)$. Thus we get an identity
for the multi-bracket alone:
 \begin{cor}\po\label{multi-Jacobi}
For the above cyclic permutation $\t$ and arbitrary vector
differential  operators $\nabla_i\in\mathcal{A}^m$ it holds:
 $$
\sum(-1)^k\{\nabla_1,\dots,\check\nabla_k,\dots,\nabla_{m+2}\}^\t_m
\nabla_k= \sum(-1)^k \nabla_{k,m}
\{\nabla_1,\dots,\check\nabla_k,\dots,\nabla_{m+2}\}.
 $$
 \end{cor}

We readily generalize the above formulae for non-linear operators
(with the help of linearization operator as in the previous
section). For instance, the latter formula becomes:
 $$
\sum(-1)^k\Bigl(\ell_{\{F_1,\dots,\check F_k,\dots,F_{m+2}\}^\t_m}
F_k- \ell_{F_{k,m}} \{F_1,\dots,\check F_k,\dots,F_{m+2}\}\Bigr)=0.
 $$

Notice that this formula for $m=1$ becomes the standard Jacobi
identity. In this scalar case our multi-bracket becomes bi-bracket
and it coincides with the classical Jacobi bracket $\{F,G\}$ of
scalar (non-linear) differential operators
$F,G\in\op{diff}(\1,\1)$ (in the linear case $F,G\in\mathcal{A}$
it is the commutator, for the non-linear case see \cite{KLV}).
Indeed the formula is:
 $$
\sum_\text{cyclic}\bigl(\ell_F\{G,H\}-\ell_{\{G,H\}}F\bigr)=
\sum_\text{cyclic}\{F,\{G,H\}\}=0.
 $$

Thus the multi-bracket identities could be considered as
generalized Jacobi identities (but not in the sense of
\cite{APB}). We however call them non-commutative Pl\"ucker
identities by the following reason. Consider for simplicity the
case $m=2$.

In this case the multi-bracket of operators
$\nabla_i=(\nabla_{i,1},\nabla_{i,2})$, $i=1,2,3$ correspond to
the composition $\vp_0\vp_1$ of the $\mathcal{A}$-homomorphisms
 $$
0\to\mathcal{A}\stackrel{\vp_1}\longrightarrow\mathcal{A}^3
\stackrel{\vp_0}\longrightarrow\mathcal{A}^2,
 $$
where (the determinant is $\op{Ndet}$)
 $$
\vp_1=\left(\begin{vmatrix} \nabla_{2,1}\!\! & \!\!\nabla_{2,2} \\
\nabla_{3,1}\!\! & \!\!\nabla_{3,2} \end{vmatrix},
\begin{vmatrix} \nabla_{3,1}\!\! & \!\!\nabla_{3,2} \\
\nabla_{1,1}\!\! & \!\!\nabla_{1,2} \end{vmatrix},
\begin{vmatrix} \nabla_{1,1}\!\! & \!\!\nabla_{1,2} \\
\nabla_{2,1}\!\! & \!\!\nabla_{2,2} \end{vmatrix}\right),\
\vp_0=\left(\begin{array}{cc} \nabla_{1,1} & \nabla_{1,2} \\
\nabla_{2,1} & \nabla_{2,2} \\ \nabla_{3,1} & \nabla_{3,2}
\end{array}\right).
 $$
The above sequence is not a complex (whence the multi-bracket),
but its symbolic part is a complex and is actually a resolution of
the module $g^*$ corresponding to the system
$\{\nabla_1[u]=0,\nabla_2[u]=0,\nabla_3[u]=0\}$.

To perceive the properties of the multi-brackets we consider 4
vector differential operators $\nabla_i\in\mathcal{A}^2$, $1\le
i\le 4$, and extend the above complex to
 $$
0\to\mathcal{A}^2\stackrel{\phi_2}\longrightarrow\mathcal{A}^4
\stackrel{\phi_1}\longrightarrow\mathcal{A}^4
\stackrel{\phi_0}\longrightarrow\mathcal{A}^2,
 $$
where $\phi_0$ is a $4\times 2$ matrix with rows $\nabla_i$,
$\phi_2$ is the transposed matrix and $\phi_1$ is a skew-symmetric
$4\times 4$ matrix with entries
$\begin{vmatrix} \nabla_{i,1}\!\! & \!\!\nabla_{i,2} \\
\nabla_{j,1}\!\! & \!\!\nabla_{j,2} \end{vmatrix}$ being
non-commutative determinants.

Again the symbolic sequence is exact, but the general sequence is
not a complex and the composition $\phi_0\phi_1\phi_2$ gives us
the desired properties of the multi-bracket. Clearly the above
"resolution" is built on a certain determinental identity, which
is exactly the Pl\"ucker identity in Grassmaninan $G(2,4)$.

For $m>2$ we see that our non-commutative identities model the
standard Pl\"ucker identities in other Grassmaninans.

\subsection{\hpss Reduced brackets and coordinates}\label{S31d}

 \abz
Let $\E=\{F_1[u]=0,\dots,F_r[u]=0\}$ be an overdetermined system
of PDEs, where $F_i\in\op{diff}(m\cdot\1,\1)$ and
$u=(u_1,\dots,u_m)^t\in C^\infty(M,\R^m)$. As in the introduction
we denote by $\mathcal{J}_s(\E)=\langle \ell_\Delta\circ
F_i\,|\,\op{ord}(\Delta)+\op{ord}{F_i}\le s\rangle$ the submodule
generated by $F_1,\dots,F_r$.

We let $\E^*_s=\mathcal{A}^m_s/\mathcal{J}_s(\E)$ in the linear
case and $\E^*_s=\op{diff}_s(m\cdot\1,\1)/\mathcal{J}_s(\E)$ in
the non-linear. In this way we obtain the reduced bracket
 $$
[f_1,\dots,f_{m+1}]_\E= \{f_1,\dots,f_{m+1}\}\,\op{mod}
\mathcal{J}_{l-1}(\E)\in\E_{l-1}^*.
 $$
for $l=l(f_1)+\dots+l(f_{m+1})$.

 \begin{rk}\po
This multi-bracket appears due to the fact, that we do not have a
unique non-commutative determinant. If we consider determinants
with the values in the reduced (quotient) module, as it is done in
the case of Dieudonne determinant, we will arrive to this reduced
multi-bracket.
 \end{rk}

Every $k$-th order scalar differential operator
$G\in\op{diff}_k(\cdot\1,\1)$ induces via linearization a map
$\hat G:\E_s^*\to\E_{s+k}^*$. With respect to this map Theorem
\ref{Plucker} implies the reduced identities:
 \begin{theorem}\po
Let $\E$ be an over-determined system ($r>m$) defined by
(non-linear) differential operators $F_1,\dots,F_r$. Then for any
subset $\{i_1,\dots,i_{m+2}\}\subset\{1,\dots,r\}$ and any
$j\in[1,m]$ we have:
 $$
 \sum(-1)^k \hat F_{i_k,j}
[F_{i_1},\dots,\check F_{i_k},\dots,F_{i_{m+2}}]_\E=0.
 \text{\phantom{\underline{y}}\hfill\hfill}
 $$
\vskip-32pt\qed
 \end{theorem}

Notice that we do not assume integrability and so the
multi-brackets occurring in the Main Theorem are not arbitrary,
but vanishing of some of them gives certain restrictions for the
rest.

 \begin{rk}\po
In general it is not true that
 $$
\{\mathcal{J}_{l_1}(\E),\dots,\mathcal{J}_{l_{m+1}}(\E)\}\subset
\mathcal{J}_{l_1+\dots+l_{m+1}-1}(\E),
 $$
but for formally integrable $\E$ it is. Then we can define the
bracket
 $$
[\E^*_s,\mathcal{J}_{l_1}(\E),\dots,\mathcal{J}_{l_m}(\E)]_\E
\subset\E^*_{l_1+\dots+l_m+s-1}.
 $$
Then elements $\theta\in\E_s^*$ such that
$[\theta,\mathcal{J}_{l_1}(\E),\dots,\mathcal{J}_{l_m}(\E)]_\E=0$
with respect to this bracket, can be interpreted as another
generalization of the classical notion of symmetry.
 \end{rk}

Finally we can give a coordinate representation of the introduced
multi-bracket. For calculational purposes it is however more
convenient to work with the following multi-bracket:
 $$
[F_1,\dots,F_{m+1}]=\dfrac1{m!}\!\!\sum_{\begin{array}{c}
\scriptstyle \z\in S_{m+1}\\ \scriptstyle \nu\in S_m \end{array}}
\!\!\!\dfrac{\op{sgn}(\z)}{\op{sgn}(\nu)}\hspace{-5pt}
\sum_{\begin{array}{c} \scriptstyle 1\le i\le m \\ \scriptstyle
|\t_i|=l(\z(i))
\end{array}}\hspace{-5pt}
\prod_{j=1}^m\dfrac{\p F_{\z(j)}}{\p p^{\nu(j)}_{\t_j}}\
\D_{\t_1+\dots+\t_m}F_{\z(m+1)},
 $$
where $F_i\in\op{diff}_{l(i)}(m\cdot\1,\1)$. For $m=1$ this gives
Mayer brackets instead of Jacobi brackets \cite{KL$_1$}. The
following statement is straightforward:

 \begin{prop}\po
Restrictions of the two multi-brackets to the system $\E$
coincide:
 $$
[F_1,\dots,F_{m+1}]_\E\equiv[F_1,\dots,F_{m+1}]\mod
\mathcal{J}_{l-1}(\E),
 $$
where $l=\sum_{i=1}^{m+1}l(i)$ as before. \qed
 \end{prop}

\section{\hps Compatibility criterion}

 \abz
In this section we prove our compatibility criterion. Its particular
cases are theorems from \cite{KL$_1$,KL$_2$,KL$_3$}, where we used
geometric theory of PDEs and the obstructions to compatibility were
identified with certain curvatures (Weyl tensors). Here we propose
an approach based on the construction of symbolic compatibility
complex, which uses the dual algebraic approach.

\subsection{\hpss Syzygies for modules of linear differential operators}\label{S32}

 \abz
Consider the filtered $\mathcal{A}$-module $\op{Diff}(\pi,\1)$ and
let
 $$
\op{Diff}(\nu,\1)\stackrel{\phi^\Delta}\longrightarrow
\op{Diff}(\pi,\1)\to\E^*\to 0
 $$
be a representation of the dual to a system of linear equations
$\E$. Set $\mathcal{I}=\op{Im}\bigl(\phi^\Delta\bigr)$. In other
words, we let $\mathcal{I}_k\subset\op{Diff}_k(\pi,\1)$ denote the
sequence of submodules
$\mathcal{J}_k(\E)=\mathcal{J}_k(F_1,\dots,F_r)$ as in Introduction,
$F_i\in\op{Diff}_{l(i)}(\pi,\1)$. Notice that these $\mathcal{I}_k$
define our equations: $\E_k=\{x_k:h(x_k)=0\,\forall
h\in\mathcal{I}_k\} \subset J^k(\pi)$.

 \begin{prop}\po\label{567}
The formal integrability of the system $\E$ is equivalent to the
requirement that
$\mathcal{I}_{k+1}\cap\op{Diff}_k(\pi,\1)=\mathcal{I}_k$ for all
$k$.
 \end{prop}

 \begin{rk}\po
The last condition means it is not possible to get new
relations of order $k$ in $\mathcal{I}$ from the relations of
order $k+1$ via linear combinations over $\mathcal{A}_{k+1}$.
 \end{rk}

 \begin{proof}
In fact, surjectivity of the projections
$\pi_{k+1,k}:\E_{k+1}\to\E_k$ (formal integrability) is equivalent
to injectivity of the dual maps
$\pi_{k+1,k}^*:\E_k^*\to\E_{k+1}^*$. The claim follows from the
natural isomorphism
$\E_k^*\simeq\op{Diff}_k(\pi,\1)/\mathcal{I}_k$.
 \end{proof}

Now if $\mathcal{I}$ is generated by linear differential operators
$F_1,\dots,F_r$, every element in $\I_{k+1}$ is represented in the
form $\sum A_i F_i$, where $A_i=\sum a_i^\t\D_\t$ are scalar
differential operators of $\op{ord}(A_i)\le k+1-\op{ord}(F_i)$.
The condition of proposition \ref{567} is equivalent to the
following relation on the $k$-symbols:
 \begin{equation}\label{eq3}
\z\bigl(\sum A_iF_i\bigr)=\sum \z(A_i)\z(F_i)=0\in S^{k+1}T\ot
N^*.
 \end{equation}
To describe all such relations in the submodule of symbolic
relations $I\subset ST\ot N^*$ we use the syzygy approach of
\S\ref{S23}.

We also denote $V^*=ST^*\ot \pi$ and $U^*=ST^*\ot\nu$, where the
bundle $\nu=\oplus\nu_i$ is, in general, graded by the degrees of
operators $\Delta_i$. We have symbols $f_i=\z(F_i)$ (previously
denoted $\z^{F_i}$) of our differential operators $F_i$,
$i=1,\dots,r$. The map $\psi=\vp^*:V^*\to U^*$, which can be
represented in bases as $\psi:R^m\to R^r$,
$(u^1,\dots,u^m)\mapsto(f^1,\dots,f^r)$, has the kernel $g\subset
ST^*\ot N$.

If $\Pi=\Pi(g)\subset ST\ot N^*$ is the annihilator submodule, then
$ST\ot N^*/\Pi(g)$ is the symbolic module $g^*=\op{Coker}(\vp)$.
Consider $\R$-dual to the Buchsbaum-Rim complex from \S\ref{S23}
(as the Spencer complex is $\R$-dual to the Koszul complex):
 $$
0\to g\to V^*\stackrel{\psi}\to
U^*\stackrel{\omega}\to\La^{m+1}U^*\stackrel{\d}\to
V^\times\ot\La^{m+2}U^*\stackrel{\d}\to\dots
S^{r-m-1}V^\times\ot\La^rU^* \to0,
 $$
where $\d=\p^*$, $\oo=\varepsilon^*$ and $V^\times=(V^\star)^*$.
Choosing a basis $e_1,\dots,e_m$ of $V^*\simeq R^m$ we can
describe informally $\oo(\x)=\x\we\psi(e_1)\we\dots\we\psi(e_m)$.
However as the symbolic differential operator $\vp$ decreases
degrees, we need to change it to $\oo(h_1\dots h_m\x)=\x\we
\psi(h_1e_1)\we\dots\we\psi(h_me_m)$ for certain elements $h_i\in
R$ of sufficiently high degrees.

Thus all the relations between the symbolic differential operators
$\z(\Delta_i)$ are given by the explicit formula $\oo\circ\psi=0$
from the above complex, which we eventually call generalized
Spencer $\d$-complex.

\subsection{\hpss Proof of the main theorem for linear systems}\label{S33}

 \abz
At first let us consider two partial cases.

{\bf 1.} {\it Scalar equations}. In this case $m=1$ and the
condition of definition \ref{dfn1} says that $g^*=ST/I(g)$ is the
usual complete intersection. So $I(g)=\langle f_1,\dots,
f_r\rangle$, where $\{f_i\}$ form a regular sequence of length
$r\le n$. Then the Koszul complex
 $$
0\to I\ot \La^r R^r\stackrel{\p}\to \dots\to
I\ot\La^2R^r\stackrel{\p}\to I\ot R^r\to I\to 0
 $$
is exact. In particular, vanishing of the 1st homology yields:
 $$
\sum_{i=1}^r a_if_i=0\ \Longrightarrow\ a_i=\sum_{j=1}^r
c_{ij}f_j,\quad c_{ij}+c_{ji}=0\qquad (a_i,c_{ij}\in R),
 $$
and so the 1-syzygy module is generated by the relations
$f_if_j-f_jf_i=0$. Thus we need to check the condition of
proposition \ref{567} only for combinations $F_iF_j-F_jF_i$
(multiplication in algebra $\mathcal{A}$ is non-commutative), see
(\ref{eq3}).

As a consequence we obtain that the commutator $[F_i,F_j]$, which
is an operator of order $l_{ij}=l(i)+l(j)-1$ should belong to the
space $\mathcal{J}_{l_{ij}}$ generated by $F_1,\dots,F_r$ and
their total derivatives up to the order $l_{ij}$. This is the
compatibility condition for the system $\E$, exactly as theorem
\ref{th A} states.

{\bf 2.} {\it Systems on two-dimensional manifolds}. In this case
$n=2$ and condition 1 of definition \ref{dfn1} gives $r=m+1$. Now
instead of Koszul complex we use the following approximation to a
resolution
 $$
0\to g\to R^m\stackrel{\psi}\to R^{m+1}\stackrel{\t}\to\hat I\to
0,
 $$
where $\hat I=aJ_m(\psi)\subset R$ for a non-zero divisor $a$. Let
$f_i:R^m\to R$ be symbols of the defining equations for the system
$\E$ and $f_i(e_j)=f_{ij}$ their values on a basis. Denote by
$A(f)=\|f_{ij}\|$ the $(m+1)\times m$ matrix of the operator
$\psi$. Then the map $\t$ is given by the formula
$\t(\x)=\sum_{i=1}^{m+1}(-1)^{i-1}\x_i\det A_i(f)$, with $A_i(f)$
is obtained from $A(f)$ by deleting the $i$-th row (Laplace
decomposition).

The Hilbert-Burch theorem states that the above complex is exact
whenever $\op{depth}J_m(\psi)\ge 2$. By theorem \ref{090} this
follows from the conditions of Definition \ref{dfn1} (in
particular $\op{Char}^\C(g)=\emptyset$) and also we see that we
have equality. Thus the only generator of 1-syzygy is the
relation:
 $$
\sum_{i=1}^{m+1}(-1)^{i-1}\det A_i(f)f_i=0.
 $$
Let $\mathcal{A}_i(f)$ be some differential operators with the
symbols $A_i(f)$. As a consequence of the 1-syzygy description we
obtain that the expression
 \begin{equation}\label{eq4}
\sum_{i=1}^{m+1}(-1)^{i-1}\det \mathcal{A}_i(f)F_i,
 \end{equation}
which is an operator of order $l=l(1)+\dots+l(r)-1$, should belong
to the space $\mathcal{J}_{l}$ generated by $F_1,\dots,F_r$ and
their total derivatives up to the order $l$. Noticing that the
expression in (\ref{eq4}) is the multi-bracket
$\{F_1,\dots,F_r\}$, we get the compatibility condition from
theorem \ref{th A}.

Now, as we have clarified the simple situations, we will study

{\bf 3.} {\it The general linear case}.

Here we should use the generalized Spencer complex described in
\S\ref{S32} ($\R$-dualization of the complex $\mathcal{C}^1$ from
\S\ref{S23}). If the assumptions of definition \ref{dfn1} are
satisfied, we conclude that the relations generating the 1-syzygy
are of the type (\ref{eq4}), namely
 $$
\sum_{j=1}^{m+1}(-1)^{j-1}\det A_{\varrho,i_j}(f)f_{i_j},
 $$
where $A_\varsigma(f)$ is obtained from the $r\times m$ matrix
$A(f)$ of the map $\vp$ by deleting rows with the numbers
$\varsigma\subset\{1,\dots,r\}$ and $\varrho$ is the complement to
the subset\linebreak
 $\{i_1,\dots,i_{m+1}\}$, so that the
resulting matrix $A_{\varrho,i}(f)$ is square of size $m\times m$.

It's clear that changing
$A_{\varrho,i_j}(f)\mapsto\mathcal{A}_{\varrho,i_j}(f)$ and
$f_{i_j}\mapsto F_{i_j}$ in the above expressions we obtain the
multi-brackets $\{F_{i_1},\dots,F_{i_{m+1}}\}$ (or
$[F_{i_1},\dots,F_{i_{m+1}}]$ depending on the manner we extend
the symbol $A_{\varrho,i_j}(f)$ to a differential operator
$\mathcal{A}_{\varrho,i_j}(f)$) for various (ordered) subsets
$\{i_1,\dots,i_{m+1}\}\subset\{1,\dots,r\}$, which by proposition
\ref{567} should belong to the subspace
$\mathcal{J}_{l(i_1)+\dots+l(i_{m+1})-1}$ iff $\E$ is formally
integrable. The claim is proved.

\subsection{\hpss The proof for non-linear operators and generalization}\label{S34}

 \abz
Now consider the general case of non-linear systems $\E$. The main
theorem in this case can be proved as follows: We linearize the
system. The new linear system is of generalized complete
intersection type as well. Compatibility of non-linear system $\E$
is equivalent to compatibility for each linearization on a
jet-solution. The multi-brackets are also preserved (this is
well-studied in the case of usual brackets, see \cite{KLV} for a
relation between universal linearization operators, Jacobi
brackets and evolutionary differentiations), see (\ref{lin-mult}).

We however will not develop these vague ideas, but show instead
how to modify the proof from \S\ref{S33}. Let the system $\E$ be
given by a (matrix) differential operator
$\Delta=(F_1,\dots,F_r):C^\infty(\pi)\to C^\infty(\nu)$. The
bundle $\nu$ can be locally trivialized ($\R^r$ above) or more
generally we can split it $\nu=\oplus\nu_s$ according to different
orders of the differential operators $F_i$. Since we showed in
\S\ref{S23} how to reduce this case to the pure order by imposing
weights on $\nu_s$, we just assume the operator is of pure order
$k$ (can be set 1) $\Delta:J^{i}(\pi)\to J^{i-k}(\nu)$. Its symbol
we denote by $\vp$.

Let $t=i-k(m+1)$. Denote by $\La^{m+1}_{[t]}$ the $t$-th graded
component of the $(m+1)$-st exterior degree of the corresponding
module in the context of symbolic systems and the $t$-th filtered
submodule for the case of differential operators (with the above
weight-convention). Let us also use the short notations
$\op{Diff}_i(\pi)=\op{Diff}_i(\pi,\1)$,
$\op{Diff}_i(\nu)=\op{Diff}_i(\nu,\1)$ [only for large diagrams].

The following diagram commutes:

 $$
 \begin{CD}
 0 @. 0 @. 0 @. 0\\
 @AAA @AAA @AAA @AAA @. \\
\hspace{-4pt}
\La^{m+1}_{[t]}ST\ot\nu^*@>{\varepsilon}>>S^{i-k}T\ot\nu^*
@>{\vp^*}>>S^iT\ot\pi^*@>>>g_i^*@>>>0\\
 @AAA @AAA @AA{\varrho}A @AAA @. \\
\hspace{-7.5pt}
\La^{m+1}_{[t]}\op{Diff}(\nu)@>{\nabla^*_t}>>\op{Diff}_{i-k}(\nu)
@>{\Delta^*}>>\op{Diff}_i(\pi)@>>>\E_i^*@>>>0\\
 @. @AAA @A{j}AA @AA{(\pi_{i,i-1}^\E)^*}A @. \\
\hspace{-12pt}  @. \op{Diff}_{i-k-1}(\nu)
@>{\Delta^*}>>\op{Diff}_{i-1}(\pi)@>>>\E_{i-1}^*@>>>0\\
 @. @AAA @AAA @AAA @. \\
 @. 0 @. 0 @. \hspace{-16pt}\op{Ker}(\pi_{i,i-1}^\E)^*\hspace{-20pt}\\
 \end{CD}
 $$
 \vspace{1pt}\hspace{-15pt}
 \begin{picture}(0,0)
 {\thicklines\dottedline[$\cdot$]{5}(63,133)(105,112)}
 \put(107,111){\vector(2,-1){5}}
 {\thicklines\dottedline[$\cdot$]{5}(155,87)(191,69)}
 \put(193,68){\vector(2,-1){5}}
 {\thicklines\dottedline[$\cdot$]{5}(235,47)(271,29)}
 \put(273,28){\vector(2,-1){5}}
 \end{picture}

The first row is a part of the Buchsbaum-Rim complex $\mathcal{C}^1$
and the operator $\nabla_t^*$ in the second one is a differential
operator with the symbol $\varepsilon$.

 \begin{prop}\po\label{prp10}
Let the system $\E$ be generalized complete intersection. Its
compatibility on the level of $(i-1)$ jets, i.e.
$\op{Ker}(\pi_{i,i-1}^\E)^*=0$, is equivalent to existence of an
operator $\nabla_t^*$, such that the second row is a complex.
 \end{prop}

 \begin{proof}
By the assumption the first row is exact. The diagram chase, not
involving $\nabla_t^*$, yields a homomorphism
 \begin{equation}\label{897}
\zeta_i:\La^{m+1}_{[t]}ST\ot\nu^*/\op{Ker}\varepsilon
\simeq\op{Im}\varepsilon\longrightarrow
\op{Ker}(\pi_{i,i-1}^\E)^*.
 \end{equation}
It is always an epimorphism, but it is a monomorphism iff the map
$\Delta^*\circ\nabla_t^*$ is a monomorphism (the precise form of
$\nabla_t^*$ is inessential here, we may consider its construction
to be inductive by $t$ or refer to our variety of multi-brackets).
Vanishing of our homomorphism is equivalent to vanishing of
$\Delta^*\circ\nabla_t^*$.
 \end{proof}

Since $\varrho\circ\Delta^*\circ\nabla_t^*=0$, the map
$\Delta^*\circ\nabla_t^*:
\La^{m+1}_{[t]}\op{Diff}(\nu,\1)\to\op{Diff}_i(\pi,\1)$ from the
proof can be identified as follows:
 $$
j^{-1}\circ\Delta^*\circ\nabla_t^*:
\La^{m+1}_{[t]}\op{Diff}(\nu,\1)\to\op{Diff}_{i-1}(\pi,\1),
 $$
which can be varied by a map of the form
 $$
\Delta^*\circ\nabla_{t-1}^*:
\La^{m+1}_{[t-1]}\op{Diff}(\nu,\1)\to\op{Diff}_{i-1}(\pi,\1).
 $$

Since every differential operator from $\op{Diff}(\pi,\1)$
determining our system (i.e. vanishing on $\E$) can be represented
by operators from $\op{Diff}(\nu,\1)$ composed with the operator
$\Delta^*$, the above operator represents our multi-bracket
(indeed one of them due to non-uniqueness of $\nabla_t^*$), or
more exactly by the collection of reduced multi-brackets
$[F_{i_1},\dots,F_{i_{m+1}}]_\E=\{F_{i_1},\dots,F_{i_{m+1}}\}
\op{mod}\mathcal{J}'_{i-1}$ (with an appropriate index $i$).

The condition $\op{Ker}(\pi_{i,i-1}^\E)^*=0$ $\forall i$ is
equivalent to vanishing of the maps
$j^{-1}\circ\Delta^*\circ\nabla_t^*$ on
$\La^{m+1}_{[t]}\op{Diff}(\nu,\1)$ for all $t$. Thus formal
integrability is equivalent to vanishing of all multi-brackets of
$(m+1)$-tuples of differential operators $F_1,\dots,F_r$ due to
the system $\E$ (as explained in theorem \ref{th A}).

Now this line of arguments works well for the case of non-linear
vector differential operators $F=(F_1,\dots,F_r)$:

 \medskip

 \begin{proof}[Proof of theorem \ref{th A}]
The above diagram over $M$ should be considered\- over the
equation $\E$ and the $\op{Diff}(\1,\1)$-module
$\op{Diff}(\pi,\1)$ of linear differential operators should be
changed to the $\Cc\op{Diff}^\E(\1,\1)$-module
$\Cc\op{Diff}^\E(\pi,\1)$ of $\Cc$-differential operators with
coefficients in $\F^\E$ (see \S\ref{S24}).

From the arising diagram of non-linear complexes we also obtain
the multi-brackets as the obstructions to formal integrability and
observe that their vanishing due to the system ensures this
integrability.

Namely when the conditions of definition \ref{dfn1} are satisfied,
we obtain the following sequence of non-linear differential
operators, which is exact at the terms ${\Cc}\op{Diff}^\E(\pi,\1)$
and $\E^*$ (in the sense specified in \S\ref{S24}):
 \begin{equation}\label{nlnEF}
\La^{m+1}{\Cc}\op{Diff}^\E(\nu,\1)\stackrel{\nabla^*}
\longrightarrow{\Cc}\op{Diff}^\E(\nu,\1)
\stackrel{\ell_F}\longrightarrow
{\Cc}\op{Diff}^\E(\pi,\1)\longrightarrow\E^*\to0.
  \end{equation}

The composition $\ell_F\circ\nabla^*$ corresponds to
multi-brackets of non-linear differential operators and this
sequence is a complex iff the brackets are zero due to the system.
Thus vanishing of the multi-brackets yields formal integrability.

Let us give more details. Compatibility on the level of $i$-jets,
which is equivalent to injectivity of the map
$\pi_{i,i-1}^*:\F^\E_{i-1}\to\F^\E_i$, can be expressed by saying
that no relation of functions from $\mathcal{J}_i(F_1,\dots,F_r)$
over $\F_i$ is an operator from $\op{diff}_{i-1}(\pi,\1)$, see
(\ref{diffE}).

Linearization reduces this to the following claim (cf. Proposition
\ref{567}):
 $$
\Cc\op{Diff}^\E_i(\pi,\1)\cap\Cc\op{Diff}_{i-1}(\pi,\1)=
\Cc\op{Diff}^\E_{i-1}(\pi,\1).
 $$
Thus due to exact sequences (\ref{0JEF0}) and (\ref{diffC}) the
formal integrability can be expressed via 1-syzygy of the module
$\E^*$, which is given by sequence (\ref{nlnEF}).

Indeed any non-linear relation can be evaluated at points of the
equation $x_i\in\E_i$, which gives a relation for the symbolic
system $g(x_i)$. Since all such relations are described via the
Buchsbaum-Rim complex, we describe 1-syzygy of $\E^*$ in exactly
the same way as in \S\ref{S32}.

This yields multi-brackets $[F_{i_1},\dots,F_{i_{m+1}}]_\E$ as
obstructions to compatibility and theorem \ref{th A} is proved.
 \end{proof}

{\bf A generalization.} What happens to the compatibility
criterion, when the system is not of generalized complete
intersection type?

Then the first row of the diagram from \S\ref{S34} needs not to be
exact and homomorphism (\ref{897}) does not tell much. In
proposition \ref{prp10} we assumed the system is of generalized
complete intersection type, which gives
$\op{Im}\varepsilon=\op{Ker}\vp^*$.

If this does not hold, we get a cohomology group $H^{i-k}_{\nu^*}$
at the term $S^{i-k}T\ot\nu^*$ of the Buchsbaum-Rim complex
$\mathcal{C}^1$. Then homomorphism (\ref{897}) will be changed to
 $$
\hat\zeta_i=\zeta_i\oplus\zeta_i':\op{Ker}\vp^*=
\op{Im}\varepsilon\oplus H^{i-k}_{\nu^*}\longrightarrow
\op{Ker}(\pi_{i,i-1}^\E)^*.
 $$
The first component $\zeta_i$ is represented by multi-brackets of
differential operators as before, but the second component
$\zeta_i'$ is of completely different nature. Again the map
$\hat\zeta_i$ is an epimorphism, but neither $\zeta_i$ nor
$\zeta_i'$ needs to be a monomorphism.

Nevertheless, since the multi-brackets together with
$\op{Im}\zeta_i'$ span $\op{Ker}(\pi_{i,i-1}^\E)^*$, they give all
compatibility conditions of the system $\E$.

The image of the projection $\pi_{i,i-1}^\E:\E_i\to\E_{i-1}$ is the
locus of the Weyl tensors $W_{i-1}(\E)\in H^{i-2,2}(\E)$ (we will
not describe them here, see \cite{L$_1$}; for geometric structures
they are also called structural functions \cite{St}, curvature
tensors and sometimes torsion tensors \cite{BCG$^3$}). We let
$W(\E)=\oplus W_i(\E)$. Note that the relation $W(\E)=0$ is
precisely the set of all compatibility conditions.

Epimorphism $\hat\zeta_i$ leads to the decomposition
 $$
\op{Ker}(\pi_{i,i-1}^\E)^*\simeq
(\zeta_i\circ\varepsilon)(\La^{m+1}_{[t]}ST\ot\nu^*)\oplus
\zeta_i'(H^{i-k}_{\nu^*}).
 $$
and so we get:

 \begin{theorem}\po
The Weyl tensor can be decomposed $W(\E)=W_B(\E)+W_H(\E)$. The
bracket part $W_B(\E)$ can be expressed via multi-brackets by the
same formula as in corollary \ref{cor2}. The second part $W_H(\E)$
is a homological term. \qed
 \end{theorem}

\subsection{\hpss An exact sequence of non-linear differential operators}

 \abz
We would like now to make the complex from the proof of theorem
\ref{th A} into an exact sequence by interchanging the functors
$\La^{m+1}$ and $\op{Diff}$:

 \begin{theorem}\po\label{exact_seq_DO}
Assume that the system $\E$ is of generalized complete
intersection type and the multi-brackets vanish as in theorem
\ref{th A}. Then there is an exact complex of differential
operators, which is a resolution of the module $\E^*$:
 \begin{multline}\label{975}
\dots\to{\Cc}\op{Diff}^\E(\La^{m+3}\nu\ot S^2\pi^*,\1)
\stackrel{\square_2^*}\longrightarrow{\Cc}\op{Diff}^\E(\La^{m+2}\nu\ot\pi^*,\1)
\stackrel{\square_1^*}\longrightarrow\\
\to{\Cc}\op{Diff}^\E(\La^{m+1}\nu,\1)
\stackrel{\nabla^*}\longrightarrow{\Cc}\op{Diff}^\E(\nu,\1)
\stackrel{\ell_F}\longrightarrow{\Cc}\op{Diff}^\E(\pi,\1)\longrightarrow\E^*\to0.
 \end{multline}
 \end{theorem}

 \begin{proof}
Consider at first the linear case, where we shall use the notation
$\Delta=(F_1,\dots,F_r)$ for a given vector-valued operator,
defining $\E$.

We have the following diagram with the exact columns:
 $$
 \begin{CD}
 0 @. 0 @. 0 @. 0\\
 @AAA @AAA @AAA @AAA @. \\
\hspace{-7pt} \stackrel{\d}\to
\La^{m+1}_{[t]}ST\ot\nu^*@>{\varepsilon}>>S^{i-k}T\ot\nu^*
@>{\vp^*}>>S^iT\ot\pi^*@>>>g_i^*@>>>0\\
 @A{q_3\downarrow}A{\rho_3}A @A{q_2\downarrow}A{\rho_2}A
@A{q_1\downarrow}A{\rho_1}A @AAA @. \\
\hspace{-10.5pt} \dots\to\!\! \
\op{Diff}_t(\La^{m+1}\nu)@>{\nabla^*_t}>>\op{Diff}_{i-k}(\nu)
@>{\Delta^*_i}>>\op{Diff}_i(\pi)@>{\a_i}>>\E_i^*@>>>0\\
 @AA{j_3}A @AA{j_2}A @AA{j_1}A @AA{(\pi_{i,i-1}^\E)^*}A @. \\
\hspace{-15pt}
\dots\op{Diff}_{t-1}(\La^{m+1}\nu)@>{\nabla^*_{t-1}}>>\op{Diff}_{i-k-1}(\nu)
@>{\Delta^*_{i-1}}>>\op{Diff}_{i-1}(\pi)@>{\a_{i-1}}>>\E_{i-1}^*@>>>0\\
 @AAA @AAA @AAA @AAA @. \\
 0 @. 0 @. 0 @. 0\\
 \end{CD}
 $$
The first row is an exact sequence. We will now inductively
construct the operators $\nabla_t^*$, $\square_1^*$ etc, so that
all the rows are exact and the diagram commutes.

Let us choose a connection on the base manifold $M$ and on the
bundles $\pi$ and $\nu$. Then we can canonically embed the symbol
bundle into the bundle of differential operators (\cite{P,L$_3$})
via some maps $q_i$, whence the splitting of the middle row of the
diagram into direct sum of the first and the third rows.

Now elements of the middle complex are 2-vectors, with the
upper and lower components being elements of the first and the
third complex respectively, so that the maps $\nabla_t^*$ and
$\Delta_i^*$ must have the following matrix form
 $$
\nabla_t^*=\begin{pmatrix} \varepsilon & 0 \\ X_t & \nabla_{t-1}^*
\end{pmatrix},\ \
\Delta_i^*=\begin{pmatrix} \vp^* & 0 \\ \Gamma_i & \Delta_{i-1}^*
\end{pmatrix},
 $$
where the map $\Gamma_i:S^{i-k}T\ot\nu^*\to\op{Diff}_{i-1}(\pi,\1)$
is given by the system and the splittings and the map
$X_t:\La^{m+1}_{[t]}ST\ot\nu^*\to\op{Diff}_{i-k-1}(\nu,\1)$ is
unknown.

By inductive assumptions, the equality
$\Delta_i^*\circ\nabla_t^*=0$ is equivalent to the following:
 $$
\Gamma_i\varepsilon+\Delta_{i-1}^*X_t=0.
 $$
Thus existence of $X_t$ is equivalent to the fact that
$\Gamma_t\varepsilon$ belongs to the image of $\Delta^*_{i-1}$. By
the inductive assumption this in turn is equivalent to the
equality $\a_{i-1}\Gamma_i\varepsilon=0$.

The map $\Gamma_i$ is given by the condition
$j_1\Gamma_i=\Delta_i^*q_2-q_1\vp^*$. Thus we get the required
equality:
 $$
(\pi_{i,i-1}^\E)^*\a_{i-1}\Gamma_i\varepsilon=\a_ij_1\Gamma_i\varepsilon=
\a_i\Delta_i^*q_2\varepsilon-\a_iq_1\vp^*\varepsilon=0
 $$
because the map $(\pi_{i,i-1}^\E)^*$ is a monomorphism and the the
second row is exact on the level of the right three arrows.
Therefore we constructed the map $\nabla_t^*$ in the second row,
so that the constructed part (right four arrows) is a complex in
the commutative diagram.

Now by the same procedure we construct successively the arrows
$\square_1^*,\square_2^*$ etc with the condition that the second
row is a complex and the diagram commutes.

Since the first and the third complexes are exact, the diagram
chase yields exactness of the middle complex. This proves the
theorem in the linear case.

In the non-linear case the proof is basically the same, but we need
to work over the infinitely prolonged equation $\E^{(\infty)}\subset
J^\infty(\pi)$, use linearizations of differential operators and
change the rings as in \S\ref{S24}. Then we get the complex similar
to the above, evaluation of which at the corresponding jet equals
the complex for the case of the linearized system.
 \end{proof}

We will now turn the exact sequence from the theorem into
compatibility complex for the differential operator $\Delta$. For
the linear differential operators it is known \cite{S,KLV,V} that
the cohomology of the compatibility complex coincide with the
cohomology of the stable Spencer D-complex. We can generalize this
now to the non-linear case, whenever the system is of generalized
complete intersection type and is formally integrable.

 \begin{rk}\po\label{prababa}
In fact, the general non-linear Spencer D-complex is unknown, but in
the considered situation we can construct it via the splitting as in
the proof above. Indeed, the symbolic part is given by the diagram
from the following section \ref{S41} and then the complexes can be
constructed inductively in the filtration number by the technique of
the above proof.
 \end{rk}

\subsection{\hpss The compatibility complex}

 \abz
The homomorphisms of the spaces of differential operators as in
theorem \ref{exact_seq_DO} are not always coming from the actual
differential operators (as we treated in \S\ref{S12}), but in our
situation we can arrange this.

Indeed, in the proof from the preceding section we need to
construct the operators $\nabla_t^*,\square_1^*,\square_2^*$ etc
only when they appear for the first time (the corresponding index
is zero or equivalently, when the defining differential operators
are introduced). Afterwards, we can let
$\nabla_t=\nabla_{t-1}^{(1)}$ to be the prolongation,
$\square_{1,s}=\square_{1,s-1}^{(1)}$ and so on.

Thus we get actual differential operators between vector bundles,
forming {\em compatibility complex\/} for the differential
operator $\Delta$ (or the system $\E$).

Let us consider at first the case of linear system $\E$. Denote by
$\mathcal{S}_\E=\mathcal{S}ol_\E$ the sheaf of local solutions of
the PDE system $\E$. Then by the preceding discussion we
constructed the complex
 $$
0\to\mathcal{S}_\E\to
C^\infty_\text{loc}(\pi)\stackrel{\Delta}\longrightarrow
C^\infty_\text{loc}(\nu)\stackrel{\nabla}\longrightarrow
C^\infty_\text{loc}(\La^{m+1}\nu)\stackrel{\square_1}\longrightarrow
C^\infty_\text{loc}(\pi^*\ot\La^{m+2}\nu)\to\dots,
 $$
which is formally exact and thus is a compatibility complex of the
PDE system $\E$ (in some cases this complex is required to be
constructed to both sides, but we will work only with right
differential resolutions).

 \begin{rk}\po
In the classical geometric theory of PDEs, the compatibility
resolution is constructed for any formally integrable system of
PDEs (see the section about Kuranishi theorem in \cite{S}), but
then the construction is quite sophisticated (as in the second
Spencer complex \cite{S,KLV}). In our case of generalized complete
intersection the spaces in the complex are explicit and simple.

In fact, the same construction will work for any formally
integrable system, where we have managed to write explicitly a
resolution of the symbolic module.
 \end{rk}

What is more interesting is that we can write the whole
compatibility complex explicitly. Consider at first two examples.

{\bf 1.} Let us write exact form of the compatibility conditions
for linear non-homogeneous scalar PDEs of complete intersection
type ($m=1$ and $1<r\le n$):
 \begin{equation}\label{nl-sys}
\Delta_1(u)=f_1,\dots,\Delta_r(u)=f_r.
 \end{equation}
We solve non-homogeneous linear system under assumption that the
linear homogeneous system is compatible. By theorem \ref{th A}
this means certain commutation relation between operators
$\Delta_i$ and their differential corollaries (we use the standard
summation rule by the repeated indices):
 \begin{gather}
[\Delta_i,\Delta_j]=C_{ij}^k\Delta_k,\label{gath1}\\
\!\!\!\sum_{\text{cyclic}:\ ijk}
[\Delta_i,[\Delta_j,\Delta_k]]=0\Leftrightarrow
\sum_{\text{cyclic}:\ ijk} [\Delta_i,C_{jk}^\a]=
\!\!\!\sum_{\text{cyclic}:\ ijk}\!\!\!C_{jk}^\b C_{\b
i}^\a+D_{ijk}^{\b\g}R_{\b\g}^\a,\nonumber\\
\text{ }\hspace{70pt}\text{ where }
R_{\b\g}^\a=\Delta_\b\d_\g^\a-\Delta_\g\d_\b^\a- C_{\b\g}^\a,
\qquad\text{etc}.\label{gath2}
 \end{gather}
Here $C_{ij}^k$ are certain scalar differential operators on the
base manifold $M^n$ of $\op{ord}C_{ij}^k<l_i+l_j-l_k$,
$l_s=\op{ord}\Delta_s$, and $D_{ijk}^{\b\g}$ are scalar
differential operators of $\op{ord}D_{ijk}^{\b\g}
<l_i+l_j+l_k-\min[l_s+l_t:s,t\in\{\a,\b,\g\},s\ne t]-1$ etc.

We may assume that $C_{ij}^k$ is skew-symmetric by $ij$ and
$D_{ijk}^{\b\g}$ is skew-symmetric by $ijk$ and by $\b\g$. Similar
conditions will be imposed on other differential operators with
multi-indices, which occur in the higher Jacobi identities.

Because of (\ref{gath1}) the compatibility conditions for
non-linear system (\ref{nl-sys}) are
 $$
[\Delta_i-f_i,\Delta_j-f_j]-C_{ij}^k(\Delta_k-f_k)=
\Delta_jf_i-\Delta_if_j+C_{ij}^kf_k=0.
 $$
Thus the differential 1-syzygy operator is
$\nabla=(\nabla_{ij})_{1\le i<j\le n}$, where
$\nabla_{ij}=\Delta_iP_j-\Delta_jP_i-C_{ij}^kP_k$ with $P_k$ being
the projector to the $k$-th component. In other words, for
$f=(f_1,\dots,f_r)^t$ we have:
 $$
\nabla(f_1,\dots,f_r)=(\dots,\Delta_if_j-\Delta_jf_i-C_{ij}^kf_k,\dots).
 $$

To obtain differential 2-syzygy we write compatibility conditions
for the system $\nabla_{ij}(f)=\l_{ij}$, which are:
 \begin{multline*}
\sum_{\text{cyclic}:\ ijk}(\Delta_i\nabla_{jk}
-C_{ij}^\a\nabla_{\a
k}+D_{ijk}^{\a\b}\nabla_{\a\b})\\
=\sum_{\text{cyclic}:\ ijk} (C_{ij}^\b C_{\b
k}^\a-[\Delta_k,C_{ij}^\a]+D_{ijk}^{\b\g}R_{\b\g}^\a)P_\a=0
 \end{multline*}
due to (\ref{gath2}). Thus the next operator in the differential
resolution is $\square=(\square_{ijk})_{1\le i<j<k\le n}$, where
for $\l=(\dots,\l_{ij},\dots)$ we have
 $$
\square_{ijk}(\l)=\sum_{\text{cyclic}:\
ijk}(\Delta_i\l_{jk}-C_{ij}^\a\l_{\a k}+D_{ijk}^{\a\b}\l_{\a\b}).
 $$
Continuing in this way we arrive at the following compatibility
complex:
 \begin{equation}\label{deRham}
0\to\mathcal{S}_\E\longrightarrow
C^\infty(\1)\stackrel{\Delta}\longrightarrow
C^\infty(\nu)\stackrel{\nabla}\longrightarrow
C^\infty(\La^2\nu)\stackrel{\square}\longrightarrow
C^\infty(\La^3\nu)\to\dots,
 \end{equation}
where $\nu$ is a trivial bundle over $M$ with the fiber $\R^r$
(naturally $\La^0\nu\simeq\1$).

The operator $\Theta:C^\infty(\La^k\nu)\to C^\infty(\La^{k+1}\nu)$
in the above complex has the following form. Let
$\oo=\{i_1<\dots<i_k\}$, $\varsigma=\{j_1<\dots<j_{k+1}\}$ be
ordered multi-indices, so that the elements of the domain for
$\Theta$ are $\theta=\{\theta_\oo\}$, while the elements for the
range are $\Theta=\{\Theta_\varsigma\}$. The map is
 \begin{multline*}
\Theta_\varsigma(\theta)=\frac1{k!}\sum_{\z\in
S_{k+1}}(-1)^\z\Bigl[ \Delta_{\z(j_1)}\theta_{\z(j_2\dots
j_{k+1})}+
\varkappa_{1,k}C_{\z(j_1j_2)}^\a\theta_{\a,\z(j_3\dots j_{k+1})}\\
+\varkappa_{2,k}D_{\z(j_1j_2j_3)}^{\a_1\a_2}\theta_{\a_1\a_2\z(j_4\dots
j_{k+1})}+\dots+\varkappa_{k,k}K_{\z(j_1\dots
j_{k+1})}^{\a_1\dots\a_k}\theta_{\a_1\dots\a_k}\Bigr].
 \end{multline*}
($D,\dots,K$ are some differential operators with multi-indices
mentioned above and $\varkappa_{i,k}$ some real numbers calculated
recursively, $1\le i\le k$; for instance the first numbers are
$\varkappa_{1,1}=-\frac12,\varkappa_{1,2}=-1,
\varkappa_{2,2}=1,\varkappa_{1,3}=-\frac32$ and so on).

Notice that on the symbolic level (\ref{deRham}) is the usual de
Rham complex, but due to non-commutativity of the components
$\Delta_i$ of the operator $\Delta$, defining the system $\E$, we
need to compensate the syzygy maps with the tails of (lower order)
differential operators. Horizontal de Rham complex (\cite{KLV}) is a
particular case of (\ref{deRham}), when the operators $\Delta_i$
commute.

{\bf 2.} Consider now a vector system of $(m+1)$ equations on $m$
unknown functions, given by the operators
$\Delta_i=(\Delta_{i1},\dots,\Delta_{im})$, which is of
generalized complete intersection type:
 \begin{equation*}
\left\{\begin{array}{l}
\Delta_{11}(u_1)+\dots+\Delta_{1m}(u_m)=f_1\\
\qquad\dots\quad\dots\quad\dots\quad\dots\quad\dots\quad\dots\quad \\
\Delta_{(m+1)1}(u_1)+\dots+\Delta_{(m+1)m}(u_m)=f_{m+1}\ .
\end{array}\right.
 \end{equation*}

In this case the compatibility condition for the homogeneous
system is given by the condition from Theorem \ref{th A}
 $$
\{\Delta_1,\dots,\Delta_{m+1}\}=B^i\Delta_i,
 $$
where $\{\dots\}$ is the multi-bracket from \S\ref{S31b} and $B^i$
some scalar differential operators on $M$ of
$\op{ord}(B^i)<\sum_{j\ne i}\op{ord}(\Delta_j)$.

Assuming formal integrability for the homogeneous system, we get
the compatibility complex:
 $$
0\to\mathcal{S}_\E\hookrightarrow C^\infty(M,\R^m)
\stackrel{\Delta}\longrightarrow C^\infty(M,\R^{m+1})
\stackrel{\nabla}\longrightarrow C^\infty(M,\R)\to0,
 $$
where $\Delta=(\Delta_{ij})$ is the $(m+1)\times m$ matrix
determining the system $\E$ and
 $$
\nabla(f_1,\dots,f_{m+1})^\text{t}=
\sum_{k=1}^{m+1}\left((-1)^{k-1} \op{Ndet}[\Delta_{ij}]_{i\ne
k}^{1\le j\le m}\,-B^k\right)f_k,
 $$

For $m=2$ this formula looks as follows (the determinants are
non-commutative):
 \begin{multline*}
\nabla\left(\begin{array}{c}f_1\\ f_2\\ f_3\end{array}\right)=
\left(\begin{vmatrix}\Delta_{21} & \!\!\!\Delta_{22}\\ \Delta_{31}
& \!\!\!\Delta_{32}\end{vmatrix}-B^1\right)f_1\\+
\left(\begin{vmatrix}\Delta_{31} & \!\!\!\Delta_{32}\\ \Delta_{11}
& \!\!\!\Delta_{12}\end{vmatrix}-B^2\right)f_2+
\left(\begin{vmatrix}\Delta_{11} & \!\!\!\Delta_{12}\\ \Delta_{21}
& \!\!\!\Delta_{22}\end{vmatrix}-B^3\right)f_3.
 \end{multline*}

{\bf 3.} Now we combine the first case of arbitrary dimension of
base with the second case of arbitrary dimension of the fiber of
$\pi$ as in \S\ref{S33}. Then we can write the compatibility
complex explicitly (in the generalized complete intersection case;
in general it is more complicated to write the
$W_H(\E)$-component).

Let $\pi$ and $\nu$ be trivial vector bundles (locally this is not a
restriction). Moreover we will assume $\nu$ is graded in the various
orders case for the operator $\Delta:C^\infty(\pi)\to
C^\infty(\nu)$.

 \begin{theorem}\po
The compatibility complex for a collection of linear operators
$\Delta=\{\Delta_i\}$, defining a system $\E$ of generalized
complete intersection type, has the following form:
 \begin{multline*}
0\to\mathcal{S}_\E\to
C^\infty_\text{loc}(\pi)\stackrel{\Delta}\longrightarrow
C^\infty_\text{loc}(\nu)\stackrel{\nabla}\longrightarrow
C^\infty_\text{loc}(\La^{m+1}\nu)\stackrel{\square_1}
\longrightarrow\\
\to C^\infty_\text{loc}(\pi^*\ot\La^{m+2}\nu)
\stackrel{\square_2}\longrightarrow
C^\infty_\text{loc}(S^2\pi^*\ot\La^{m+3}\nu)\to\dots
 \end{multline*}
The differential 1-syzygy operator $\nabla:C^\infty(M;\R^r)\to
C^\infty(M;\La^{m+1}\R^r)$ is given by the formula
 $$
\nabla_\tau(f)=\sum_{k=1}^{m+1}(-1)^{k-1}
\op{Ndet}[\Delta_{i_sj}]_{s\ne k}^{1\le j\le
m}\,f_{i_k}-\sum_{j=1}^rB^j_\tau f_j,
 $$
where $f=(f_1,\dots,f_r)$ and $\nabla=\{\nabla_\tau\}$, $\tau$
being a multi-index $1\le i_1<\dots<i_{m+1}\le r$. The
coefficients $B_\tau^k$ are obtained from the multi-bracket
relation (compatibility of the linear homogeneous system)
 $$
\{\Delta_{i_1},\dots,\Delta_{i_{m+1}}\}=B^j_{i_1\dots
i_{m+1}}\Delta_j.
 $$
The higher syzygy operators are given by the method above.\qed
 \end{theorem}

Non-linear case is treated by the same formulas via the
linearization operator. We need to restrict to infinitely
prolonged equation $\E^{(\infty)}\subset J^\infty(\pi)$. Denote
the projection from this equation to the base $M$ by $\pi_\E$.
Then changing the operator $\Delta$ to its linearization
$\ell_\Delta$ (the operator itself on the equation vanishes) and
coupling this to the arguments above we will arrive at the
non-linear compatibility complex
 \begin{multline*}
0\to\mathcal{S}ym_\E\to
C^\infty_\text{loc}\bigl(\pi_\E^*(\pi)\bigr)\stackrel{\ell_\Delta}\longrightarrow
C^\infty_\text{loc}\bigl(\pi_\E^*(\nu)\bigr)\stackrel{\hat\nabla}\longrightarrow
C^\infty_\text{loc}\bigl(\La^{m+1}\pi_\E^*(\nu)\bigr)\stackrel{\hat\square_1}
\longrightarrow\\
\to
C^\infty_\text{loc}\bigl(\pi_\E^*(\pi^*)\ot\La^{m+2}\pi_\E^*(\nu)\bigr)
\stackrel{\hat\square_2}\longrightarrow
C^\infty_\text{loc}\bigl(S^2\pi_\E^*(\pi^*)\ot\La^{m+3}\pi_\E^*(\nu)\bigr)\to\dots
 \end{multline*}
(note the change of the solution sheaf
$\mathcal{S}_\E=\mathcal{S}ol_\E$ in the linear case to the
symmetries sheaf $\mathcal{S}ym_\E=\op{Ker}(\ell_\Delta)$ in the
non-linear case: For a linear equation $\E$ each shift by a
solution is a symmetry).

This complex is formally exact and its cohomology can be identified,
as we have noticed in the previous section, with what we call
non-linear Spencer cohomology.

This makes an effective representation for this important
invariant $H^i_D(\E)$ of the system $\E$ of differential equations
(recall that in general non-linear case this invariant is not
defined).

\section{\hps Integrability theory}

 \abz
In this section we consider certain topics closely related to the
main subject of this paper. In fact, explicit
compatibility/solvability criteria can be applied to integrability
of determined PDEs. Indeed, one can find (simple Frobenius or
sometimes more sophisticated) compatibility schemes in most
well-known integrability approaches, in particular in the following
theories:

 \begin{itemize}
 \item[$-$] Symmetry calculus;\quad B\"acklund transformations;
 \item[$-$] Lax pairs; \quad Zero-curvature representations;
 \item[$-$] Darboux integrability; \quad Sato theory etc.
 \end{itemize}

Thus it is important to understand formal integrability and
related Poisson geometry, method of differential constraints via
multi-brackets or other criteria. Some of these topics will be
considered in this section.

\subsection{\hpss Spencer cohomology and curvature tensors}\label{S41}

 \abz
Let $\E\subset J^\infty(\pi)$ be a system of differential
equations defined by a set of differential operators
$\Delta:C^\infty(\pi)\to C^\infty(\nu)$, which we allow to be of
different orders.

\vspace{4pt}
 \begin{proof}[Proof of theorem \ref{th B}]
Below we denote by $g=\oplus g_i$ the symbolic system associated
to $\E$.

Consider the following commutative diagram:
 $$
 \begin{array}{cccccccccc}
 &\hspace{-8pt}&0&&0&&0&&0\\
 &\hspace{-8pt}&\downarrow&&\downarrow&&\downarrow&&\downarrow\\
\hspace{-3pt}0&\hspace{-8pt}\to\hspace{-7pt}& g &\hspace{-8pt}\to
&\hspace{-8pt} ST^*\ot\pi &\hspace{-8pt}\to &\hspace{-8pt}
ST^*\ot\nu&\hspace{-8pt}\to &\hspace{-8pt}
\La^{m+1}_\diamond &\hspace{-8pt}\to\dots\\
 &\hspace{-8pt}&\downarrow&&\downarrow&&\downarrow&&\downarrow\\
\hspace{-3pt}0&\hspace{-8pt}\to\hspace{-7pt}& g\ot T^*
&\hspace{-8pt}\to &\hspace{-8pt} ST^*\ot\pi\ot T^*
&\hspace{-8pt}\to &\hspace{-8pt} ST^*\ot\nu\ot
T^*&\hspace{-8pt}\to &\hspace{-8pt}
\La^{m+1}_\diamond\ot T^* &\hspace{-8pt}\to\dots\\
 &\hspace{-8pt}&\downarrow&&\downarrow&&\downarrow&&\downarrow\\
\hspace{-3pt}0&\hspace{-8pt}\to\hspace{-7pt}& g\ot\La^2T^*
&\hspace{-8pt}\to &\hspace{-8pt} ST^*\ot\pi\ot\La^2T^*
&\hspace{-8pt}\to &\hspace{-8pt}
ST^*\ot\nu\ot\La^2T^*&\hspace{-8pt}\to &\hspace{-8pt}
\La^{m+1}_\diamond\ot\La^2T^* &\hspace{-8pt}\to\dots\\
 &\hspace{-8pt}&\vdots&&\vdots&&\vdots&&\vdots&\hspace{-20pt}\ddots\\
 &\hspace{-8pt}&\downarrow&&\downarrow&&\downarrow&&\downarrow\\
\hspace{-3pt}0&\hspace{-8pt}\to\hspace{-7pt}& g\ot\La^nT^*
&\hspace{-8pt}\to &\hspace{-8pt} ST^*\ot\pi\ot\La^nT^*
&\hspace{-8pt}\to &\hspace{-8pt}
ST^*\ot\nu\ot\La^nT^*&\hspace{-8pt}\to &\hspace{-8pt}
\La^{m+1}_\diamond\ot\La^nT^* &\hspace{-8pt}\to\dots\\
 &\hspace{-8pt}&\downarrow&&\downarrow&&\downarrow&&\downarrow\\
 &\hspace{-8pt}&0&&0&&0&&0
 \end{array}
 $$
where we write $\La^{m+1}_\diamond=ST^*\ot\La^{m+1}\nu$ for
brevity. The rows are generalized Spencer complexes (i.e.
$\R$-dual to Buchsbaum-Rim) and are exact.

The columns are the usual Spencer $\d$-complexes and all except
the first one are acyclic, i.e. have cohomology in the bi-degree
$(0,0)$ only. Thus by the diagram chase the cohomology of the
first column at the $j$-th term $H^{*,j}(g)$ equals the zero
cohomology group of the $(j+1)$-th column.
 \end{proof}

Denote $h^{i,j}=\dim H^{i,j}(g)$ the Betti numbers for Spencer
cohomology.

 \begin{cor}\po\label{cor1}
For a generalized complete intersection $g$ we have: $h^{*,0}=m$,
$h^{*,1}=r$, $h^{*,j}=\binom{m+j-3}{j-2}\binom{r}{m+j-1}$ for
$2\le j\le r+1-m$ and $h^{*,j}=0$ for $r+1-m<j\le n$. \qed
 \end{cor}

Note that the Euler characteristic vanishes as it should:
$\chi=\sum(-1)^ih^{*,i}=0$.

It is possible to specify the bi-grades, where the $\d$-cohomology
does not vanish. If $g$ is a system of pure order $k$, then the
non-zero Betti numbers are: $h^{0,0}$, $h^{k-1,1}$ and
$h^{km+kj-j-k,j}$ for $2\le j\le r$.

For a system of different orders $l(1),\dots,l(r)$ the above
complex allows to specify the bi-grades too. The precise
combinatorics is straightforward from the weighted Buchsbaum-Rim
complex. However, since the formulas are involved, we indicate
only what happens to first and second $\d$-cohomology (this latter
space contains curvature, i.e. Weyl tensors, see Corollary
\ref{cor2}).

For them the only non-zero Betti numbers are: $h^{l(i)-1,1}$,
$1\le i\le r$ and $h^{l(i_1)+\dots+l(i_{m+1})-2,2}$, $1\le
i_1<\dots<i_{m+1}\le r$.

 \begin{rk}\po
Formulas for $h^{*,j}$ from corollary \ref{cor1} suggest some
polynomial relations between multi-brackets. In fact, for $m=1$ we
can see the basis of $H^{*,j}$ for $j>2$ as the set of power $j$
subsets in $\{1,\dots,r\}$: for each such subset we associate all
possible iterated Mayer-Jacobi brackets, which is an analog of
$j$-form. The Jacobi identity and its higher analogs yield
relations between these iterated brackets
$\{...\{\{F_{i_1},F_{i_2}\},F_{i_3}\},\dots,F_{i_j}\}$.

When $m>1$ the number of subsets of power $m+j-1$ in
$\{1,\dots,r\}$ is $\binom{r}{m+j-1}$. Another factor of $\dim
H^{*,j}$ for $j\ge2$ is related to the fact that iteration of
brackets is now arranged in a multiplicative manner, so that the
relations are generalized Plucker identities as in \S\ref{S31c}.
 \end{rk}

\subsection{\hpss Integrability of characteristics and 
multi-brackets}\label{S42}

 \abz
Consider the space of scalar linear differential operators
$\op{Diff}(\1,\1)$ on the manifold $M$. It bears the structure of
infinite-dimensional Lie algebra when equipped with the Jacobi
bracket $\{,\}$ (this is also true for the space
$\op{Diff}(\pi,\pi)$). This bracket induces the classical Poisson
bracket on $T^*M$ as follows.

Let $\z:\op{Diff}(\1,\1)\to C^\infty(T^*M)$ be the symbol map. It
associates to an order $k$ differential operator $F$ a polynomial
$\z(F)=\op{smbl}_k(F)$ of degree $k$ in momenta $p$ in canonical
coordinates $(x,p)$ on $T^*M$. Then we have (with an ambiguity in
notations):
 $$
\{\z(F),\z(G)\}=\z(\{F,G\}),
 $$
where the brackets on the left are Poisson and on the right are
Jacobi brackets. In order to make a distinction we will write
$\{,\}_\z$ for the symbolic (2-)bracket.

Consider a system of linear scalar PDEs $\E=\{F_1=0,\dots,F_r=0\}$,
$F_i\in\op{Diff}(\1,\1)$. Let $\I=\I(F)$ be the differential ideal
generated by $F_i$, $1\le i\le r$.

 \begin{prop}\po\label{aua}
Let the scalar system $\E$ be formally integrable. Then the
corresponding ideal $\I$ is closed with respect to the higher
Jacobi bracket.
 \end{prop}

 \begin{proof}
The system $\E$ is defined by the ideal $\I$ whenever formally
integrable. Jacobi bracket of two operators is their differential
corollary. Since the ideal $\I$ is differentially closed the claim
follows.
 \end{proof}

 \begin{cor}\po
The corresponding symbolic (characteristic) ideal $I=\z(\I)$ is
closed with respect to Poisson bracket. \qed
 \end{cor}

A stronger statement was proved in \cite{GQS}: Namely the radical
$\sqrt{I}$ is Poisson-closed too (in fact, the claim was justified
only for the components of $\op{Char}^\C(\E)$ of maximal
dimension). This is the celebrated integrability of
characteristics: The affine characteristic variety
$\op{Char}^\C_\text{aff}(\E)\subset T^*M$ is integrable in the
Frobenius sense, i.e. if $F$ and $G$ vanish on it, their Poisson
bracket $\{F,G\}$ vanishes on it as well.

This was applied in \cite{GQS} to the case when $\I$ (resp. $I$)
is the annihilator of the module $\E^*$ defining a formally
integrable PDE system (resp. the symbolic module $g^*$). However,
for the system involving several unknown functions the
characteristic variety does not bear the complete information
about the dynamics in generalization of Hamilton-Jacobi theory.
Indeed, as we described in \S\ref{S13} the characteristic variety
$\op{Char}^\C(\E)$ of the system $\E$ defined by a differential
operator $\Delta:C^\infty(\pi)\to C^\infty(\nu)$ is the support of
the characteristic sheaf $\mathcal{K}=\op{Ker}\z(\Delta)$.

Dualizing this we can equally work with $\op{Coker}\z(\Delta^*)$. We
describe at first the picture with differential operators. The
$\op{Diff}(\1,\1)$-module $\E^*$ is given by the exact sequence
 $$
\op{Diff}(\nu,\1)\stackrel{\vp^\Delta}\longrightarrow\op{Diff}(\pi,\1)\to\E^*\to0.
 $$
Denote $\mathcal{J}=\op{Im}(\vp^\Delta)$ the submodule in
$\op{Diff}(\pi,\1)$.

The following statement is obtained similarly to proposition
\ref{aua}:

 \begin{prop}\po\label{uau}
Let the vector system $\E$ be formally integrable. Then the
corresponding module $\mathcal{J}$ is closed with respect to the
multi-bracket $\{\cdot,\dots,\cdot\}$.\qed
 \end{prop}

This is the necessary condition for formal integrability from
theorem \ref{th A}. It is sufficient only in special cases, see
\S\ref{S44} below. We now define the symbolic multi-bracket
generalizing the above Poisson 2-bracket.

 \begin{dfn}\po\label{smbl_mltbr}
The symbolic multi-bracket of sections $f_i=\z(F_i)$ is given by
the formula:
 $$
\{\z(F_1),\dots,\z(F_{m+1})\}_\z=\z(\{F_1,\dots,F_{m+1}\})
 $$
(we use the index $\z$ just to keep distinction between
multi-brackets).
 \end{dfn}

It is easy to see that the right-hand side does not depend on the
lower order terms of $F_i$, so that the left-hand side is the
well-defined expression $\{f_1,\dots,f_{m+1}\}_\z$. However the
above formula gives the symbolic multi-bracket only for
vector-valued polynomials $f\in ST\ot\pi^*$. The standard trick
extends this to formal series and analytic functions, but we can
define the symbolic multi-bracket $\{\cdots\}_\z$ to all smooth
vector-valued functions on the cotangent bundle $f_i\in
C^\infty(T^*M;\pi^*)$, $1\le i\le m+1$. This follows from the
following proposition.

Let $e_i$ be a basis in the bundle $\pi^*$ (which we assume
trivial or make a localization). Then a vector-valued function
$f\in C^\infty(T^*M;\pi^*)$ can be identified with the collection
of functions $f^j\in C^\infty(T^*M)$ via the decomposition
$f=f^je_j$. We use the components $f^j$ below.

 \begin{prop}\po
The symbolic multi-bracket of the vector-valued functions $f_i$
can be expressed via the product and the standard Poisson bracket
of their components $f_i^j$ ($1\le i\le m+1$, $1\le j\le m$).
 \end{prop}

 \begin{proof}
The multi-brackets of differential operators
$F_i=(F_i^1,\dots,F_i^m)$ has the following form in components:
 $$
\{F_1,\dots,F_{m+1}\}^l=\sum_{\tau\in S_{m+1}}(-1)^\tau
F^1_{\tau(1)}F^2_{\tau(2)}\cdots F^m_{\tau(m)}F^l_{\tau(m+1)}.
 $$
Taking the symbol and making elementary transformations we get
(both indices $j$ and $s$ vary between 1 and $m$) the
multi-brackets of $f_i=(f_i^1,\dots,f_i^m)$:
 \begin{multline}\label{mult-via-2}
\{f_1,\dots,f_{m+1}\}_\z^l=(-1)^{m-l}\sum_{\tau\in
S_{m+1}}\Bigl[\frac12\prod_{s\ne l}f^s_{\tau(s)}\cdot
\{f^l_{\tau(l)},f^l_{\tau(m+1)}\}_\z \\
+\sum_{j>l}\prod_{s\ne j,l}f^s_{\tau(s)}\cdot
\{f^l_{\tau(l)},f^j_{\tau(j)}\}_\z\cdot f^l_{\tau(m+1)}\Bigr]
 \end{multline}
This represents components of the symbolic multi-bracket via the
2-bracket.\!\!
 \end{proof}

Remark that our symbolic multi-bracket differs from other
multi-versions of the Poisson bracket, like Fillipov-Nambu or
generalized Poisson bracket (\cite{APB,MV,N}). For instance, the
multi-bracket with $m>1$ is not a derivation in its arguments.
However it is not pretty far from this:

 \begin{cor}\po
The symbolic multi-bracket $\{f_1,\dots,f_{m+1}\}_\z$ is a
differential operator of the first order in each of its arguments
$f_i\in C^\infty(T^*M,\pi^*)$. \qed
 \end{cor}

Indeed, we can write for $g\in C^\infty(T^*M)$
 \begin{multline}\label{diff_1-for-mult}
\{f_1,\dots,gf_i,\dots,,f_{m+1}\}_\z-g\{f_1,\dots,f_i,\dots,,f_{m+1}\}_\z\\
=\sum\nolimits_{j,k}\{g,f^j_k\}_\z\cdot
c_{ij}^k(f_1,\dots,f_i,\dots,,f_{m+1}),
 \end{multline}
where the exact form of the functions $c_{ij}^k$ can be obtained
from (\ref{mult-via-2}).

Let $J=\z(\mathcal{J})$ be the symbol of the submodule
$\mathcal{J}$. Denoting as in \S\ref{S22}-\ref{S23}
$\vp=\z_\Delta$ the symbol of the differential operator $\Delta$,
we can define $J$ via the exact sequence
 $$
ST\ot\nu^*\stackrel{\vp}\longrightarrow ST\ot\pi^*\to g^*\to0
 $$
as $J=\op{Im}(\vp)\subset ST\ot\pi^*$, so that the symbolic module
is $g^*=ST\ot\pi^*/J$.

Then proposition \ref{uau} yields:

 \begin{cor}\po
The (characteristic) submodule $J=\z(\mathcal{J})$ is closed with
respect to the symbolic multi-bracket. \qed
 \end{cor}

Beside identity (\ref{diff_1-for-mult}) the symbolic multi-bracket
satisfies the same properties as the multi-bracket of differential
operators, for instance we have the generalized Pl\"ucker identity
as a corollary of \S\ref{S31c} (the upper index means
component):
 $$
\sum(-1)^k(\{f_1,\dots,\check f_k,\dots,f_{m+2}\}_\z^\dag)^i\cdot
f_k= \sum(-1)^k \{f_1,\dots,\check f_k,\dots,f_{m+2}\}_\z\cdot
f_k^i,
 $$
where $\{\cdots\}_\z^\dag$ is the symbolic multi-bracket
associated to the opposite multi-bracket of differential operators
$\{\cdots\}^\dag$ of \S\ref{S31b} by the same rule as in
definition~\ref{smbl_mltbr}.

Consider again the case (of importance to PDEs), when $f_i\in
ST\ot\pi^*$ are polynomial vector-valued functions on $T^*M$.

Let $I_0(J)\subset R=ST$ be the ideal generated by the functions
$f_i^j$ and $C_0$ its Poisson center. The following statement
follows directly from (\ref{diff_1-for-mult}):

 \begin{cor}\po
The multi-bracket is a homomorphism in its arguments over the
ideal $C_0\subset ST$:
$\{\cdots\}_\z\in\op{Hom}_{C_0}(\La^{m+1}J,J)$. \qed
 \end{cor}

The symbolic multi-bracket $\{\dots\}_\z$ can be used to formulate
a Hamiltonian formalism as generalized Poisson brackets are used.
Namely, let $f_1,\dots,f_m\in J$ be Hamiltonians and $h\in
ST\ot\a^*$ a polynomial vector-function. A multi-Hamiltonian
operator
 $$
X_{f_1\we\dots\we f_m}:g^*\to g^*,\qquad h\,\op{mod}
J\mapsto\{h,f_1,\dots,f_m\}_\z\,\op{mod} J
 $$
determines transport on $\op{Char}^\C(g)$ (symbolic module $g^*$
is supported on characteristics). Caustics of solutions (wave
fronts) develop according to it. More details will be provided
elsewhere.

\subsection{\hpss Applications to smooth integrability of PDEs}\label{S43}

 \abz
The classical Lagrange-Charpit method \cite{Gou,Gu} is designed
for first order scalar PDEs and is as follows. Let
$F=F_1(x^1,\dots,x^n,u,p_1,\dots,p_n)=0$ be a differential
equation. To solve it one searches for functions $F_2,\dots,F_n$
on $J^1(\R^n)$ such that $[F_i,F_j]=0\mod(F_1,\dots,F_n)$, where
the $[,]$ is the classical Mayer bracket. Then the system
$F_1=0,\dots,F_n=0$ gives a finite dimensional family of solutions
of the PDE $F=0$ (this family is equivalent to an ODE).

If in addition $F_2,\dots,F_n$ are symmetries of the system,
$[F,F_i]=\l\cdot F$, then one can obtain a complete integral of
the PDE $F=0$ as the system $F_1=0,F_2=c_2,\dots,F_n=c_n$, which
in addition can be found in quadratures by the symmetry method of
S.Lie (\cite{Lie}).

 \begin{examp}\po
Consider the PDE $up_1\cdots p_n=x^1\cdots x^n$ on $\R^n$. It
possesses a collection of {\em auxiliary integrals\/}:
$F_i=\dfrac{p_iu^{1/n}}{x^i}$. This gives a complete integral of
the differential equation.
 \end{examp}

Now basing on our theorem \ref{th A} we can formulate a generalized
Lagrange-Charpit method. For manifolds of dimension 2 this was done
in \cite{KL$_2$}. To obtain the general version we start with the
following idea.

 \begin{dfn}\po
Let $\E$ be a formally integrable system of PDEs. Call a system
$\tilde\E$ an {\em auxiliary integral\/} (or a set of integrals)
for the system $\E$ if the joint system $\E\cap\tilde\E$ is also
formally integrable (= compatible).
 \end{dfn}

We proved in \cite{KL$_2$} that classical objects, such as point
symmetries, contact symmetries and intermediate integrals
(\cite{Gou,Gu,LE}) as well as higher symmetries (\cite{KLV}) are
partial cases of this notion. Moreover, some of the newly introduced
generalized symmetries are also auxiliary integrals.

Another traditional method for finding exact solutions of PDEs is a
method of differential constraints. Then one considers an
overdetermination $\ti\E$ on a system $\E$, such that the system
$\E\cap\tilde\E$ is {\em solvable\/}. The solvability is very
non-trivial to check in practice. Thus our notion of auxiliary
integral is more constructive, since we can use an effective
criterion from theorem \ref{th A} to check compatibility (see
\cite{KL$_1$,KL$_2$} for examples).

{\em Generalized Lagrange-Charpit method\/} is the following
special form of an auxiliary integral. Consider a determined
system of $m$ PDEs $\E=\{F_1=0,\dots,F_m=0\}$ on $m$ unknowns
functions $u_1,\dots,u_m$ (we can even start with an
underdetermined system). We search to add to it $n-1$ differential
equations $\ti\E=\{F_{m+1}=0,\dots,F_{m+n-1}=0\}$, so that the
resulting overdetermined systems
$\E\cap\tilde\E=\{F_1=0,\dots,F_{m+n-1}=0\}$ is:
 \begin{enumerate}
 \item[a)] generalized complete intersection;
 \item[b)] compatible.
 \end{enumerate}
Of course, one can add less functions, but advantage of $(n-1)$ is
that the system becomes of finite type (being compatible it
constitutes an integrable distribution by the Frobenius theorem
\cite{St,KL$_3$}) and thus reduces to a system of ODEs.

By theorem \ref{th A} compatibility of $\E\cap\tilde\E$ is given
by the conditions
 $$
\{F_{i_1},\dots,F_{i_{m+1}}\}=\sum\nolimits_j A_{i_1\dots
i_{m+1}}^j\circ F_j,
 $$
for some differential operators $A$ of orders
$\op{ord}(A^j_{i_1\dots i_{m+1}})\le i_1+\dots+i_{m+1}-i_j-1$.

On the symbolic level the first step is to include $f=\z(F)$ into
a submodule $J=\langle f_1,f_2,\dots,f_{m+n-1}\rangle$ of
$ST\ot\pi^*$, which is closed under the symbolic multi-bracket, as
integrability of characteristics from \S\ref{S42} claims:
 $$
\{f_{i_1},\dots,f_{i_{m+1}}\}_\z=\sum\nolimits_j a_{i_1\dots
i_{m+1}}^j\cdot f_j.
 $$
Then we shall adjust sub-principal symbols of $F_i$.

 \begin{examp}\po
Consider the Cauchy-Riemann system on the plane:
 \begin{equation}\label{CReq}
\E:\quad u_x=v_y,\ u_y=-v_x.
 \end{equation}
Then the differential equation
 $$
\ti\E:\quad 
\det\begin{bmatrix}u_x & v_x\\ u_y & v_y\end{bmatrix}=G(u,v)
 $$
is an auxiliary integral suitable for generalized Lagrange-Charpit
method iff
 \begin{equation}\label{detHes}
\Delta G=\|\nabla G\|^2/G,
 \end{equation}
where $\Delta$ is the standard Laplacian and $\nabla$ is the
standard gradient on the Euclidean plane $\R^2(u,v)$. Thus our
transformation can be seen as a kind of Backl\"und transformation,
which for any solution $G(u,v)$ of equation (\ref{detHes})
associates a 3-dimensional family of solutions of system
(\ref{CReq}) (for dimensional calculus see the next section).
 \end{examp}

\subsection{\hpss Formal dimension of the solutions space}\label{S43+e}

 \abz
Consider a symbolic system $g=\{g_l\subset S^lT^*\ot N\}$ and let
$V^*\subset T^*$ be a subspace. Then we can define another symbolic
system $\tg=\{g_l\cap S^lV^*\ot N\}\subset SV^*\ot N$. It is called
the {\it $V^*$-reduction\/} of $g$ (\cite{KL$_2$}). Denote
$W^*=T^*/V^*$.

We take $g$ to be the symbolic system of a generalized complete
intersection $\E$ of formal codimension $r$, the same number as in
definition \ref{dfn1}. Note that this definition can be
reformulated for symbolic systems as well. Then the characteristic
variety $\op{Char}^\C(g)=\op{Char}^\C(\E)$ has codimension $r-m+1$
in $P^\C T^*$.

 \begin{theorem}\po\label{fromKL2}
Let $g$ be a symbolic system of generalized complete intersection
type and the subspace $V^*\subset T^*$ of dimension $(r-m+1)$ be
transversal to the affine characteristic variety of $g$:
$\op{Char}(g)\cap P^\C V^*=\emptyset$. Then the reduced symbolic
system $\tg$ is also a generalized complete intersection, Spencer
$\d$-cohomology of the system $g$ and of its $V^*$-reduction $\tg$
are isomorphic and $g\simeq \tg\ot SW^*$.
 \end{theorem}

 \begin{proof}
This statement was proved for more general Cohen-Macaulay symbolic
systems in \cite{KL$_2$}. Since we proved in \S\ref{S23} that a
generalized complete intersection is Cohen-Macaulay, the claim
follows (the part that the $V^*$-reduction form a generalized
complete intersection is straightforward). The last claim of the
theorem $g_l\simeq \sum_j\tg_{l-j}\ot S^jW^*$ is not a part of
theorem A from \cite{KL$_2$}, but is contained in the proof, see
remark 8 loc.cit.
 \end{proof}

Notice that due to Noether normalization lemma \cite{E} a generic
subspace $V^*\subset T^*$ is transversal to the characteristic
variety over $\C$. The functions on $W=\op{Ann}(V^*)$ are those on
which a general solution of $\E$ depends, as is refereed to in the
discussion of functional dimension before Theorem \ref{th C}.

\vspace{4pt}
 \begin{proof}[Proof of theorem \ref{th C}]
Let us consider at first the case, when $r=n+m-1$. Then the system
$\E$ is of finite type, namely
$\pi_{i,i-1}:\E_i\stackrel\sim\to\E_{i-1}$ for $i\ge\sum_{i=1}^r
k_i$. The equation $\E_\infty$ is thus a finite-dimensional manifold
equipped with the Cartan distribution $\mathcal{C}_\E$
(\cite{KLV,KL$_1$}). Since the system is compatible, the local
solutions of $\E$ are integral manifolds of the distribution
$\mathcal{C}_\E$ of dimension $n$.

Thus the dimension of the solutions space is
 $$
\dim\mathcal{S}_\E=\dim\E_\infty-n=\sum_{i=0}^\infty\dim g_i
 $$
(the sum is indeed finite). We can calculate dimensions of the
symbol spaces $g_i$ explicitly, since the system is of generalized
complete intersection type.

For brevity sake we omit the immense combinatorics and provide
details of this step only for the case $n=\dim M=2$. The general
case is similar (we refer to \cite{KL$_5$}, where arbitrary systems
of pure order $k$ are scrutinized).

If $n=2$, we get $r=m+1$. Let $k_1\le\dots\le k_r$ be the orders
of the system $g$. We have:
 $$
\dim g_i=\left\{\begin{array}{ll}
 m(i+1) & \text{ if } i<k_1\\
 (m-j)(i+1)+\sum_{s=1}^jk_s & \text{ if } k_j\le i<k_{j+1}\\
 \sum_{s=1}^rk_s-1-i & \text{ if } k_r\le i<\sum_{s=1}^rk_s\\
 0 & \text{ else. }
\end{array}\right.
 $$
Thus
 \begin{multline*}
\sum\dim g_i=\sum_{i=0}^{k_1-1}m(i+1)+
\sum_{i=k_1}^{k_2-1}\Bigl((m-1)(i+1)+k_1\Bigr)\\
+\sum_{i=k_2}^{k_3-1}\Bigl((m-2)(i+1)+(k_1+k_2)\Bigr)+\dots
+\sum_{i=k_{r-1}}^{k_r-1}\bigl(k_1+\dots+k_{r-1}\bigr)\\
+\sum_{i=k_r}^{k_1+\dots+k_r-1}\bigl(k_1+\dots+k_r-1-i\bigr)=
\frac{(k_1+\dots+k_{r-1}-1)(k_1+\dots+k_{r-1})}2\\
+\sum_{i=1}^{r-1}\frac{k_i(k_i+1)}2+
\sum_{i=1}^r(k_i-k_{i-1})(k_1+\dots+k_{i-1})= \sum_{i<j}k_ik_j
 \end{multline*}
and the result follows.

Consider now the case $m<r<n+m-1$. Then the characteristic variety
$\op{Char}^\C(\E)$ is non-empty and the system $\E$ is of infinite
type. By theorem \ref{fromKL2} the symbolic system $g$ has a free
factor $SW^*$ and the dimension
 $$
p=\op{dim}W=n-\op{dim}V=n+m-r-1
 $$
is clearly the formal functional dimension of $\mathcal{S}ol_\E$.

Thus the quantity $\sum_{i=0}^t g_i$ when $t\to\infty$ grows as
$\ti d\cdot \dim S^tW^*$, where $\ti d$ is the formal functional
rank for the system $\tg$, which is of finite type (no complex
characteristics). The reduction $\tg$ is of generalized complete
intersection type and it has the same orders $k_1\le\dots\le k_r$
because the first Spencer $\d$-cohomology coincide (theorem
\ref{fromKL2}). Thus we can use the calculations above to conclude
 $$
\ti d=\sum\dim\tg_i=\sum_{i_1<\dots<i_l}k_{i_1}\cdots
k_{i_l}=d,\qquad l=\dim V=r-m+1.
 $$
Since $\dim W^*=n+m-r-1$, the  asymptotic of the Hilbert
polynomial for the symbolic module $g^*$ is
 $$
d\cdot\dim S^tW^*\sim d\frac{t^{n+m-r-2}}{(n+m-r-2)!}.
 $$
This proves the theorem.
 \end{proof}

 \begin{cor}\po
If a generalized complete intersection $\E$ of formal codimension
$r$ is formally integrable and has equations of the same order
$k$, then its formal functional dimension and functional rank are
$p=n+m-r-1$ and $d=\binom{n+m-1}{n}k^n$ (i.e. a general formal
solution depends on $d$ functions of $p$ arguments). \qed
 \end{cor}

\subsection{\hpss Generalizations of compatibility via brackets}\label{S44}

 \abz
Vanishing of multi-brackets is necessary, but not a sufficient
condition for compatibility. By theorem \ref{th A} it is
sufficient for generalized complete intersections, however this
does not generalize to more general class of Cohen-Macaulay
systems.

To see it consider a system of finite type. It is a Cohen-Macaulay
system. In fact, the finite type condition means
$\op{Char}^\C(g)=\emptyset$, so that $\dim g^*=0$. But this
condition also implies $\op{depth}g^*=0$ (alternatively
$0\le\op{depth}g^*\le\dim g^*$).

Now for a completely determined system of order $k$ (Frobenius
type), i.e.\ $g_k=0$, the compatibility conditions do not coincide
with these of our theorem. Actually, let us write the equation in
the orthonomic form:
 $$
\E=\Bigl\{\frac{\p^{|\z|}u^i}{\p
x^\z}=F_\z^i\Bigl(x,\frac{\p^{|\t|}u^j}{\p x^\t}\Bigr)\,\big|\,
1\le i,j\le m,\ |\z|=k,\ 0\le|\t|<k\Bigr\}.
 $$
The compatibility condition is
$\Phi_{a,b,\t}^i=\D_aF^i_{\t+1_b}-\D_bF^i_{\t+1_a}=0\mod\mathcal{J}_k(\E)$,
$|\t|=k-1$, while the multi-brackets are:
 $$
\{p^{i_1}_{\z_1}-F^{i_1}_{\z_1},\dots,
p^{i_{m+1}}_{\z_{m+1}}-F^{i_{m+1}}_{\z_{m+1}}\}=H_T+(\text{smaller
order terms}),
 $$
where the higher order term $H_T$ (with order equal to
$\sum_{j=1}^{m+1}|\z_j|-1$) is non-zero iff for some indices
$\a,\b$ we have (in the display below the braces mean a set):
 $$
i_\a=i_\b=k\in\{1,\dots,n\} \text{ and }
\{i_1,\dots,\check{i_\a},\dots,\check{i_\b},\dots,i_{m+1}\}=\{1,\dots,\check{k},\dots,n\},
 $$
in which case $H_T^i=\d_k^i(\pm\Pi_{j\ne\a,\b}\D_{\z_j})
[\D_{\z_{\a}}F^k_{\z_b}-\D_{\z_{\b}}F^k_{\z_a}]$. This provides
more conditions and they are of higher order.

 \begin{rk}\po
It is possible however to give explicit compatibility conditions for
some systems different from generalized complete intersection. Let
$\E$ be a Cohen-Macaulay system and $g$ its symbolic module.
Choosing a subspace $V^*\subset T^*$ not meeting the characteristic
variety $\op{Char}^\C(g)$ and of complimentary dimension, we get the
reduction $\tg$, which is also Cohen-Macaulay by theorem A of
\cite{KL$_2$}.

If we have a compatibility criterion for systems $\ti\E$ of type
$\tg$, we can transform it to get a criterion for the system $\E$
of type $g$ (cf. theorem \ref{fromKL2}). For instance, if the
system $\tg$ is completely determined (Frobenius type), we can use
the above formulas to get similar compatibility conditions for the
system $\E$ (which is not of Frobenius type!).
 \end{rk}

Thus usually for systems different from generalized complete
intersections, multi-brackets do not provide a basis of
compatibility conditions (though multi-brackets are part of them).
Indeed, in these other cases the obstructions to formal
integrability (Weyl tensors) belong to different Spencer
cohomology groups.

However in some cases the system not of generalized complete
intersection type can have compatibility conditions in a form of
multi-brackets. Usually this happens when the Spencer cohomology
is of the type described in theorem \ref{th B}. For instance, for
a (skew-)product of a generalized complete intersection and some
involutive system. We give two examples.

{\bf 1.} Let $J$ be an almost complex structure on a manifold $M$.
It defines the Cauchy-Riemann operator
$\bar\p_J:\Omega^{p,0}(M)\to\Omega^{p,1}(M)$, so that the system
$\E=\op{Ker}(\bar\p_J)$ is not a generalized complete
intersection. However we can represent the CR-operator locally as
the product of the operator
$\bar\p_J:C^\infty(M)\to\Omega^{p,1}(M)$ and the identity on
$\Omega^{p,0}(M)$. The compatibility condition (which is
equivalent to integrability of the structure $J$) has now the form
of vanishing multi-brackets.

{\bf 2.} Let $\nabla$ be a connection on the bundle $\pi:E\to M$
and $\Omega^{p,q}(E)$ the bundle of $p$-vertical, $q$-horizontal
forms. The horizontal de Rham operator
$d_\nabla:\Omega^{p,0}(E)\to\Omega^{p,1}(E)$ is locally a product
of $\Omega^{p,0}(E)$ and the horizontal de Rham operator
$d_\nabla:C^\infty(\pi)\to\Omega^{0,1}(E)$. The latter has
compatibility condition (which is equivalent to flatness of the
connection $\nabla$) again in the form of vanishing multi-brackets
as in theorem \ref{th A}.

\section{\hps Applications}

 \abz
In this section we apply the compatibility criterion to solve
certain problems arising in differential geometry. Many of them
address solvability of overdetermined systems of PDEs.

To decide if the system is solvable we add compatibility
conditions to the system (here we use theorem \ref{th A}),
investigate the new system, add its compatibility conditions etc.
In other words we apply prolongation-projection scheme and either
close up the system or get a contradiction (empty equation
$\E\subset J^k(\pi)$).

For example, linearization problem of 3-webs on the plane is
equivalent to solvability of a system of 2 scalar second order
equations of two variables. The system is a complete intersection
and the compatibility condition is given by vanishing of
Jacobi-Mayer bracket, as was sketched in \cite{KL$_1$}. This
condition is equivalent to vanishing of the Chern curvature and
yields parallelizable webs.

To linearize the web one adds the bracket to the system and further
investigates compatibilities. This was done in \cite{GL$_1$} and
thus the long standing Blaschke problem was solved.

Another example is the problem of finding the number of Abelian
relations, which is equivalent to solvability of a system of $(m+1)$
differential equations on $m$ unknown functions on the plane. The
system is a generalized complete intersection and the compatibility
conditions is given by vanishing of multi-brackets. This method was
applied in \cite{GL$_2$} and the rank problem, addressed by Lie,
Poincar\'e and Bol, was solved.

In this section we solve with our technique some other pending
problems of classical differential geometry.

\subsection{\hpss Killing vector fields on the plane}\label{51}

 \abz
If a Riemannian metric $g$ on a surface $M^2$ possesses a Killing
vector field, it has the following local form:
$ds^2=g_{11}(x)dx^2+2g_{12}(x)dxdy+g_{22}(x)dy^2$ (near the point,
where the field does not vanish) and vise versa, so that this is a
surface of revolution.

Now we address the following question: How to recognize if a metric
can be brought to such a form? This classical question was studied
by Darboux in \cite{D}. We however did not find a clearly formulated
answer in the literature. Here we give a criterion in differential
invariants using our compatibility technique.

The problem is equivalent to solvability of the equation
$L_\xi ds^2=0$, where $\xi=u\p_x+v\p_y$ is the required vector
field and $ds^2=g_{ij}(x,y)dx^idx^j$ is the metric, $x^1=x,x^2=y$.
The tensor equation is equivalent to the following 3 first order
linear PDEs on the functions $u(x,y),v(x,y)$:
 $$
2u_xg_{11}+2v_xg_{12}+u(g_{11})_x+v(g_{11})_y=0,
2u_yg_{12}+2v_yg_{22}+u(g_{22})_x+v(g_{22})_y=0,
 $$
 $$
u_yg_{11}+u_xg_{12}+v_yg_{12}+v_xg_{22}+u(g_{12})_x+v(g_{12})_y=0.
 $$
Denote them by $E_1,E_3,E_2$ respectively. We get the linear system
$\E=\E_1\subset J^1(2;2)$ of codimension 3 (we shall write
$J^k(n,m)$ instead of $J^k(\R^n,\R^m)$), so that $\dim g_1=1$, $\dim
g_2=0$, whence the isomorphism
$\pi_{2,1}:\E_2\stackrel{\sim}\to\E_1$.

The compatibility condition is equivalent to the Frobenius condition
on the corresponding distribution $L(\pi_{2,1}^{-1})$ on $\E_1$ and
is given by the condition $E'_4=[E_1,E_2,E_3]_\E=0\in\E_1^*$. This
differential operator $E_4'$ has order 2, but due to the above
isomorphism can be considered as a function on $J^1(2;2)$. However
if we consider it modulo $E_1=0$, $E_2=0$, $E_3=0$, it becomes a
function on $J^0(2;2)$ (this is not automatical and is a peculiarity
of the system) and has the form:
 $$
E_4=E_4'(\op{mod} E_1,E_2,E_3)=4|g|^2(K_xu+K_yv),
 $$
where $|g|=g_{11}g_{22}-g_{12}^2$ is the determinant of the metric
and $K$ is the Gaussian curvature. Thus compatibility condition is
equivalent to the claim that $(M^2,g)$ is a spacial form:
$K=\op{const}$. Note that this is the case, when the solutions
space has dimension 3.

Let us study solvability, then we need to add the equation $E_4=0$
to the system. This means $u=K_y w$, $v=-K_x w$ and we obtain the
following system on one function $w(x,y)$:
 $$
\begin{pmatrix}
2\a & 0 & \g_1\\ 0 & 2\b & \g_2\\ \b & \a & \g_3
\end{pmatrix}\cdot
\begin{bmatrix}
w_x \\ w_y \\ w
\end{bmatrix}=0,
 $$
where $\a=g_{11}K_y-g_{12}K_x$, $\b=g_{12}K_y-g_{22}K_x$,
$\g_1=(g_{11})_xK_y-(g_{11})_yK_x+2g_{11}K_{xy}-2g_{12}K_{xx}$,
$\g_2=(g_{22})_xK_y-(g_{22})_yK_x+2g_{12}K_{yy}-2g_{22}K_{xy}$,
$\g_3=(g_{12})_xK_y-(g_{12})_yK_x+g_{11}K_{yy}-g_{22}K_{xx}$. Note
that $\a,\b$ do not vanish simultaneously unless $K_x=K_y=0$.

Denoting by $S_1$ the determinant of the above matrix we obtain
two necessary and sufficient conditions for non-trivial
solvability ($w=0$ is always a solution): $S_1=0$ and $S_2=0$,
where:
 \begin{itemize}
 \item[--] If $\a\b\ne0$, then
$S_2=(\a(\g_1)_y-\a_y\g_1)\b^2-(\b(\g_2)_x-\b_x\g_2)\a^2$.
 \item[--] If $\a=0$, $\b\ne0$,
$S_2=\b(\g_2)_x-\b_x\g_2-\b(\g_3)_y+\b_y\g_3$
($S_1=0\Rightarrow\g_1=0$).
 \item[--] If $\a\ne0$, $\b=0$,
$S_2=\a(\g_1)_y-\a_y\g_1-\a(\g_3)_x+\a_x\g_3$
($S_1=0\Rightarrow\g_2=0$).
 \end{itemize}

Thus the criterion for existence of Killing vector field becomes
the following two non-linear differential relations $S_1=0$ and
$S_2=0$, having orders 4 and 5 in the coefficients of the metric
$g$ respectively.

For a tensor $T$ denote $d_\nabla^{\ot2}T=d_\nabla(d_\nabla T)$
the covariant derivative of the tensor $d_\nabla T$
($d_\nabla^{\ot2}$ differs from $d_\nabla^2$, which is equal to
multiplication by the curvature tensor). In particular, we obtain
the forms $d_\nabla^{\ot i}K\in C^\infty(\ot^i\,T^*M)$. Note that
the form $d_\nabla^{\ot 2}$ is symmetric, but the higher
covariant derivatives lack this property.

Let also $\op{grad}K$ be the $g$-gradient of the curvature and
$\op{sgrad}K=J\op{grad}K$ be its rotation by $\pi/2$ (fix
orientation). The preceding calculations imply:

 \begin{theorem}\po\label{Killing}
The space of local Killing vector fields can have dimension 3, 1 or
zero. A Riemannian metric $g$ possesses a local Killing vector
field iff
 $$
d_\nabla^{\ot2}K(\op{grad}K,\op{sgrad}K)=0 \text{ and }
d_\nabla^{\ot3}K(\op{sgrad}K,\op{sgrad}K,\op{sgrad}K)=0.
 $$
There are 3 independent Killing fields iff $K$ is constant.
 \end{theorem}

 \begin{rk}
The main claim of the theorem is the formula (sufficiency of which
is obvious). Other statements were known to Darboux \cite{D}. In
fact, even formulas can be attributed to him, though no precise
statement was made in \cite{D}; see \cite{K$_2$} for details.

Note also that global implications are straightforward, but
the dimension of the space of Killing vector fields can differ.
For instance, for the standard flat torus it is 2.
 \end{rk}

 \begin{proof}
Let us note that if $K\not=\op{const}$, then the system
 $$
\E'=\{E_1=E_2=E_3=E_4=0\}
 $$
has symbol dimensions: $\dim g'_0=1,\dim g'_1=0$ and so in the
compatible case the solution space is one-dimensional. Thus we
need only to prove the existence part of the theorem.

For this we express the above $S_1$ and $S_2$ via differential
invariants. Direct calculation shows:
 \begin{multline*}
S_1=-4|g|^5d_\nabla^{\ot2}K(\op{grad}K,\op{sgrad}K)\quad\text{ and }\\
S_2=A\,d_\nabla^{\ot3}K(\op{sgrad}K,\op{sgrad}K,\op{sgrad}K)
+\frac{2|g|^2}{|\op{grad}K|^4}\square(|g|^3S_1),
 \end{multline*}
where in isothermal coordinates, when $ds^2=e^\lambda(dx^2+dy^2)$,
we have:
 \begin{multline*}
A=\frac{-2K_x^2K_y^2|g|^5}{|\op{grad}K|^2}\quad\text{ and }\quad
\square=-K_xK_y^4D_x-K_x^4K_yD_y+\\
+\bigl(K_y^4K_{xx}+K_x^4K_{yy}+2K_xK_y(K_x^2+K_y^2)K_{xy}-2K_x^2K_y^2(\l_xK_x+\l_yK_y)\bigr).
 \end{multline*}
Thus $|g|^{-5}S_1$ is a     differential invariant, while
$A^{-1}S_2$ is a differential invariant relative the condition
$S_1=0$.
 \end{proof}

In the next section we'll need to enumerate differential
invariants of a Riemannian metric on a surface. It is a known fact
(see \cite{T}), that the space of scalar differential invariants
of order $k$ of a Riemannian metric on a surface is generated by
$(k-1)$ differential invariants for all $k>0$ except $k=3$, where
there is only one invariant.

The first invariants are: $I_2=K$ and $I_3=|\nabla K|^2$ (the
index refers to the order of differential invariant).

To fix a basis in invariants of order $i=2+k$ we consider the form
$d_\nabla^{\ot i}K$ and substitute $\op{grad}K$ as first $(i-j)$
arguments and $\op{sgrad}K$ as the next $j$ arguments ($0\le j\le
i$). We denote the resulting function $I_{ij}$ and enumerate the
index $j$ by letters (so we write $I_{4b}$ instead of $I_{41}$,
$I_{5d}$ instead of $I_{53}$ etc).

In these invariants the criterion of Theorem \ref{Killing} writes:
$I_{4b}=0,I_{5d}=0$.

\subsection{\hpss Higher order integrals of plane metrics}\label{52}

 \abz
Killing vector field on a surface $M$ can be represented as a linear
(in momenta) integral of the geodesic flow on $T^*M$. It is
important to know when the flow admits a polynomial integrals.
Locally geodesic flows are integrable, but the corresponding
integrals are usually analytic only on $T^*M\setminus M$. So in
general polynomial integrability requires certain conditions even
locally (here and throughout the standard regularity assumption
should be imposed).

Let $(x,y)$ be local coordinates on $M^2$ and $p_x,p_y$ be the
corresponding momenta on $T^*M$. Since every homogeneous term of the
integral is obviously an integral, we consider a function
$F_d=\sum_{i+j=d}a_{ij}(x,y)p_x^ip_y^j$ of degree $d$ on $T^*M$ ($i$
in $p_z^i$ is a power, not index).

The Hamiltonian of the geodesic flow is
$H=g^{11}p_x^2+2g^{12}p_xp_y+g^{22}p_y^2$ (matrix $g^{ij}$ is
inverse to the matrix $g_{ij}$ of the metric). Let $\{H,F_d\}$ be
the Poisson bracket of $H$ and $F_d$. It is a polynomial in momenta
of degree $d+1$.

Thus involutivity condition $\{H,F_d\}=0$ is equivalent to $(d+2)$
equations $E_1=0,\dots,E_{d+2}=0$ on $(d+1)$ unknown function
$a_{d0}(x,y),\dots,a_{0d}(x,y)$. These equations form the first
order system $\E$ of generalized complete intersection type and so
the compatibility condition can be expressed via the multi-bracket
 $$
E_{d+3}=[E_1,\dots,E_{d+2}]_\E=0.
 $$
If this condition is not satisfied we add $E_{d+3}$ to the system
and continue with investigation of solvability.

In this section we consider the case $d=2$. This is the classical
case, studied since Darboux. It is known (\cite{D,Ko}, see also
\cite{Bi}) that existence of an additional integral, quadratic in
momenta, is locally equivalent to the possibility of transforming
the metric to the Liouville form
 $$
ds^2=(f(x)+h(y))(dx^2+dy^2).
 $$

However, no effective criterion for recognizing Liouville metric was
obtained despite many attempts. The only visible success was a note
\cite{Su}. The solution to the problem was sketched there, but the
answer was not written in invariant terms (notwithstanding the
title), and the number of differential invariants characterizing
Liouville surfaces was not given (in fact, it is difficult to pursue
what the proposed set of compatibility conditions actually is and
why it is complete, so that we choose another approach below).

We describe a criterion basing on our compatibility criterion. Let
us write the metric in isothermal coordinates:
$ds^2=e^{\l(x,y)}(dx^2+dy^2)$ (the approach works with the general
form as well, but the expressions become too complicated; However
since the answer will be given in differential invariants, the
method plays no role).

The function $I=u(x,y)dx^2+2v(x,y)dxdy+w(x,y)dy^2$ is a quadratic
integral of the geodesic flow iff the following system $\E$
(coefficients of $\{H,I\}=0$) is satisfied:
 \begin{align*}
u_x&+\l_x u+\l_y v=0, & u_y+2v_x&+\l_x v+\l_y w=0, \\
2v_y+w_x&+\l_x u+\l_y v=0, & w_y&+\l_x v+\l_y w=0.
 \end{align*}
Denoting the equations by $E_1,E_2,E_3,E_4$ we obtain their
compatibility condition
 $$
E_5=\frac12[E_1,E_2,E_3,E_4]_\E=0.
 $$
From the general theory it might be expected that $E_5$ has order 2
(3 in non-reduced form), but in fact it is of the first order and
has the following form (after cancelation by $2\sqrt{\det(g)}$):
 \begin{multline*}
E_5=5K_xv_x-5K_yv_y\\
-(K_{xy}+2\l_yK_x+2\l_xK_y)(u-w)+
(K_{xx}-K_{yy}+4\l_xK_x-4\l_yK_y)v.
 \end{multline*}

Thus the system $\E$ is integrable iff $K=\op{const}$. In this
case dimension of the solutions space is $\sum\dim g_k=6$ ($\dim
g_i=\max\{3-i,0\}$) and the space of quadratic integrals is the
symmetric square of the 3-dimensional space of linear integrals (a
basis of the former is the pair-wise product of a basis of the
latter).

Suppose that $K\ne\op{const}$, so that at least one of the
functions $K_x,K_y$ is not zero. We add the equation $E_5=0$ and
get a system $\E'\subset J^1(2,3)$ of formal codimension 5.

Its symbols $g_i\subset S^iT^*\ot\R^3$ have $\dim g_0'=3$, $\dim
g_1'=1$, $\dim g_2'=0$ and thus the only non-zero second
$\d$-cohomology are $H^{0,2}(\E')\simeq\R^1$, $H^{1,2}(\E'')\simeq
\R^1$. There are two obstructions to compatibility -- Weyl tensors
$W_1'$ and $W_2'$. The tensor $W_1'$ is proportional to
  $$
E_6'=K_xE_{5x}+K_yE_{5y}-\frac52K_x^2(E_{2x}-E_{1y})
+\frac52K_y^2(E_{3y}-E_{4x})
 $$
Dividing this by $5K_y$ and simplifying modulo
$E_1,E_2,E_3,E_4,E_5$ we obtain the following expression:
 $$
E_6= -35\,|g|\,I_{4b}\cdot v_x+Q_1\cdot(u-w)+Q_2\cdot v,
 $$
where $Q_1$, $Q_2$ are certain differential expression of 5th order
in the coefficients of the metric. The coefficients of $E_6$, as
well as others $E_i$, are not invariant, but the condition of their
vanishing is invariant, and so can be expressed in terms of
differential invariants. Indeed,
 \begin{gather*}
Q_1=J_5\cdot I_3^{-3}\sqrt{|g|}(-K_x)(I_3\sqrt{|g|}-K_x^2),\\
Q_2=J_5\cdot I_3^{-3}\sqrt{|g|} (I_3\sqrt{|g|}-K_x^2)^{1/2}
(I_3\sqrt{|g|}-2K_x^2),
 \end{gather*}
where
 $$
J_5=5I_3(I_{5a}-I_{5c})+(I_{4a}-I_{4c})(I_{4c}-6I_{4a})-25I_2I_3^3
 $$
is a differential invariant. Thus the equations $Q_1=Q_2=0$ are
equivalent to one condition $J_5=0$.

It is possible to show that this condition together with $I_{4b}=0$
implies $I_{5d}=0$, which gives another proof of Darboux theorem
\cite{D} (proof in \cite{Ko}) that a Riemannian surface with 4
quadratic integrals is a surface of revolution.

The second obstruction to existence of 4 integrals -- tensor $W_2'$
-- can be calculated similarly. Its vanishing is given by a scalar
differential invariant of order 6 in metric, but it can be
simplified modulo the conditions $I_{4b}=I_{5d}=J_5=0$ to the
following expression:
 $$
J_4=3(I_{4a}-I_{4c})(I_{4a}+4I_{4c})I_{4c}-15I_2I_3^3(I_{4a}+4I_{4c})+25I_3^5.
 $$
Thus we obtain the following statement:

 \begin{theorem}\po\label{thdeg+}
The condition of exactly 4 quadratic integrals can be expressed as 3
differential conditions on the metric: $I_{4b}=0$, $J_5=0$, $J_4=0$.
 \end{theorem}

If the compatibility condition $E_6(\op{mod}
E_1,E_2,E_3,E_4,E_5)=0$ is satisfied, then the system $\E'$ is
integrable. Otherwise we add this new equation and get (again in
generic case, when the corresponding matrix of coefficients of
derivatives is non-degenerate) the system $\E''$ with symbol
$g_1''=0$, i.e. it is of Frobenius type.

Its Spencer cohomology group $H^{0,2}(\E'')\simeq\R^3$, so the
obstruction to integrability -- curvature tensor -- $W_1''$ has 3
components, represented by 3 linear equations relations on
$J^0(2,3)$:
 $$
E_{7j}=A_{1j}(u-w)+A_{2j}v=0,\qquad j=1,2,3.
 $$
The expressions $A_{ij}$ are not invariant, but their vanishing is
invariant and can be expressed via four differential invariants of
order 6 of the metric:
 $$
J_{6k}=I_{6k}-P(I_2,I_3,I_{4i},I_{5j}),\qquad k=a,b,c,d,
 $$
where $P$ is a quadratic function in $I_{5j}$ with rational
coefficients in other variables (note that $I_{6e}$ does not enter
the formulae). All the expressions are rather long and shall be
provided elsewhere. Let us indicate only equation $J_{6a}=0$:

 \begin{multline*}
I_{6a}=\frac1{175I_3^2I_{4b}}
 \bigl( 700I_3^5I_{4b} - 825I_2I_3^4I_{5b} +
50I_2I_3^3I_{4b} (31I_{4a} - 18I_{4c})\\
 + 6I_{4b} (I_{4a} - I_{4c}) (6I_{4a}^2 + 49I_{4b}^2
- 37I_{4a}I_{4c} + 6I_{4c}^2 ) - 25I_3^2I_{5b} (-8I_{5a} +I_{5c})\\
 - 5I_3 (48I_{4a}^2I_{5b} - 27I_{5b}I_{4c}^2 +
2I_{4b}I_{4c} (-11I_{5a} + 46I_{5c} )\\
 + I_{4a} (-43I_{5a}I_{4b} - 21I_{5b}I_{4c} +
8I_{4b}I_{5c}) + 7I_{4b}^2 (4I_{5b} - 11I_{5d}) )
 \bigr).
 \end{multline*}

 \begin{theorem}\po\label{thdeg3}
The condition of exactly 3 quadratic integrals can be expressed as
4 differential conditions on the metric:
$J_{6a}=J_{6b}=J_{6c}=J_{6d}=0$.
 \end{theorem}

Finally if $E_{7j}$ are non-zero, we add these equations to the
system. Compatibility condition of the new system $\E'''$ are
$(E_{7j})_x=0$, $(E_{7j})_y=0$, when expressed as linear functions
on $J^0(2,3)$ via the system $\E''$:
 $$
E_{8l}=B_{1l}u+B_{2l}v+B_{3l}w=0.
 $$
Consider the matrix of coefficients of equations $E_{7j},E_{8l}$:
$U=U(A,B)$. It always satisfies the condition $\op{rank}(U)<3$,
because $H$ is an integral of the geodesic flow. Also
$\op{rank}(U)>0$ if conditions of Theorem \ref{thdeg3} are not
fulfilled.

Thus we have only two possibilities: If $\op{rank}U=2$, the flow
does not possesses local quadratic in momenta integrals. Otherwise
$E_{8l}(\op{mod} E_{7j})=0$. This means $\op{rank}U=1$ and
expressing this condition in differential invariants we get one
condition of order 6 and four conditions of order 7 in coefficients
of the metric (these long expressions will be omitted; note though
that the above 4 scalar differential invariants of order 7 involve 6
basic invariants $I_{7k}$, but the last of them $I_{7f}$ does not
enter the formulae): $\tilde J_6=0$, $J_{7i}=0$.

These conditions give solvability of the system $\E$, which yields
us a 2-dimensional linear space of solutions generated by $H$ and
$I$ -- an independent integral of degree 2. We set $\square=(\tilde
J_6,J_{7a},J_{7b},J_{7c},J_{7d})$.

Denote by $\mathcal{S}$ the singular locus of functions $\l(x,y)$,
that are non-generic w.r.t.\ at least one one of the above steps
(it consists of functions of one variable -- metrics with non-zero
Killing fields, and certain finite-dimensional families), which
corresponds to zero denominators in $J_{7i}$.

Let also $\mathcal{L}_+$ be the set of $\l$ corresponding to the
metrics with more than one additional quadratic integral (constant
curvature or the conditions of Theorems \ref{thdeg+} and
\ref{thdeg3}). We have proved:

 \begin{theorem}\po
There exists a polynomial vector-valued differential operator of
order 7 \/ $\square:C^\infty_\text{loc}(\R^2)\to
C^\infty_\text{loc}(\R^2;\R^5)$ and a residual subset
$\mathcal{S}\subset C^\infty_{loc}(\R^2)$ such that for $\l\in
C^\infty_\text{loc}(\R^2)\setminus\mathcal{S}$ the metric
$e^{\l(x,y)}ds^2_\text{Eucl}$ is Liouville (has quadratically
integrable geodesic flow) iff\/ $\l\in\mathcal{L}_+$ or
$\square(\l)=0$. Moreover, Liouville metrics in $\mathcal{S}$ are
residual among all Liouville metrics.
 \end{theorem}

The singular locus $\mathcal{S}$ in the space of germs of Riemannian
metrics (which we identified with $C^\infty_\text{loc}(\R^2)$ only
for convenience) is given by the condition $I_{4b}=0$. The
expressions of the above invariants together with a more detailed
argumentation have appeared now in \cite{K$_2$}.

\subsection{\hpss Gaussian curvature of minimal surfaces}\label{53}

 \abz
Consider a minimal surface $M^2\subset\R^3$. The Gauss map defined
on it depends on the curvature function and this function is
unrestricted (i.e. can be arbitrary in a certain open domain; of
course, it is non-positive, but there are no
equality-restrictions) if the surface is considered abstractly
(non-parametrized). But it's quite known that the Gauss map is not
arbitrary, which is manifested by the fact, that the Gaussian
curvature on the immersed (parametrized) surface is not arbitrary.
We will describe precisely, which functions $K$ on
$M^2\subset\R^3$ are realized locally.

So let $M^2$ be given as the graph $z=u(x,y)$. Then
$\nabla_1(u)=\frac{u_{xx}u_{yy}-u_{xy}^2}{(1+u_x^2+u_y^2)^2}$ is
the Gaussian curvature operator and
$\nabla_2(u)=\frac{(1+u_x^2)u_{yy}-2u_xu_yu_{xy}+(1+u_y^2)u_{xx}}
{(1+u_x^2+u_y^2)^{3/2}}$ is the operator of mean curvature. Let
$\mathcal{H}^\infty=\{u\in C^\infty_\text{loc}(\R^2)\,|\,
\nabla_2(u)=0\}$ be the sheaf of minimal surfaces. We define
$\nabla_1:\mathcal{H}^\infty\to C^\infty_\text{loc}(\R^2)$ and
denote the image by $\mathcal{K}^\infty$. Now we want to resolve
this term:

 \begin{theorem}\po
There exists an algebraic differential operator
$\square:C^\infty_\text{\rm loc}(\R^2)\to C^\infty_\text{\rm
loc}(\R^2;\R^2)$ of order 4 and a finite-dimensional stratified
submanifold $\mathcal{S}\subset C^\infty_\text{\rm loc}(\R^2)$
with $\mathcal{K}^\infty\setminus\mathcal{S}= \op{Ker}(\square)$.

More precisely, there exist 4 polynomials $F_3,F_6,F_7,F_8$ on the
plane with coefficients depending differentially on $K$ such that
the function is realized as the curvature of a minimal surface iff
they have a common root.
 \end{theorem}

The form of the operator $\square=(\square_1,\square_2)$ will be
clear from the proof, though we suppress the formulas because of
their size. Since the operator has singularities, we obtain cases
and $K\in\mathcal{K}^\infty\setminus\mathcal{S}$ iff
$\square_1(K)=\square_2(K)=0$, while description of
$\mathcal{K}^\infty\cap\mathcal{S}$ is given by some other
operators, which we omit.

 \begin{proof}
Let us find the solvability criterion for the system
$\nabla_1(u)=K(x,y)$, $\nabla_2(u)=0$. Denote
 $$
F_1=u_{xx}u_{yy}-u_{xy}^2-K\cdot(1+u_x^2+u_y^2)^2,
F_2=(1+u_y^2)u_{xx}-2u_xu_yu_{xy}+(1+u_x)^2u_{yy}.
 $$
The Mayer bracket of these operators is
 $$
F_3=[F_1,F_2]=a_{11}u_x^2+2a_{12}u_xu_y+a_{22}u_y^2+b,
 $$
where $a_{11}=K^2(\ln|K|)_{yy}-4K^3$, $a_{12}=-K^2(\ln|K|)_{xy}$,
$a_{22}=K^2(\ln|K|)_{xx}-4K^3$,
$b=K^2((\ln|K|)_{xx}+(\ln|K|)_{yy})-4K^3$. Notice that $F_3$ has
order 1, while generically the bracket of 2 second order operators
after reduction is also of 2nd order. Therefore the system is
compatible if $F_3=0$, which is equivalent to $K=0$, i.e. the
surface is a plane.

Denote $F_4=\D_x(F_3)$, $F_5=\D_y(F_3)$. The system
$F_2=F_4=F_5=0$ has the form:
 $$
 \begin{pmatrix}
(1+u_y^2) & -2u_xu_y & (1+u_x^2) \\
2(a_{11}u_x+a_{12}u_y) & 2(a_{12}u_x+a_{22}u_y) & 0 \\
0 & 2(a_{11}u_x+a_{12}u_y) & 2(a_{12}u_x+a_{22}u_y)
 \end{pmatrix}
 \begin{bmatrix}
u_{xx} \\ u_{xy} \\ u_{yy}
 \end{bmatrix}
 =
 \begin{bmatrix}
0 \\ b_1 \\ b_2
 \end{bmatrix},
 $$
where
$b_1=-b_x-(a_{11})_xu_x^2-2(a_{12})_xu_xu_y-(a_{22})_xu_y^2$,
$b_2=-b_y-(a_{11})_yu_x^2-2(a_{12})_yu_xu_y-(a_{22})_yu_y^2$.

Resolving this for the second derivatives and substituting to
$F_1$ we get after multiplication by the square of the determinant
of the above matrix a first order polynomial differential operator
$F_6'$. It has degree 8 by variables $u_x,u_y$, but reduction due
to the system is of degree 6:
 \begin{multline*}
F_6=-F_6'+4K(1+u_x^2+u_y^2)^2(2bF_3-F_3^2)\\
\!\!\!= b_1^2+b_2^2+(b_2u_x-b_1u_y)^2+4K(1+u_x^2+u_y^2)^2
\bigl(b^2+(a_{11}u_x+a_{12}u_y)^2+(a_{12}u_x+a_{22}u_y)^2\bigr).
 \end{multline*}

Denote $E_7=[E_3,E_6]$ the Mayer bracket. It is the third 1st
order PDE, which can be equivalently written as
 $$
\frac{\p(F_3,F_6)}{\p(x,u_x)}+\frac{\p(F_3,F_6)}{\p(y,u_y)}=0.
 $$

Let us prolong the equation $F_6=0$: $\D_x(F_6)=0$, $\D_y(F_6)=0$.
Resolving these equations coupled with $F_4=F_5=0$ with respect to
second derivatives ($F_7=0$ guarantees compatibility) and
substituting to $F_2=0$ we get:
 $$
F_8=(1+u_y^2)\dfrac{\p(F_3,F_6)}{\p(x,u_y)}+
u_xu_y\left(\dfrac{\p(F_3,F_6)}{\p(x,u_x)}
-\dfrac{\p(F_3,F_6)}{\p(y,u_y)}\right)
-(1+u_x)^2\dfrac{\p(F_3,F_6)}{\p(y,u_x)}.
 $$
Thus we obtain four 1st order PDEs $F_3=0$, $F_6=0$, $F_7=0$,
$F_8=0$, which are polynomials by $u_x,u_y$ of degrees 2, 6, 7, 9
(they do not depend on $u$) and are differential operators by $K$
of orders 2, 3, 4, 4 respectively.

Solvability of the system $F_1=F_2=0$ is equivalent to the claim
that polynomial by $u_x,u_y$ system $F_3=F_6=F_7=F_8=0$ has a
solution. The latter is equivalent to 2 conditions $\square_1=0$,
$\square_2=0$, algebraic by the corresponding coefficients, which
are differential operators of $K$. Therefore we set
$\square=(\square_1,\square_2)$.

The set $\mathcal{S}$ is formed by functions $K$ for which some of
$F_3,F_6,F_7,F_8$ become dependent or singular. This set is given
by a collection of overdetermined systems of PDEs of finite type
and hence is stratified finite-dimensional.

Let us demonstrate a geometric idea behind this proof. $F_1$ and
$F_2$ are Monge-Amp\`ere operators on the plane, which means
(\cite{L$_2$}) that they are given by 2-forms $\Omega_1,\Omega_2$ on
$J^1(\R^2)$: Equations $F_i(u)=0$ can be rewritten as
$\Omega_i|_{j_1(u)}=0$, where $j_1(u)\subset J^1(\R^2)$ is the
jet-section determined by $u$ and
 \begin{gather*}
\Omega_1=du_x\we du_y-K(1+u_x^2+u_y^2)^2dx\we dy,\\
\Omega_2=(1+u_y^2)du_x\we dy+u_xu_y(dx\we du_x+du_y\we dy)+
(1+u_x^2)dx\we du_y.
 \end{gather*}
Since $\Omega_i$ do not depend on $u$, the construction descends
onto $T^*\R^2=J^1(\R^2)/\R^1$, where we have in addition the
canonical symplectic 2-form $\Omega_0$.

We search for a Lagrangian surface, which is isotropic w.r.t.
$\Omega_1,\Omega_2$ and lies in the hypersurface
$\Sigma^3=\{F_3=0\}\subset T^*\R^2$. Let $\xi_i$ be the kernels of
the 2-forms $\Omega_i$ restricted to $\Sigma^3$. They fail to be
in general position only along a surface $\Sigma^2=\{F_3=F_6=0\}$,
which is equivalently determined by the equation
 $$
\Sigma^2=\{x\in\Sigma^3\,|\,\op{rank}(\xi_0(x),\xi_1(x),\xi_2(x))<3\}.
 $$
This surface is ramified over $\R^2$ with the fiber consisting of
no more than 12 points (precisely this number if counted over $\C$
and with multiplicities).

Local solvability means that over a neighborhood in $\R^2$ at
least one component $\Sigma^2_\a$ of $\Sigma^2$ is isotropic for
all $\Omega_i$. By characteristic property of this surface we need
to require only two restrictions $\Omega_i|_{\Sigma^2_\a}$ to
vanish. For instance, we take the conditions $\xi_0,\xi_2\subset
T\Sigma^2_\a$. They are given by the equations $\Omega_0\we
dF_3\we dF_6=0$ and $\Omega_2\we dF_3\we dF_6=0$, which correspond
to the operators $F_7$ and $F_8$.
 \end{proof}

 \begin{examp}\po
Consider the family $K=ke^{\a x+\b y}$. Then it realizes the
curvature of a minimal surface iff $k=0$. In fact, the quadric
$F_3-b$ is definite and $a_{11}b\ge0$. So $F_3k\le0$ and the
inequality is strict if $k\ne0$.
 \end{examp}

Notice that $F_3$ is non-singular if $K\ne0$. The set
$\mathcal{K}^\infty$ is a subset of the solutions to $F_3=F_6=0$.
Let us describe non-holonomic solutions to this system (ignoring
the compatibility condition $F_7=0$).

Let's parameterize the quadric $F_3=0$: Consider a line
$u_x=t,u_y=\l t$. It meets $F_3=0$ at the points with
$t^2=-b/\rho(\l)$, $\rho(\l)=a_{11}+2a_{12}\l+a_{22}\l^2$.
Substituting this into $\rho(\l)^3F_6=0$ we get a polynomial of
degree $6$:
 \begin{multline*}
q_6(\l)=\rho(\l)^3b^4\Bigl[\Bigl(\frac{\rho(\l)}b\Bigr)_x^2+
\Bigl(\frac{\rho(\l)}b\Bigr)_y^2\Bigr]-
b^5\Bigl[\l\Bigl(\frac{\rho(\l)}b\Bigr)_x-
\Bigl(\frac{\rho(\l)}b\Bigr)_y\Bigr]^2\\
+4K\bigl[\rho(\l)-b(1+\l^2)\bigr]^2
\bigl[b^2\rho(\l)-b\bigl((a_{11}+a_{12}\l)^2
+(a_{12}+a_{22}\l)^2\bigr)\bigr].
 \end{multline*}
Thus $\{F_3=0\}\cap\{F_6=0\}$ corresponds (2-to-1) to the roots of
$q_6(\l)=0$.

 \begin{examp}\po
Consider the family $K=\vp(x)$. For an open set of such functions
the equation $q_6(\l)=0$ has 6 positive roots and so the system
$F_3=F_6=0$ has 12 non-holonomic solutions $u_x,u_y$. None of them
satisfies the other two equations $F_7=F_8=0$, save for the case
$u_x=\op{const}_1$, $u_y=\op{const}_2$. Actually, the solutions
depend on $x$ only and so $u_{xy}=u_{yy}=0$, whence $u_{xx}=0$ and
$K=0$.
 \end{examp}

Let $\mathcal{S}_0\subset\mathcal{S}$ be the set of functions $K$,
such that $\{F_3=0\}\equiv\{F_6=0\}$. It is given by the condition
$q_6(\l)\equiv0$, which is a system of 7 third order PDEs on $K$.
As a by-product of the proof we get the following statement:

 \begin{prop}\po\label{thrm9}
For $K\in\mathcal{K}^\infty\setminus\mathcal{S}_0$ the set of
minimal surfaces through the origin with this curvature
$\nabla_1^{-1}(K)\cap\mathcal{H}^\infty\cap\{u(0,0)=0\}$ has
cardinality at most 12, though generically this number is 2
(respectively 6 and 1 if we the function $u$ is considered up to
the sign).
 \end{prop}

 \begin{proof}
In fact, all the equations are symmetric to the change
$u\mapsto-u$. The sub-system $F_3=F_6=0$ with respect to $u_x,u_y$
has as maximum $\op{ord}(F_3)\cdot\op{ord}(F_6)=12$ algebraic
solutions. Some of them may not satisfy the other constraints
$F_7=F_8=0$ and generically (in
$\mathcal{K}^\infty\setminus\mathcal{S}_0$) only 2 do satisfy.
 \end{proof}

The above statement does not hold for $K=0\in\mathcal{S}_0$. But
we suggest there are no other examples. To see the reason let us
consider the overdetermined system of 7 polynomial differential
equations on $K$ of order 3:
$q_6(0)=0,q_6'(0)=0,\dots,q_6^{(6)}(0)=0$. One can expect that it
has no other solutions than $K=0$, but this is not so.

For instance, $q_6(\l)\equiv0$ follows from $b=0$ and
$\left|\begin{matrix} a_{11} & a_{12} \\ a_{12} & a_{22}
\end{matrix}\right|=0$. This latter is equivalent to the
system $\op{det}\op{Hess}(\ln|K|)=0$, $\Delta\ln|K|=4K$, which
though incompatible is solvable and has the solution $K=\vp(x\pm
y)$, where $\vp'=2\vp\sqrt{\vp+c}$.

But for this choice of $K$ the other two equations $F_7=F_8=0$ are
not satisfied unless $\vp=0$. In fact, then $F_3=0$ and $F_6=0$
are both equivalent to $u_x+u_y=0$, which coupled with $F_1=F_2=0$
gives $K=0$.

The computer programs do not give us other solutions to the above
system of 7 third order PDEs on $K$ and so {\it we conjecture
that}
 $$
\mathcal{K}^\infty\cap\mathcal{S}_0=0.
 $$

 \begin{rk}\po
This our conjecture that except the plane case a minimal surface is
restored (up to $\pm$) from its Gauss image in maximal 6 different
ways is similar to the known Gronwall conjecture about webs on the
plane. It says that a linearizable 3-web has maximally 11
linearizations (and generically only 1), cf. \cite{GL$_1$}. Both
problems are basically algebraic and have equal complexities.
 \end{rk}

Similarly one investigates the problem, when two functions $K$ and
$H$ on a paramet\-rized surface $M^2$ can be realized as Gaussian
and mean curvatures, with the surface realized as a graph
(projection $\R^3\to\R^2$ yields the parametrization).

This is an analog of the classical Bonnet problem of realizing two
quadrics as the first and the second quadratic forms on a surface.
Bonnet theorem states that compatibility and solvability of this
problem is the system of one Gauss and two Kodazzi equations.

For realization of the curvatures $K,H$ compatibility is equivalent
to the condition $K=H=0$, i.e. the surface is plane (notice that the
solutions space is 3-dimensional, not 4-dimensional as one can
expect after \S\ref{S43+e}, but this is due to non-genericity of the
condition). Solvability leads to an operator of order 4, similar to
the above $\square$. Thus we get solution to generalized Bonnet
problem.

\subsection{\hpss Quantum integration}\label{63}
 \abz
Consider the algebra of scalar linear differential operators
$\mathcal{A}=\op{Diff}(\1,\1)$ on the manifold $M$ filtered by the
$C^\infty(M)$-modules $\mathcal{A}_k$ of order $\le k$ differential
operators. Let
$\mathcal{P}=\oplus_{k\ge0}\mathcal{P}_k=\op{gr}(\mathcal{A})$ be
the corresponding graded module. Here
$\mathcal{P}_k=\mathcal{A}_k/\mathcal{A}_{k-1}$ consists of degree
$k$ homogeneous in momenta polynomials on $T^*M$. Thus
$\mathcal{P}=S\mathcal{D}=\oplus S^k\mathcal{D}$, where
$\mathcal{D}$ is the $C^\infty(M)$-module of vector fields on $M$.
The canonical Poisson structure on $\mathcal{P}$ is given by the
formula
 $$
\{\nabla_1\op{mod}\mathcal{A}_k,\nabla_2\op{mod}\mathcal{A}_l\}
=[\nabla_1,\nabla_2]\op{mod}\mathcal{A}_{k+l-1},
 $$
where the bracket in the r.h.s. is the usual commutator (or Jacobi
bracket). In other words, the mapping $\z=\op{smbl}:\mathcal{A}\to
\mathcal{P}$ is a homomorphism of Lie algebras.

By quantization one understands an inverse map $q_*:\mathcal{P}
\to\mathcal{A}$, i.e. a collection of morphisms
$q_k:\mathcal{P}_k\to\mathcal{A}_k$ splitting the sequence
 $$
0\to\op{Diff}_{k-1}(\1,\1)\hookrightarrow\op{Diff}_k(\1,\1)
\stackrel{\dashleftarrow}\longrightarrow S^k\mathcal{D}\to0.
 $$

This map allows to introduce a new non-commutative associative
product on $\mathcal{P}$:
 $$
a\star b=q_*^{-1}(q_*(a)\circ q_*(b)).
 $$
These kinds of products are important in the deformation
quantization. Moyal \cite{Mo} and other star-products \cite{L$_3$}
are obtained by specifying the morphism $q_*$.

Denoting by $p_k$ the homogeneous $\mathcal{P}_k$-component of a
polynomial $p\in\mathcal{P}$ we observe the relation
 $$
\{a,b\}=(a\star b-b\star
a)_{k+l-1}=\op{smbl}_{k+l-1}([q_k(a),q_l(b)]),\quad
a\in\mathcal{P}_k,b\in\mathcal{P}_l
 $$
between the commutator, Poisson bracket and the star-product.

Consider a mechanical system with a Hamiltonian $h\in\mathcal{P}$.
Its quantization is given by $H=q_*(h)$. If the choice of $q_*$ is
subject to certain connections (\cite{L$_3$}), then Riemannian
metric produces the Laplace operator, choice of potential --
Shr\"odinger operator etc.

Suppose the classical system is integrable in Liouville sense,
i.e. there exist functions $f_1=h,f_2,\dots,f_n\in\mathcal{P}$
($n=\dim M$) functionally independent a.e. which Poisson-commute
$\{f_i,f_j\}=0$. We wish to quantize this picture.

 \begin{dfn}\po
Differential operator $H$ is called quantum completely integrable
if there exist commuting differential operators
$F_1=H,F_2,\dots,F_n\in\mathcal{A}$, which are independent a.e.
 \end{dfn}
Clearly then the system $h=\op{smbl}(H)$ is Liouville-integrable
with integrals $f_i=\op{smbl}(F_i)$. The quantization poses the
inverse problem: To find quantum integrable system $(H,F_i)$ by
the given classical $(h,f_i)$.

This problem was solved for many classically integrable
Hamiltonian systems (\cite{Pe}) using different approaches:
analytical, Dunkl's differential-difference operator \cite{Du},
Moyal quantization, via geodesic equivalence \cite{MT} and others.

We discuss one of them, which is closely related to our
integration method. It was proposed in \cite{He} and is based on
universal enveloping algebras.

Consider the rigid body equations, which is the Hamiltonian system
on $T^*SO(3)$ with Hamiltonian $h=\frac12\sum
\t_i^{-1}p_i^2+\sum\g_ix_i$, where $x_i$ are the base coordinates
and $p_i$ are the corresponding momenta.

Let $X_i=q_*(x_i)$ at $R\in\op{SO}(3)$ be equal to $(Re_i,e)$,
where $e_i$ is an orthonormal basis and $e$ some unit vector in
$\R^3$, and $P_i=q_*(p_i)$ be the left-invariant fields
$\exp(E_i)$ generated by the basis $E_i\in\op{so}(3)$ given by the
relations $E_i(e_i)=0$, $E_i(e_{i\pm1})=\pm e_{i\mp1}$,
$i\in\Z_3$. Then the subalgebra of $\mathcal{A}(\op{SO}(3))$
generated by $P_i,X_i$ is isomorphic to the universal enveloping
algebra of $\op{so}(3)\ltimes\R^3$: $[P_i,P_j]=\epsilon_{ijk}P_k$,
$[P_i,X_j]=\epsilon_{ijk}X_k$, $[X_i,X_j]=0$.

The quantized Hamiltonian has the form $H=q_*(h)=\frac12\sum
\t_i^{-1}P_i^2+\sum\g_iX_i$ and we denote it also by $F_1$. We
have two Casimir functions in $U(\op{so}(3)\ltimes\R^3)$:
$F_2=\sum X_i^2$, $F_3=\sum P_iX_i$. To achieve complete
quantization of all known integrable Euler equations we must
quantize the forth integral. It is as follows:

{\bf Euler case:} $\g_i=0$. Then $F_4=\sum P_i^2$.

{\bf Lagrange case:} $\t_1=\t_2$, $\g_1=\g_2=0$. Then $F_4=P_3$.

{\bf Kovalevskaya case:} $\t_1=\t_2=\t_3/2$, $\g_3=0$. Then
$F_4=K\bar K+\bar K K-g_4$, where
$K=\t_1(P_1+iP_2)^2-2(\g_1+i\g_2)(X_1+iX_2)$ and
$g_4=8\t_1^2(P_1^2+P_2^2)$. In all these cases $[H,F_4]=0$.

{\bf Goryachev-Chaplygin case (conditional integrability):}
$\t_1=\t_2=\t_3/4$, $\g_3=0$. Then
$F_4=\t_1(P_1^2+P_2^2)P_3-X_3(\g_1P_1+\g_2P_2)-g_3$, where
$g_3=\frac12(\g_2X_1-\g_1X_2)+\frac14\t_1P_3$. In this case
$[H,F_4]=\t_1(\g_2P_1-\g_1P_2)F_3$.

Note that according to \cite{L$_3$} any quantization $q_*$ is
determined by two linear connections (electromagnetic field and
gravity) and a collection of tensors
$g_k:S^k\mathcal{D}\to\mathcal{A}_{k-2}$. The first two of the above
integrable cases are obtained from the classical scheme via the
operator $q_*$ and trivial tensors $g_2$ and $g_1$. In the
Kovalevskaya case one needs a second order correction $g_4$ and in the
Goryachev-Chaplygin case a first order $g_3$ (these corrections were
found previously in \cite{He}, but the explanation of actual orders
meaning was lacking).

 \begin{theorem}\po
In each of the classical integrable cases the obtained quantum
integrals allow to integrate the Shr\"odinger operator
$L[u]=u_t-H(u)-\l_1u$ (where $\l_1$ is the spectral parameter)
classically: The system
 $$
L[u]=0,\ F_2(u)=\l_2u,\ F_3(u)=\l_3u,\ F_4(u)=\l_4u
 $$
is of finite type and compatible ($\l_3=0$ in the
Goryachev-Chaplygin case).
 \end{theorem}

 \begin{proof}
In the first three cases, where we have a commutative collection
of integrals $F_1,F_2,F_3,F_4$, the statement is rather known.
However, the Goryachev-Chaplygin case seems to be quantization by
analogy and its meaning is given by the above statement, which
follows from our compatibility criterion.
 \end{proof}


 \begin{rk}\po
The above result is constructive, not mere an existence statement.
In fact, the considered system is of Frobenius type and one
reduces its integration to a certain system of ODEs, which may be
integrated via symmetry approach.

So far we obtained only local solutions, but if we are given a
global Shr\"odinger equation on a manifold $M$, then we get a
topological restriction for a compatibility/solvability. Namely
the monodromy operator determines which spectral parameters $\l$
are admissible. For them and only for them we get a closed leaf of
the foliation corresponding to the above Frobenius system.
 \end{rk}

Note that in the definition of quantum integrability we used
$\op{smbl}$-map, not $q_*^{-1}$. In fact, $q_*^{-1}(F_i)$ need not
to commute, only their top-components. Denote by
$\rho_k:\mathcal{P}\to\mathcal{P}_k$ the natural projection. Then
commutation of $q_k(a)$, $q_l(b)$ implies that
$\{a,b\}=\rho_{k+l-1}q_*^{-1}([q_*(a),q_*(b)])=0$. So the
quantization philosophy suggests to make the construction so that
the functions commute w.r.t. the deformed bracket
 $$
\{a,b\}_q=q_*^{-1}([q_*(a),q_*(b)]).
 $$

This is equivalent to the quantum integrability problem: given a
system $f_1,\dots,f_n$ involutive w.r.t.\ the Poisson bracket
describe the quantizations $q_*$ such that the set is still
involutive w.r.t.\ the new bracket $\{,\}_q$.

However our compatibility result suggest that it is equally
important to search for deformations giving sub-algebras, i.e.
$\{f_i,f_j\}_q=\sum c_{ij}^kf_k$. In the classical case they
correspond to invariant submanifolds of the Hamiltonian system.
Thus our main result interprets as quantization of conditional
integrability in the classical mechanics.

\small

\vspace{-10pt} \hspace{-20pt} {\hbox to 12cm{ \hrulefill }}
\vspace{-1pt}

{\footnotesize \hspace{-10pt} Institute of Mathematics and
Statistics, University of Troms\o, Troms\o\ 90-37, Norway.

\hspace{-10pt} E-mails: \quad kruglikov\verb"@"math.uit.no, \quad
lychagin\verb"@"math.uit.no.} \vspace{-1pt}

\end{document}